%  processed by /home/shlhetal/bin/jakob/citeatalya.pl  0.3, July 2005 on 2006-08-12
\ifx\shlhetal\undefinedcontrolsequence\let\shlhetal\relax\fi
% To: shlhetal@math.huji.ac.il
% Subject: 861
% Date: Sat, 12 Aug 2006 16:20:44 +0300
% From: atalya <atalyaw@gmail.com>
% Mime-Version: 1.0
% X-sliced-and-diced-by: 'savemail' 2.0, Oct 2005

% this file was modified by a referee and sent back on 2006-04-20 (JK)
% 861, processed by citealice (2004-09-11) on Thu Mar 31 21:11:25 IST 2005

\ifx\shlhetal\undefinedcontrolsequence\let\shlhetal\relax\fi
% To: Saharon Shelah <shelah@math.huji.ac.il>, Andrzej Roslanowski <roslanow@member.ams.org>, %"Saharon Shelah's Office" <shlhetal@math.huji.ac.il>, <shani@math.huji.ac.il>
% Subject: Sh861
% Date: Thu, 31 Mar 2005 13:51:55 -0500 (EST)
% From: Alice Leonhardt <leonhard@math.rutgers.edu>
% Mime-Version: 1.0
% Content-Description: revisions
% X-sliced-and-diced-by: 'savemail' 1.3b, Feb 2003

\def\smallbox#1{\leavevmode\thinspace\hbox{\vrule\vtop{\vbox
   {\hrule\kern1pt\hbox{\vphantom{\tt/}\thinspace{\tt#1}\thinspace}}
   \kern1pt\hrule}\vrule}\thinspace}

\def\fmtname{AmS-TeX}

\def\fmtversion{2.2}
\catcode`\@=11
\ifx\amstexloaded@\relax\catcode`\@=\active
  \endinput\else\let\amstexloaded@\relax\fi
\newlinechar=`\^^J
\def\W@{\immediate\write\sixt@@n}
\def\CR@{\W@{^^J\fmtname - Version \fmtversion^^J}}
\CR@ \everyjob{\CR@}
\message{Loading definitions for}
\message{misc utility macros,}
\toksdef\toks@@=2
\long\def\rightappend@#1\to#2{\toks@{\\{#1}}\toks@@
 =\expandafter{#2}\xdef#2{\the\toks@@\the\toks@}\toks@{}\toks@@{}}
\def\alloclist@{}
\newif\ifalloc@
\def\showallocations{{\def\\{\immediate\write\m@ne}\alloclist@}\alloc@true}
\def\alloc@#1#2#3#4#5{\global\advance\count1#1by\@ne
 \ch@ck#1#4#2\allocationnumber=\count1#1
 \global#3#5=\allocationnumber
 \edef\next@{\string#5=\string#2\the\allocationnumber}%
 \expandafter\rightappend@\next@\to\alloclist@}
\newcount\count@@
\newcount\count@@@
\def\FN@{\futurelet\next}
\def\DN@{\def\next@}
\def\DNii@{\def\nextii@}
\def\RIfM@{\relax\ifmmode}
\def\RIfMIfI@{\relax\ifmmode\ifinner}
\def\setboxz@h{\setbox\z@\hbox}
\def\wdz@{\wd\z@}
\def\boxz@{\box\z@}
\def\setbox@ne{\setbox\@ne}
\def\wd@ne{\wd\@ne}
\def\iterate{\body\expandafter\iterate\else\fi}
\def\err@#1{\errmessage{AmS-TeX error: #1}}
\newhelp\defaulthelp@{Sorry, I already gave what help I could...^^J
Maybe you should try asking a human?^^J
An error might have occurred before I noticed any problems.^^J
``If all else fails, read the instructions.''}
\def\Err@{\errhelp\defaulthelp@\err@}
\def\eat@#1{}
\def\in@#1#2{\def\in@@##1#1##2##3\in@@{\ifx\in@##2\in@false\else\in@true\fi}%
 \in@@#2#1\in@\in@@}
\newif\ifin@
\def\space@.{\futurelet\space@\relax}
\space@. %
\newhelp\athelp@
{Only certain combinations beginning with @ make sense to me.^^J
Perhaps you wanted \string\@\space for a printed @?^^J
I've ignored the character or group after @.}
{\catcode`\~=\active % just in case
 \lccode`\~=`\@ \lowercase{\gdef~{\FN@\at@}}}
\def\at@{\let\next@\at@@
 \ifcat\noexpand\next a\else\ifcat\noexpand\next0\else
 \ifcat\noexpand\next\relax\else
   \let\next\at@@@\fi\fi\fi
 \next@}
\def\at@@#1{\expandafter
 \ifx\csname\space @\string#1\endcsname\relax
  \expandafter\at@@@ \else
  \csname\space @\string#1\expandafter\endcsname\fi}
\def\at@@@#1{\errhelp\athelp@ \err@{\Invalid@@ @}}%%
\def\atdef@#1{\expandafter\def\csname\space @\string#1\endcsname}%%
\newhelp\defahelp@{If you typed \string\define\space cs instead of
\string\define\string\cs\space^^J
I've substituted an inaccessible control sequence so that your^^J
definition will be completed without mixing me up too badly.^^J
If you typed \string\define{\string\cs} the inaccessible control sequence^^J
was defined to be \string\cs, and the rest of your^^J
definition appears as input.}
\newhelp\defbhelp@{I've ignored your definition, because it might^^J
conflict with other uses that are important to me.}
\def\define{\FN@\define@}
\def\define@{\ifcat\noexpand\next\relax
 \expandafter\define@@\else\errhelp\defahelp@                               %1
 \err@{\string\define\space must be followed by a control
 sequence}\expandafter\def\expandafter\nextii@\fi}                          %2
\def\undefined@@@@@@@@@@{}
\def\preloaded@@@@@@@@@@{}
\def\next@@@@@@@@@@{}
\def\define@@#1{\ifx#1\relax\errhelp\defbhelp@                              %1
 \err@{\string#1\space is already defined}\DN@{\DNii@}\else
 \expandafter\ifx\csname\expandafter\eat@\string                            %2
 #1@@@@@@@@@@\endcsname\undefined@@@@@@@@@@\errhelp\defbhelp@
 \err@{\string#1\space can't be defined}\DN@{\DNii@}\else
 \expandafter\ifx\csname\expandafter\eat@\string#1\endcsname\relax          %3
 \global\let#1\undefined\DN@{\def#1}\else\errhelp\defbhelp@
 \err@{\string#1\space is already defined}\DN@{\DNii@}\fi
 \fi\fi\next@}

\def\predefine#1#2{\let#1#2}
\def\undefine#1{\let#1\undefined}
\message{page layout,}
\newdimen\captionwidth@
\captionwidth@\hsize
\advance\captionwidth@-1.5in
\def\pagewidth#1{\hsize#1\relax
 \captionwidth@\hsize\advance\captionwidth@-1.5in}
\def\pageheight#1{\vsize#1\relax}
\def\hcorrection#1{\advance\hoffset#1\relax}
\def\vcorrection#1{\advance\voffset#1\relax}
\message{accents/punctuation,}

\let\graveaccent\`
\let\acuteaccent\'
\let\tildeaccent\~
\let\hataccent\^
\let\underscore\_
\let\B\=
\let\D\.
\let\ic@\/
\def\/{\unskip\ic@}
\def\textfonti{\the\textfont\@ne}
\def\t#1#2{{\edef\next@{\the\font}\textfonti\accent"7F \next@#1#2}}
\def~{\unskip\nobreak\ \ignorespaces}
\def\.{.\spacefactor\@m}
\atdef@;{\leavevmode\null;}
\atdef@:{\leavevmode\null:}
\atdef@?{\leavevmode\null?}
\edef\@{\string @}
\def\relaxnext@{\let\next\relax}
\atdef@-{\relaxnext@\leavevmode
 \DN@{\ifx\next-\DN@-{\FN@\nextii@}\else
  \DN@{\leavevmode\hbox{-}}\fi\next@}%
 \DNii@{\ifx\next-\DN@-{\leavevmode\hbox{---}}\else
  \DN@{\leavevmode\hbox{--}}\fi\next@}%
 \FN@\next@}
\def\srdr@{\kern.16667em}
\def\drsr@{\kern.02778em}
\def\sldl@{\drsr@}
\def\dlsl@{\srdr@}
\atdef@"{\unskip\relaxnext@
 \DN@{\ifx\next\space@\DN@. {\FN@\nextii@}\else
  \DN@.{\FN@\nextii@}\fi\next@.}%
 \DNii@{\ifx\next`\DN@`{\FN@\nextiii@}\else
  \ifx\next\lq\DN@\lq{\FN@\nextiii@}\else
  \DN@####1{\FN@\nextiv@}\fi\fi\next@}%
 \def\nextiii@{\ifx\next`\DN@`{\sldl@``}\else\ifx\next\lq
  \DN@\lq{\sldl@``}\else\DN@{\dlsl@`}\fi\fi\next@}%
 \def\nextiv@{\ifx\next'\DN@'{\srdr@''}\else
  \ifx\next\rq\DN@\rq{\srdr@''}\else\DN@{\drsr@'}\fi\fi\next@}%
 \FN@\next@}

\def\textfontii{\the\textfont\tw@}
\def\lbrace@{\delimiter"4266308 }
\def\rbrace@{\delimiter"5267309 }
\def\{{\RIfM@\lbrace@\else{\textfontii f}\spacefactor\@m\fi}
\def\}{\RIfM@\rbrace@\else
 \let\@sf\empty\ifhmode\edef\@sf{\spacefactor\the\spacefactor}\fi
 {\textfontii g}\@sf\relax\fi}
\let\lbrace\{
\let\rbrace\}
\def\AmSTeX{{\textfontii A\kern-.1667em%
  \lower.5ex\hbox{M}\kern-.125emS}-\TeX\spacefactor1000 }
\message{line and page breaks,}
\def\vmodeerr@#1{\Err@{\string#1\space not allowed between paragraphs}}
\def\mathmodeerr@#1{\Err@{\string#1\space not allowed in math mode}}
\def\linebreak{\RIfM@\mathmodeerr@\linebreak\else
 \ifhmode\unskip\unkern\break\else\vmodeerr@\linebreak\fi\fi}

\newskip\saveskip@
\def\allowlinebreak{\RIfM@\mathmodeerr@\allowlinebreak\else
 \ifhmode\saveskip@\lastskip\unskip
 \allowbreak\ifdim\saveskip@>\z@\hskip\saveskip@\fi
 \else\vmodeerr@\allowlinebreak\fi\fi}
\def\nolinebreak{\RIfM@\mathmodeerr@\nolinebreak\else
 \ifhmode\saveskip@\lastskip\unskip
 \nobreak\ifdim\saveskip@>\z@\hskip\saveskip@\fi
 \else\vmodeerr@\nolinebreak\fi\fi}
\def\newline{\relaxnext@
 \DN@{\RIfM@\expandafter\mathmodeerr@\expandafter\newline\else
  \ifhmode\ifx\next\par\else
  \expandafter\unskip\expandafter\null\expandafter\hfill\expandafter\break\fi
  \else
  \expandafter\vmodeerr@\expandafter\newline\fi\fi}%
 \FN@\next@}
\def\dmatherr@#1{\Err@{\string#1\space not allowed in display math mode}}
\def\nondmatherr@#1{\Err@{\string#1\space not allowed in non-display math
 mode}}
\def\onlydmatherr@#1{\Err@{\string#1\space allowed only in display math mode}}
\def\nonmatherr@#1{\Err@{\string#1\space allowed only in math mode}}
\def\mathbreak{\RIfMIfI@\break\else
 \dmatherr@\mathbreak\fi\else\nonmatherr@\mathbreak\fi}
\def\nomathbreak{\RIfMIfI@\nobreak\else
 \dmatherr@\nomathbreak\fi\else\nonmatherr@\nomathbreak\fi}
\def\allowmathbreak{\RIfMIfI@\allowbreak\else
 \dmatherr@\allowmathbreak\fi\else\nonmatherr@\allowmathbreak\fi}
\def\pagebreak{\RIfM@
 \ifinner\nondmatherr@\pagebreak\else\postdisplaypenalty-\@M\fi
 \else\ifvmode\removelastskip\break\else\vadjust{\break}\fi\fi}
\def\nopagebreak{\RIfM@
 \ifinner\nondmatherr@\nopagebreak\else\postdisplaypenalty\@M\fi
 \else\ifvmode\nobreak\else\vadjust{\nobreak}\fi\fi}
\def\nonvmodeerr@#1{\Err@{\string#1\space not allowed within a paragraph
 or in math}}
\def\vnonvmode@#1#2{\relaxnext@\DNii@{\ifx\next\par\DN@{#1}\else
 \DN@{#2}\fi\next@}%
 \ifvmode\DN@{#1}\else
 \DN@{\FN@\nextii@}\fi\next@}
\def\newpage{\vnonvmode@{\vfill\break}{\nonvmodeerr@\newpage}}
\def\smallpagebreak{\vnonvmode@\smallbreak{\nonvmodeerr@\smallpagebreak}}
\def\medpagebreak{\vnonvmode@\medbreak{\nonvmodeerr@\medpagebreak}}
\def\bigpagebreak{\vnonvmode@\bigbreak{\nonvmodeerr@\bigpagebreak}}
\def\NoBlackBoxes{\global\overfullrule\z@}
\def\BlackBoxes{\global\overfullrule5\p@}
\def\Invalid@#1{\def#1{\Err@{\Invalid@@\string#1}}}
\def\Invalid@@{Invalid use of }
\message{figures,}
\Invalid@\caption
\Invalid@\captionwidth
\newdimen\smallcaptionwidth@
\def\topspace{\mid@false\ins@}
\def\midspace{\mid@true\ins@}
\newif\ifmid@
\def\captionfont@{}
\def\ins@#1{\relaxnext@\allowbreak
 \smallcaptionwidth@\captionwidth@\gdef\thespace@{#1}%
 \DN@{\ifx\next\space@\DN@. {\FN@\nextii@}\else
  \DN@.{\FN@\nextii@}\fi\next@.}%
 \DNii@{\ifx\next\caption\DN@\caption{\FN@\nextiii@}%
  \else\let\next@\nextiv@\fi\next@}%
 \def\nextiv@{\vnonvmode@
  {\ifmid@\expandafter\midinsert\else\expandafter\topinsert\fi
   \vbox to\thespace@{}\endinsert}
  {\ifmid@\nonvmodeerr@\midspace\else\nonvmodeerr@\topspace\fi}}%
 \def\nextiii@{\ifx\next\captionwidth\expandafter\nextv@
  \else\expandafter\nextvi@\fi}%
 \def\nextv@\captionwidth##1##2{\smallcaptionwidth@##1\relax\nextvi@{##2}}%
 \def\nextvi@##1{\def\thecaption@{\captionfont@##1}%
  \DN@{\ifx\next\space@\DN@. {\FN@\nextvii@}\else
   \DN@.{\FN@\nextvii@}\fi\next@.}%
  \FN@\next@}%
 \def\nextvii@{\vnonvmode@
  {\ifmid@\expandafter\midinsert\else
  \expandafter\topinsert\fi\vbox to\thespace@{}\nobreak\smallskip
  \setboxz@h{\noindent\ignorespaces\thecaption@\unskip}%
  \ifdim\wdz@>\smallcaptionwidth@\centerline{\vbox{\hsize\smallcaptionwidth@
   \noindent\ignorespaces\thecaption@\unskip}}%
  \else\centerline{\boxz@}\fi\endinsert}
  {\ifmid@\nonvmodeerr@\midspace
  \else\nonvmodeerr@\topspace\fi}}%
 \FN@\next@}
\message{comments,}
\def\newcodes@{\catcode`\\12\catcode`\{12\catcode`\}12\catcode`\#12%
 \catcode`\%12\relax}
\def\oldcodes@{\catcode`\\0\catcode`\{1\catcode`\}2\catcode`\#6%
 \catcode`\%14\relax}
\def\comment{\newcodes@\endlinechar=10 \comment@}
{\lccode`\0=`\\
\lowercase{\gdef\comment@#1^^J{\comment@@#10endcomment\comment@@@}%
\gdef\comment@@#10endcomment{\FN@\comment@@@}%
\gdef\comment@@@#1\comment@@@{\ifx\next\comment@@@\let\next\comment@
 \else\def\next{\oldcodes@\endlinechar=`\^^M\relax}%
 \fi\next}}}
\def\pr@m@s{\ifx'\next\DN@##1{\prim@s}\else\let\next@\egroup\fi\next@}
\def\prime{{\null\prime@\null}}
\mathchardef\prime@="0230
\let\dsize\displaystyle

\let\ssize\scriptstyle

\message{math spacing,}
\def\,{\RIfM@\mskip\thinmuskip\relax\else\kern.16667em\fi}
\def\!{\RIfM@\mskip-\thinmuskip\relax\else\kern-.16667em\fi}
\let\thinspace\,
\let\negthinspace\!
\def\medspace{\RIfM@\mskip\medmuskip\relax\else\kern.222222em\fi}
\def\negmedspace{\RIfM@\mskip-\medmuskip\relax\else\kern-.222222em\fi}
\def\thickspace{\RIfM@\mskip\thickmuskip\relax\else\kern.27777em\fi}
\let\;\thickspace
\def\negthickspace{\RIfM@\mskip-\thickmuskip\relax\else
 \kern-.27777em\fi}
\atdef@,{\RIfM@\mskip.1\thinmuskip\else\leavevmode\null,\fi}
\atdef@!{\RIfM@\mskip-.1\thinmuskip\else\leavevmode\null!\fi}
\atdef@.{\RIfM@&&\else\leavevmode.\spacefactor3000 \fi}
\def\and{\DOTSB\;\mathchar"3026 \;}

\message{fractions,}
\def\frac#1#2{{#1\over#2}}

\newdimen\ex@
\ex@.2326ex
\Invalid@\thickness
\def\thickfrac{\relaxnext@
 \DN@{\ifx\next\thickness\let\next@\nextii@\else
 \DN@{\nextii@\thickness1}\fi\next@}%
 \DNii@\thickness##1##2##3{{##2\above##1\ex@##3}}%
 \FN@\next@}

\def\thickfracwithdelims#1#2{\relaxnext@\def\ldelim@{#1}\def\rdelim@{#2}%
 \DN@{\ifx\next\thickness\let\next@\nextii@\else
 \DN@{\nextii@\thickness1}\fi\next@}%
 \DNii@\thickness##1##2##3{{##2\abovewithdelims
 \ldelim@\rdelim@##1\ex@##3}}%
 \FN@\next@}

\def\:{\nobreak\hskip.1111em\mathpunct{}\nonscript\mkern-\thinmuskip{:}\hskip
 .3333emplus.0555em\relax}
\def\snug{\unskip\kern-\mathsurround}
\message{smash commands,}
\def\topsmash{\top@true\bot@false\smash@}
\def\botsmash{\top@false\bot@true\smash@}
\newif\iftop@
\newif\ifbot@
\def\smash{\top@true\bot@true\smash@}
\def\smash@{\RIfM@\expandafter\mathpalette\expandafter\mathsm@sh\else
 \expandafter\makesm@sh\fi}
\def\finsm@sh{\iftop@\ht\z@\z@\fi\ifbot@\dp\z@\z@\fi\leavevmode\boxz@}
\message{large operator symbols,}
\def\LimitsOnSums{\global\let\slimits@\displaylimits}
\def\NoLimitsOnSums{\global\let\slimits@\nolimits}
\LimitsOnSums
\mathchardef\coprod@="1360       \def\coprod{\DOTSB\coprod@\slimits@}
\mathchardef\bigvee@="1357       \def\bigvee{\DOTSB\bigvee@\slimits@}
\mathchardef\bigwedge@="1356     \def\bigwedge{\DOTSB\bigwedge@\slimits@}
\mathchardef\biguplus@="1355     \def\biguplus{\DOTSB\biguplus@\slimits@}
\mathchardef\bigcap@="1354       \def\bigcap{\DOTSB\bigcap@\slimits@}
\mathchardef\bigcup@="1353       \def\bigcup{\DOTSB\bigcup@\slimits@}
\mathchardef\prod@="1351         \def\prod{\DOTSB\prod@\slimits@}
\mathchardef\sum@="1350          \def\sum{\DOTSB\sum@\slimits@}
\mathchardef\bigotimes@="134E    \def\bigotimes{\DOTSB\bigotimes@\slimits@}
\mathchardef\bigoplus@="134C     \def\bigoplus{\DOTSB\bigoplus@\slimits@}
\mathchardef\bigodot@="134A      \def\bigodot{\DOTSB\bigodot@\slimits@}
\mathchardef\bigsqcup@="1346     \def\bigsqcup{\DOTSB\bigsqcup@\slimits@}
\message{integrals,}
\def\LimitsOnInts{\global\let\ilimits@\displaylimits}
\def\NoLimitsOnInts{\global\let\ilimits@\nolimits}
\NoLimitsOnInts
\def\int{\DOTSI\intop\ilimits@}
\def\oint{\DOTSI\ointop\ilimits@}
\def\intic@{\mathchoice{\hskip.5em}{\hskip.4em}{\hskip.4em}{\hskip.4em}}
\def\negintic@{\mathchoice
 {\hskip-.5em}{\hskip-.4em}{\hskip-.4em}{\hskip-.4em}}
\def\intkern@{\mathchoice{\!\!\!}{\!\!}{\!\!}{\!\!}}
\def\intdots@{\mathchoice{\plaincdots@}
 {{\cdotp}\mkern1.5mu{\cdotp}\mkern1.5mu{\cdotp}}
 {{\cdotp}\mkern1mu{\cdotp}\mkern1mu{\cdotp}}
 {{\cdotp}\mkern1mu{\cdotp}\mkern1mu{\cdotp}}}
\newcount\intno@
\def\iint{\DOTSI\intno@\tw@\FN@\ints@}
\def\iiint{\DOTSI\intno@\thr@@\FN@\ints@}
\def\iiiint{\DOTSI\intno@4 \FN@\ints@}
\def\idotsint{\DOTSI\intno@\z@\FN@\ints@}
\def\ints@{\findlimits@\ints@@}
\newif\iflimtoken@
\newif\iflimits@
\def\findlimits@{\limtoken@true\ifx\next\limits\limits@true
 \else\ifx\next\nolimits\limits@false\else
 \limtoken@false\ifx\ilimits@\nolimits\limits@false\else
 \ifinner\limits@false\else\limits@true\fi\fi\fi\fi}
\def\multint@{\int\ifnum\intno@=\z@\intdots@                                %1
 \else\intkern@\fi                                                          %2
 \ifnum\intno@>\tw@\int\intkern@\fi                                         %3
 \ifnum\intno@>\thr@@\int\intkern@\fi                                       %4
 \int}                                                                      %5
\def\multintlimits@{\intop\ifnum\intno@=\z@\intdots@\else\intkern@\fi
 \ifnum\intno@>\tw@\intop\intkern@\fi
 \ifnum\intno@>\thr@@\intop\intkern@\fi\intop}
\def\ints@@{\iflimtoken@                                                    %1
 \def\ints@@@{\iflimits@\negintic@\mathop{\intic@\multintlimits@}\limits    %2
  \else\multint@\nolimits\fi                                                %3
  \eat@}                                                                    %4
 \else                                                                      %5
 \def\ints@@@{\iflimits@\negintic@
  \mathop{\intic@\multintlimits@}\limits\else
  \multint@\nolimits\fi}\fi\ints@@@}
\def\LimitsOnNames{\global\let\nlimits@\displaylimits}
\def\NoLimitsOnNames{\global\let\nlimits@\nolimits@}
\LimitsOnNames
\def\nolimits@{\relaxnext@
 \DN@{\ifx\next\limits\DN@\limits{\nolimits}\else
  \let\next@\nolimits\fi\next@}%
 \FN@\next@}
\message{operator names,}
\def\newmcodes@{\mathcode`\'"27\mathcode`\*"2A\mathcode`\."613A%
 \mathcode`\-"2D\mathcode`\/"2F\mathcode`\:"603A }
\def\operatorname#1{\mathop{\newmcodes@\kern\z@\fam\z@#1}\nolimits@}
\def\operatornamewithlimits#1{\mathop{\newmcodes@\kern\z@\fam\z@#1}\nlimits@}
\def\qopname@#1{\mathop{\fam\z@#1}\nolimits@}
\def\qopnamewl@#1{\mathop{\fam\z@#1}\nlimits@}
\def\arccos{\qopname@{arccos}}
\def\arcsin{\qopname@{arcsin}}
\def\arctan{\qopname@{arctan}}
\def\arg{\qopname@{arg}}
\def\cos{\qopname@{cos}}
\def\cosh{\qopname@{cosh}}
\def\cot{\qopname@{cot}}
\def\coth{\qopname@{coth}}
\def\csc{\qopname@{csc}}
\def\deg{\qopname@{deg}}
\def\det{\qopnamewl@{det}}
\def\dim{\qopname@{dim}}
\def\exp{\qopname@{exp}}
\def\gcd{\qopnamewl@{gcd}}
\def\hom{\qopname@{hom}}
\def\inf{\qopnamewl@{inf}}
\def\injlim{\qopnamewl@{inj\,lim}}
\def\ker{\qopname@{ker}}
\def\lg{\qopname@{lg}}
\def\lim{\qopnamewl@{lim}}
\def\liminf{\qopnamewl@{lim\,inf}}
\def\limsup{\qopnamewl@{lim\,sup}}
\def\ln{\qopname@{ln}}
\def\log{\qopname@{log}}
\def\max{\qopnamewl@{max}}
\def\min{\qopnamewl@{min}}
\def\Pr{\qopnamewl@{Pr}}
\def\projlim{\qopnamewl@{proj\,lim}}
\def\sec{\qopname@{sec}}
\def\sin{\qopname@{sin}}
\def\sinh{\qopname@{sinh}}
\def\sup{\qopnamewl@{sup}}
\def\tan{\qopname@{tan}}
\def\tanh{\qopname@{tanh}}
\def\varinjlim{\mathop{\vtop{\ialign{##\crcr
 \hfil\rm lim\hfil\crcr\noalign{\nointerlineskip}\rightarrowfill\crcr
 \noalign{\nointerlineskip\kern-\ex@}\crcr}}}}
\def\varprojlim{\mathop{\vtop{\ialign{##\crcr
 \hfil\rm lim\hfil\crcr\noalign{\nointerlineskip}\leftarrowfill\crcr
 \noalign{\nointerlineskip\kern-\ex@}\crcr}}}}
\def\varliminf{\mathop{\underline{\vrule height\z@ depth.2exwidth\z@
 \hbox{\rm lim}}}}

\newdimen\buffer@
\buffer@\fontdimen13 \tenex
\newdimen\buffer
\buffer\buffer@

\def\ResetBuffer{\fontdimen13 \tenex\buffer@\global\buffer\buffer@}
\def\shave#1{\mathop{\hbox{$\m@th\fontdimen13 \tenex\z@                     %1
 \displaystyle{#1}$}}\fontdimen13 \tenex\buffer}

\message{multilevel sub/superscripts,}
\Invalid@\\
\def\Let@{\relax\iffalse{\fi\let\\=\cr\iffalse}\fi}
\Invalid@\vspace
\def\vspace@{\def\vspace##1{\crcr\noalign{\vskip##1\relax}}}
\def\multilimits@{\bgroup\vspace@\Let@
 \baselineskip\fontdimen10 \scriptfont\tw@
 \advance\baselineskip\fontdimen12 \scriptfont\tw@
 \lineskip\thr@@\fontdimen8 \scriptfont\thr@@
 \lineskiplimit\lineskip
 \vbox\bgroup\ialign\bgroup\hfil$\m@th\scriptstyle{##}$\hfil\crcr}
\def\Sb{_\multilimits@}
\def\endSb{\crcr\egroup\egroup\egroup}
\def\Sp{^\multilimits@}

\def\spreadlines#1{\RIfMIfI@\onlydmatherr@\spreadlines\else
 \openup#1\relax\fi\else\onlydmatherr@\spreadlines\fi}
\def\Mathstrut@{\copy\Mathstrutbox@}
\newbox\Mathstrutbox@
\setbox\Mathstrutbox@\null
\setboxz@h{$\m@th($}
\ht\Mathstrutbox@\ht\z@
\dp\Mathstrutbox@\dp\z@
\message{matrices,}
\newdimen\spreadmlines@
\def\spreadmatrixlines#1{\RIfMIfI@
 \onlydmatherr@\spreadmatrixlines\else
 \spreadmlines@#1\relax\fi\else\onlydmatherr@\spreadmatrixlines\fi}
\def\matrix{\null\,\vcenter\bgroup\Let@\vspace@
 \normalbaselines\openup\spreadmlines@\ialign
 \bgroup\hfil$\m@th##$\hfil&&\quad\hfil$\m@th##$\hfil\crcr
 \Mathstrut@\crcr\noalign{\kern-\baselineskip}}
\def\endmatrix{\crcr\Mathstrut@\crcr\noalign{\kern-\baselineskip}\egroup
 \egroup\,}
\def\format{\crcr\egroup\iffalse{\fi\ifnum`}=0 \fi\format@}
\newtoks\hashtoks@
\hashtoks@{#}
\def\format@#1\\{\def\preamble@{#1}%
 \def\l{$\m@th\the\hashtoks@$\hfil}%
 \def\c{\hfil$\m@th\the\hashtoks@$\hfil}%
 \def\r{\hfil$\m@th\the\hashtoks@$}%
 \edef\preamble@@{\preamble@}\ifnum`{=0 \fi\iffalse}\fi
 \ialign\bgroup\span\preamble@@\crcr}
\def\smallmatrix{\null\,\vcenter\bgroup\vspace@\Let@
 \baselineskip9\ex@\lineskip\ex@
 \ialign\bgroup\hfil$\m@th\scriptstyle{##}$\hfil&&\thickspace\hfil
 $\m@th\scriptstyle{##}$\hfil\crcr}
\def\endsmallmatrix{\crcr\egroup\egroup\,}

\newmuskip\dotsspace@
\dotsspace@1.5mu
\def\strip@#1 {#1}
\def\spacehdots#1\for#2{\multispan{#2}\xleaders
 \hbox{$\m@th\mkern\strip@#1 \dotsspace@.\mkern\strip@#1 \dotsspace@$}\hfill}
\def\hdotsfor#1{\spacehdots\@ne\for{#1}}
\def\multispan@#1{\omit\mscount#1\unskip\loop\ifnum\mscount>\@ne\sp@n\repeat}
\def\spaceinnerhdots#1\for#2\after#3{\multispan@{\strip@#2 }#3\xleaders
 \hbox{$\m@th\mkern\strip@#1 \dotsspace@.\mkern\strip@#1 \dotsspace@$}\hfill}
\def\innerhdotsfor#1\after#2{\spaceinnerhdots\@ne\for#1\after{#2}}
\def\cases{\bgroup\spreadmlines@\jot\left\{\,\matrix\format\l&\quad\l\\}
\def\endcases{\endmatrix\right.\egroup}
\message{multiline displays,}
\newif\ifinany@
\newif\ifinalign@
\newif\ifingather@
\def\strut@{\copy\strutbox@}
\newbox\strutbox@
\setbox\strutbox@\hbox{\vrule height8\p@ depth3\p@ width\z@}
\def\topaligned{\null\,\vtop\aligned@}
\def\botaligned{\null\,\vbox\aligned@}
\def\aligned{\null\,\vcenter\aligned@}
\def\aligned@{\bgroup\vspace@\Let@
 \ifinany@\else\openup\jot\fi\ialign
 \bgroup\hfil\strut@$\m@th\displaystyle{##}$&
 $\m@th\displaystyle{{}##}$\hfil\crcr}
\def\endaligned{\crcr\egroup\egroup}

\def\alignedat#1{\null\,\vcenter\bgroup\doat@{#1}\vspace@\Let@
 \ifinany@\else\openup\jot\fi\ialign\bgroup\span\preamble@@\crcr}
\newcount\atcount@
\def\doat@#1{\toks@{\hfil\strut@$\m@th
 \displaystyle{\the\hashtoks@}$&$\m@th\displaystyle
 {{}\the\hashtoks@}$\hfil}%                                                 %1
 \atcount@#1\relax\advance\atcount@\m@ne                                    %2
 \loop\ifnum\atcount@>\z@\toks@=\expandafter{\the\toks@&\hfil$\m@th
 \displaystyle{\the\hashtoks@}$&$\m@th
 \displaystyle{{}\the\hashtoks@}$\hfil}\advance
  \atcount@\m@ne\repeat                                                     %3
 \xdef\preamble@{\the\toks@}\xdef\preamble@@{\preamble@}}

\def\gathered{\null\,\vcenter\bgroup\vspace@\Let@
 \ifinany@\else\openup\jot\fi\ialign
 \bgroup\hfil\strut@$\m@th\displaystyle{##}$\hfil\crcr}
\def\endgathered{\crcr\egroup\egroup}
\newif\iftagsleft@
\def\TagsOnLeft{\global\tagsleft@true}
\def\TagsOnRight{\global\tagsleft@false}
\TagsOnLeft
\newif\ifmathtags@
\def\TagsAsMath{\global\mathtags@true}
\def\TagsAsText{\global\mathtags@false}
\TagsAsText
\def\tagform@#1{\hbox{\rm(\ignorespaces#1\unskip)}}
\def\thetag{\leavevmode\tagform@}
\def\tag#1$${\iftagsleft@\leqno\else\eqno\fi                                %1
 \maketag@#1\maketag@                                                       %2
 $$}                                                                        %3
\def\maketag@{\FN@\maketag@@}
\def\maketag@@{\ifx\next"\expandafter\maketag@@@\else\expandafter\maketag@@@@
 \fi}
\def\maketag@@@"#1"#2\maketag@{\hbox{\rm#1}}                                %1
\def\maketag@@@@#1\maketag@{\ifmathtags@\tagform@{$\m@th#1$}\else
 \tagform@{#1}\fi}
\interdisplaylinepenalty\@M
\def\allowdisplaybreaks{\RIfMIfI@
 \onlydmatherr@\allowdisplaybreaks\else
 \interdisplaylinepenalty\z@\fi\else\onlydmatherr@\allowdisplaybreaks\fi}
\Invalid@\allowdisplaybreak
\Invalid@\displaybreak
\Invalid@\intertext
\def\allowdisplaybreak@{\def\allowdisplaybreak{\crcr\noalign{\allowbreak}}}
\def\displaybreak@{\def\displaybreak{\crcr\noalign{\break}}}
\def\intertext@{\def\intertext##1{\crcr\noalign{%
 \penalty\postdisplaypenalty \vskip\belowdisplayskip
 \vbox{\normalbaselines\noindent##1}%
 \penalty\predisplaypenalty \vskip\abovedisplayskip}}}
\newskip\centering@
\centering@\z@ plus\@m\p@
\def\align{\relax\ifingather@\DN@{\csname align (in
  \string\gather)\endcsname}\else
 \ifmmode\ifinner\DN@{\onlydmatherr@\align}\else
  \let\next@\align@\fi
 \else\DN@{\onlydmatherr@\align}\fi\fi\next@}
\newhelp\andhelp@
{An extra & here is so disastrous that you should probably exit^^J
and fix things up.}
\newif\iftag@
\newcount\and@
\def\align@{\inalign@true\inany@true
 \vspace@\allowdisplaybreak@\displaybreak@\intertext@
 \def\tag{\global\tag@true\ifnum\and@=\z@\DN@{&&}\else
          \DN@{&}\fi\next@}%
 \iftagsleft@\DN@{\csname align \endcsname}\else
  \DN@{\csname align \space\endcsname}\fi\next@}
\def\Tag@{\iftag@\else\errhelp\andhelp@\err@{Extra & on this line}\fi}
\newdimen\lwidth@
\newdimen\rwidth@
\newdimen\maxlwidth@
\newdimen\maxrwidth@
\newdimen\totwidth@
\def\measure@#1\endalign{\lwidth@\z@\rwidth@\z@\maxlwidth@\z@\maxrwidth@\z@
 \global\and@\z@                                                            %1
 \setbox@ne\vbox                                                            %2
  {\everycr{\noalign{\global\tag@false\global\and@\z@}}\Let@                %3
  \halign{\setboxz@h{$\m@th\displaystyle{\@lign##}$}%                       %4
   \global\lwidth@\wdz@                                                     %5
   \ifdim\lwidth@>\maxlwidth@\global\maxlwidth@\lwidth@\fi                  %6
   \global\advance\and@\@ne                                                 %7
   &\setboxz@h{$\m@th\displaystyle{{}\@lign##}$}\global\rwidth@\wdz@        %8
   \ifdim\rwidth@>\maxrwidth@\global\maxrwidth@\rwidth@\fi                  %9
   \global\advance\and@\@ne                                                %10
   &\Tag@
   \eat@{##}\crcr#1\crcr}}%                                                %11
 \totwidth@\maxlwidth@\advance\totwidth@\maxrwidth@}                       %12
\def\displ@y@{\global\dt@ptrue\openup\jot
 \everycr{\noalign{\global\tag@false\global\and@\z@\ifdt@p\global\dt@pfalse
 \vskip-\lineskiplimit\vskip\normallineskiplimit\else
 \penalty\interdisplaylinepenalty\fi}}}
\def\black@#1{\noalign{\ifdim#1>\displaywidth
 \dimen@\prevdepth\nointerlineskip                                          %1
 \vskip-\ht\strutbox@\vskip-\dp\strutbox@                                   %2
 \vbox{\noindent\hbox to#1{\strut@\hfill}}%                                 %3
 \prevdepth\dimen@                                                          %4
 \fi}}
\expandafter\def\csname align \space\endcsname#1\endalign
 {\measure@#1\endalign\global\and@\z@                                       %1
 \ifingather@\everycr{\noalign{\global\and@\z@}}\else\displ@y@\fi           %2
 \Let@\tabskip\centering@                                                   %3
 \halign to\displaywidth
  {\hfil\strut@\setboxz@h{$\m@th\displaystyle{\@lign##}$}%                  %4
  \global\lwidth@\wdz@\boxz@\global\advance\and@\@ne                        %5
  \tabskip\z@skip                                                           %6
  &\setboxz@h{$\m@th\displaystyle{{}\@lign##}$}%                            %7
  \global\rwidth@\wdz@\boxz@\hfill\global\advance\and@\@ne                  %8
  \tabskip\centering@                                                       %9
  &\setboxz@h{\@lign\strut@\maketag@##\maketag@}%                          %10
  \dimen@\displaywidth\advance\dimen@-\totwidth@
  \divide\dimen@\tw@\advance\dimen@\maxrwidth@\advance\dimen@-\rwidth@     %11
  \ifdim\dimen@<\tw@\wdz@\llap{\vtop{\normalbaselines\null\boxz@}}%        %12
  \else\llap{\boxz@}\fi                                                    %13
  \tabskip\z@skip                                                          %14
  \crcr#1\crcr                                                             %15
  \black@\totwidth@}}                                                      %16
\newdimen\lineht@
\expandafter\def\csname align \endcsname#1\endalign{\measure@#1\endalign
 \global\and@\z@
 \ifdim\totwidth@>\displaywidth\let\displaywidth@\totwidth@\else
  \let\displaywidth@\displaywidth\fi                                        %1
 \ifingather@\everycr{\noalign{\global\and@\z@}}\else\displ@y@\fi
 \Let@\tabskip\centering@\halign to\displaywidth
  {\hfil\strut@\setboxz@h{$\m@th\displaystyle{\@lign##}$}%
  \global\lwidth@\wdz@\global\lineht@\ht\z@                                 %2
  \boxz@\global\advance\and@\@ne
  \tabskip\z@skip&\setboxz@h{$\m@th\displaystyle{{}\@lign##}$}%
  \global\rwidth@\wdz@\ifdim\ht\z@>\lineht@\global\lineht@\ht\z@\fi         %3
  \boxz@\hfil\global\advance\and@\@ne
  \tabskip\centering@&\kern-\displaywidth@                                  %4
  \setboxz@h{\@lign\strut@\maketag@##\maketag@}%
  \dimen@\displaywidth\advance\dimen@-\totwidth@
  \divide\dimen@\tw@\advance\dimen@\maxlwidth@\advance\dimen@-\lwidth@
  \ifdim\dimen@<\tw@\wdz@
   \rlap{\vbox{\normalbaselines\boxz@\vbox to\lineht@{}}}\else
   \rlap{\boxz@}\fi
  \tabskip\displaywidth@\crcr#1\crcr\black@\totwidth@}}
\expandafter\def\csname align (in \string\gather)\endcsname
  #1\endalign{\vcenter{\align@#1\endalign}}
\Invalid@\endalign
\newif\ifxat@
\def\alignat{\RIfMIfI@\DN@{\onlydmatherr@\alignat}\else
 \DN@{\csname alignat \endcsname}\fi\else
 \DN@{\onlydmatherr@\alignat}\fi\next@}
\newif\ifmeasuring@
\newbox\savealignat@
\expandafter\def\csname alignat \endcsname#1#2\endalignat                   %1
 {\inany@true\xat@false
 \def\tag{\global\tag@true\count@#1\relax\multiply\count@\tw@
  \xdef\tag@{}\loop\ifnum\count@>\and@\xdef\tag@{&\tag@}\advance\count@\m@ne
  \repeat\tag@}%
 \vspace@\allowdisplaybreak@\displaybreak@\intertext@
 \displ@y@\measuring@true                                                   %2
 \setbox\savealignat@\hbox{$\m@th\displaystyle\Let@
  \attag@{#1}%                                                              %3
  \vbox{\halign{\span\preamble@@\crcr#2\crcr}}$}%
 \measuring@false                                                           %4
 \Let@\attag@{#1}%                                                          %5
 \tabskip\centering@\halign to\displaywidth
  {\span\preamble@@\crcr#2\crcr                                             %6
  \black@{\wd\savealignat@}}}                                               %7
\Invalid@\endalignat
\def\xalignat{\RIfMIfI@
 \DN@{\onlydmatherr@\xalignat}\else
 \DN@{\csname xalignat \endcsname}\fi\else
 \DN@{\onlydmatherr@\xalignat}\fi\next@}
\expandafter\def\csname xalignat \endcsname#1#2\endxalignat
 {\inany@true\xat@true
 \def\tag{\global\tag@true\def\tag@{}\count@#1\relax\multiply\count@\tw@
  \loop\ifnum\count@>\and@\xdef\tag@{&\tag@}\advance\count@\m@ne\repeat\tag@}%
 \vspace@\allowdisplaybreak@\displaybreak@\intertext@
 \displ@y@\measuring@true\setbox\savealignat@\hbox{$\m@th\displaystyle\Let@
 \attag@{#1}\vbox{\halign{\span\preamble@@\crcr#2\crcr}}$}%
 \measuring@false\Let@
 \attag@{#1}\tabskip\centering@\halign to\displaywidth
 {\span\preamble@@\crcr#2\crcr\black@{\wd\savealignat@}}}
\def\attag@#1{\let\Maketag@\maketag@\let\TAG@\Tag@                          %1
 \let\Tag@=0\let\maketag@=0%                                                %2
 \ifmeasuring@\def\llap@##1{\setboxz@h{##1}\hbox to\tw@\wdz@{}}%
  \def\rlap@##1{\setboxz@h{##1}\hbox to\tw@\wdz@{}}\else
  \let\llap@\llap\let\rlap@\rlap\fi                                         %3
 \toks@{\hfil\strut@$\m@th\displaystyle{\@lign\the\hashtoks@}$\tabskip\z@skip
  \global\advance\and@\@ne&$\m@th\displaystyle{{}\@lign\the\hashtoks@}$\hfil
  \ifxat@\tabskip\centering@\fi\global\advance\and@\@ne}%                   %4
 \iftagsleft@
  \toks@@{\tabskip\centering@&\Tag@\kern-\displaywidth
   \rlap@{\@lign\maketag@\the\hashtoks@\maketag@}%
   \global\advance\and@\@ne\tabskip\displaywidth}\else
  \toks@@{\tabskip\centering@&\Tag@\llap@{\@lign\maketag@
   \the\hashtoks@\maketag@}\global\advance\and@\@ne\tabskip\z@skip}\fi      %5
 \atcount@#1\relax\advance\atcount@\m@ne
 \loop\ifnum\atcount@>\z@
 \toks@=\expandafter{\the\toks@&\hfil$\m@th\displaystyle{\@lign
  \the\hashtoks@}$\global\advance\and@\@ne
  \tabskip\z@skip&$\m@th\displaystyle{{}\@lign\the\hashtoks@}$\hfil\ifxat@
  \tabskip\centering@\fi\global\advance\and@\@ne}\advance\atcount@\m@ne
 \repeat                                                                    %6
 \xdef\preamble@{\the\toks@\the\toks@@}%                                    %7
 \xdef\preamble@@{\preamble@}%                                              %8
 \let\maketag@\Maketag@\let\Tag@\TAG@}                                      %9
\Invalid@\endxalignat
\def\xxalignat{\RIfMIfI@
 \DN@{\onlydmatherr@\xxalignat}\else\DN@{\csname xxalignat
  \endcsname}\fi\else
 \DN@{\onlydmatherr@\xxalignat}\fi\next@}
\expandafter\def\csname xxalignat \endcsname#1#2\endxxalignat{\inany@true
 \vspace@\allowdisplaybreak@\displaybreak@\intertext@
 \displ@y\setbox\savealignat@\hbox{$\m@th\displaystyle\Let@
 \xxattag@{#1}\vbox{\halign{\span\preamble@@\crcr#2\crcr}}$}%
 \Let@\xxattag@{#1}\tabskip\z@skip\halign to\displaywidth
 {\span\preamble@@\crcr#2\crcr\black@{\wd\savealignat@}}}
\def\xxattag@#1{\toks@{\tabskip\z@skip\hfil\strut@
 $\m@th\displaystyle{\the\hashtoks@}$&%
 $\m@th\displaystyle{{}\the\hashtoks@}$\hfil\tabskip\centering@&}%
 \atcount@#1\relax\advance\atcount@\m@ne\loop\ifnum\atcount@>\z@
 \toks@=\expandafter{\the\toks@&\hfil$\m@th\displaystyle{\the\hashtoks@}$%
  \tabskip\z@skip&$\m@th\displaystyle{{}\the\hashtoks@}$\hfil
  \tabskip\centering@}\advance\atcount@\m@ne\repeat
 \xdef\preamble@{\the\toks@\tabskip\z@skip}\xdef\preamble@@{\preamble@}}
\Invalid@\endxxalignat
\newdimen\gwidth@
\newdimen\gmaxwidth@
\def\gmeasure@#1\endgather{\gwidth@\z@\gmaxwidth@\z@\setbox@ne\vbox{\Let@
 \halign{\setboxz@h{$\m@th\displaystyle{##}$}\global\gwidth@\wdz@
 \ifdim\gwidth@>\gmaxwidth@\global\gmaxwidth@\gwidth@\fi
 &\eat@{##}\crcr#1\crcr}}}
\def\gather{\RIfMIfI@\DN@{\onlydmatherr@\gather}\else
 \ingather@true\inany@true\def\tag{&}%
 \vspace@\allowdisplaybreak@\displaybreak@\intertext@
 \displ@y\Let@
 \iftagsleft@\DN@{\csname gather \endcsname}\else
  \DN@{\csname gather \space\endcsname}\fi\fi
 \else\DN@{\onlydmatherr@\gather}\fi\next@}
\expandafter\def\csname gather \space\endcsname#1\endgather
 {\gmeasure@#1\endgather\tabskip\centering@
 \halign to\displaywidth{\hfil\strut@\setboxz@h{$\m@th\displaystyle{##}$}%
 \global\gwidth@\wdz@\boxz@\hfil&
 \setboxz@h{\strut@{\maketag@##\maketag@}}%
 \dimen@\displaywidth\advance\dimen@-\gwidth@
 \ifdim\dimen@>\tw@\wdz@\llap{\boxz@}\else
 \llap{\vtop{\normalbaselines\null\boxz@}}\fi
 \tabskip\z@skip\crcr#1\crcr\black@\gmaxwidth@}}
\newdimen\glineht@
\expandafter\def\csname gather \endcsname#1\endgather{\gmeasure@#1\endgather
 \ifdim\gmaxwidth@>\displaywidth\let\gdisplaywidth@\gmaxwidth@\else
 \let\gdisplaywidth@\displaywidth\fi\tabskip\centering@\halign to\displaywidth
 {\hfil\strut@\setboxz@h{$\m@th\displaystyle{##}$}%
 \global\gwidth@\wdz@\global\glineht@\ht\z@\boxz@\hfil&\kern-\gdisplaywidth@
 \setboxz@h{\strut@{\maketag@##\maketag@}}%
 \dimen@\displaywidth\advance\dimen@-\gwidth@
 \ifdim\dimen@>\tw@\wdz@\rlap{\boxz@}\else
 \rlap{\vbox{\normalbaselines\boxz@\vbox to\glineht@{}}}\fi
 \tabskip\gdisplaywidth@\crcr#1\crcr\black@\gmaxwidth@}}
\newif\ifctagsplit@
\def\CenteredTagsOnSplits{\global\ctagsplit@true}
\def\TopOrBottomTagsOnSplits{\global\ctagsplit@false}
\TopOrBottomTagsOnSplits
\def\split{\relax\ifinany@\let\next@\insplit@\else
 \ifmmode\ifinner\def\next@{\onlydmatherr@\split}\else
 \let\next@\outsplit@\fi\else
 \def\next@{\onlydmatherr@\split}\fi\fi\next@}
\def\insplit@{\global\setbox\z@\vbox\bgroup\vspace@\Let@\ialign\bgroup
 \hfil\strut@$\m@th\displaystyle{##}$&$\m@th\displaystyle{{}##}$\hfill\crcr}
\def\endsplit{\crcr\egroup\egroup\iftagsleft@\expandafter\lendsplit@\else
 \expandafter\rendsplit@\fi}
\def\rendsplit@{\global\setbox9 \vbox
 {\unvcopy\z@\global\setbox8 \lastbox\unskip}%                              %1
 \setbox@ne\hbox{\unhcopy8 \unskip\global\setbox\tw@\lastbox
 \unskip\global\setbox\thr@@\lastbox}%                                      %2
 \global\setbox7 \hbox{\unhbox\tw@\unskip}%                                 %3
 \ifinalign@\ifctagsplit@                                                   %4
  \gdef\split@{\hbox to\wd\thr@@{}&
   \vcenter{\vbox{\moveleft\wd\thr@@\boxz@}}}%                              %5
 \else\gdef\split@{&\vbox{\moveleft\wd\thr@@\box9}\crcr
  \box\thr@@&\box7}\fi                                                      %6
 \else                                                                      %7
  \ifctagsplit@\gdef\split@{\vcenter{\boxz@}}\else
  \gdef\split@{\box9\crcr\hbox{\box\thr@@\box7}}\fi
 \fi
 \split@}                                                                   %8
\def\lendsplit@{\global\setbox9\vtop{\unvcopy\z@}%                          %1
 \setbox@ne\vbox{\unvcopy\z@\global\setbox8\lastbox}%                       %2
 \setbox@ne\hbox{\unhcopy8\unskip\setbox\tw@\lastbox
  \unskip\global\setbox\thr@@\lastbox}%                                     %3
 \ifinalign@\ifctagsplit@                                                   %4
  \gdef\split@{\hbox to\wd\thr@@{}&
  \vcenter{\vbox{\moveleft\wd\thr@@\box9}}}%                                %5
  \else                                                                     %6
  \gdef\split@{\hbox to\wd\thr@@{}&\vbox{\moveleft\wd\thr@@\box9}}\fi
 \else
  \ifctagsplit@\gdef\split@{\vcenter{\box9}}\else
  \gdef\split@{\box9}\fi
 \fi\split@}
\def\outsplit@#1$${\align\insplit@#1\endalign$$}
\newdimen\multlinegap@
\multlinegap@1em
\newdimen\multlinetaggap@
\multlinetaggap@1em
\def\MultlineGap#1{\global\multlinegap@#1\relax}
\def\multlinegap#1{\RIfMIfI@\onlydmatherr@\multlinegap\else
 \multlinegap@#1\relax\fi\else\onlydmatherr@\multlinegap\fi}
\def\nomultlinegap{\multlinegap{\z@}}
\def\multline{\RIfMIfI@
 \DN@{\onlydmatherr@\multline}\else
 \DN@{\multline@}\fi\else
 \DN@{\onlydmatherr@\multline}\fi\next@}
\newif\iftagin@
\def\tagin@#1{\tagin@false\in@\tag{#1}\ifin@\tagin@true\fi}
\def\multline@#1$${\inany@true\vspace@\allowdisplaybreak@\displaybreak@
 \tagin@{#1}\iftagsleft@\DN@{\multline@l#1$$}\else
 \DN@{\multline@r#1$$}\fi\next@}
\newdimen\mwidth@
\def\rmmeasure@#1\endmultline{%
 \def\shoveleft##1{##1}\def\shoveright##1{##1}%                             %1
 \setbox@ne\vbox{\Let@\halign{\setboxz@h
  {$\m@th\@lign\displaystyle{}##$}\global\mwidth@\wdz@
  \crcr#1\crcr}}}
\newdimen\mlineht@
\newif\ifzerocr@
\newif\ifonecr@
\def\lmmeasure@#1\endmultline{\global\zerocr@true\global\onecr@false
 \everycr{\noalign{\ifonecr@\global\onecr@false\fi
  \ifzerocr@\global\zerocr@false\global\onecr@true\fi}}%                    %1
  \def\shoveleft##1{##1}\def\shoveright##1{##1}%
 \setbox@ne\vbox{\Let@\halign{\setboxz@h
  {$\m@th\@lign\displaystyle{}##$}\ifonecr@\global\mwidth@\wdz@
  \global\mlineht@\ht\z@\fi\crcr#1\crcr}}}
\newbox\mtagbox@
\newdimen\ltwidth@
\newdimen\rtwidth@
\def\multline@l#1$${\iftagin@\DN@{\lmultline@@#1$$}\else
 \DN@{\setbox\mtagbox@\null\ltwidth@\z@\rtwidth@\z@
  \lmultline@@@#1$$}\fi\next@}
\def\lmultline@@#1\endmultline\tag#2$${%
 \setbox\mtagbox@\hbox{\maketag@#2\maketag@}%                               %1
 \lmmeasure@#1\endmultline\dimen@\mwidth@\advance\dimen@\wd\mtagbox@
 \advance\dimen@\multlinetaggap@                                            %2
 \ifdim\dimen@>\displaywidth\ltwidth@\z@\else\ltwidth@\wd\mtagbox@\fi       %3
 \lmultline@@@#1\endmultline$$}
\def\lmultline@@@{\displ@y
 \def\shoveright##1{##1\hfilneg\hskip\multlinegap@}%
 \def\shoveleft##1{\setboxz@h{$\m@th\displaystyle{}##1$}%
  \setbox@ne\hbox{$\m@th\displaystyle##1$}%
  \hfilneg
  \iftagin@
   \ifdim\ltwidth@>\z@\hskip\ltwidth@\hskip\multlinetaggap@\fi
  \else\hskip\multlinegap@\fi\hskip.5\wd@ne\hskip-.5\wdz@##1}%              %1
  \halign\bgroup\Let@\hbox to\displaywidth
   {\strut@$\m@th\displaystyle\hfil{}##\hfil$}\crcr
   \hfilneg                                                                 %2
   \iftagin@                                                                %3
    \ifdim\ltwidth@>\z@                                                     %4
     \box\mtagbox@\hskip\multlinetaggap@                                    %5
    \else
     \rlap{\vbox{\normalbaselines\hbox{\strut@\box\mtagbox@}%
     \vbox to\mlineht@{}}}\fi                                               %6
   \else\hskip\multlinegap@\fi}                                             %7
\def\multline@r#1$${\iftagin@\DN@{\rmultline@@#1$$}\else
 \DN@{\setbox\mtagbox@\null\ltwidth@\z@\rtwidth@\z@
  \rmultline@@@#1$$}\fi\next@}
\def\rmultline@@#1\endmultline\tag#2$${\ltwidth@\z@
 \setbox\mtagbox@\hbox{\maketag@#2\maketag@}%
 \rmmeasure@#1\endmultline\dimen@\mwidth@\advance\dimen@\wd\mtagbox@
 \advance\dimen@\multlinetaggap@
 \ifdim\dimen@>\displaywidth\rtwidth@\z@\else\rtwidth@\wd\mtagbox@\fi
 \rmultline@@@#1\endmultline$$}
\def\rmultline@@@{\displ@y
 \def\shoveright##1{##1\hfilneg\iftagin@\ifdim\rtwidth@>\z@
  \hskip\rtwidth@\hskip\multlinetaggap@\fi\else\hskip\multlinegap@\fi}%
 \def\shoveleft##1{\setboxz@h{$\m@th\displaystyle{}##1$}%
  \setbox@ne\hbox{$\m@th\displaystyle##1$}%
  \hfilneg\hskip\multlinegap@\hskip.5\wd@ne\hskip-.5\wdz@##1}%
 \halign\bgroup\Let@\hbox to\displaywidth
  {\strut@$\m@th\displaystyle\hfil{}##\hfil$}\crcr
 \hfilneg\hskip\multlinegap@}
\def\endmultline{\iftagsleft@\expandafter\lendmultline@\else
 \expandafter\rendmultline@\fi}
\def\lendmultline@{\hfilneg\hskip\multlinegap@\crcr\egroup}
\def\rendmultline@{\iftagin@                                                %1
 \ifdim\rtwidth@>\z@                                                        %2
  \hskip\multlinetaggap@\box\mtagbox@                                       %3
 \else\llap{\vtop{\normalbaselines\null\hbox{\strut@\box\mtagbox@}}}\fi     %4
 \else\hskip\multlinegap@\fi                                                %5
 \hfilneg\crcr\egroup}
\def\bmod{\mskip-\medmuskip\mkern5mu\mathbin{\fam\z@ mod}\penalty900
 \mkern5mu\mskip-\medmuskip}
\def\pmod#1{\allowbreak\ifinner\mkern8mu\else\mkern18mu\fi
 ({\fam\z@ mod}\,\,#1)}
\def\pod#1{\allowbreak\ifinner\mkern8mu\else\mkern18mu\fi(#1)}
\def\mod#1{\allowbreak\ifinner\mkern12mu\else\mkern18mu\fi{\fam\z@ mod}\,\,#1}
\message{continued fractions,}
\newcount\cfraccount@
\def\cfrac{\bgroup\bgroup\advance\cfraccount@\@ne\strut
 \iffalse{\fi\def\\{\over\displaystyle}\iffalse}\fi}
\def\lcfrac{\bgroup\bgroup\advance\cfraccount@\@ne\strut
 \iffalse{\fi\def\\{\hfill\over\displaystyle}\iffalse}\fi}
\def\rcfrac{\bgroup\bgroup\advance\cfraccount@\@ne\strut\hfill
 \iffalse{\fi\def\\{\over\displaystyle}\iffalse}\fi}
\def\gloop@#1\repeat{\gdef\body{#1}\iterate}
\def\endcfrac{\gloop@\ifnum\cfraccount@>\z@\global\advance\cfraccount@\m@ne
 \egroup\hskip-\nulldelimiterspace\egroup\repeat}
\message{compound symbols,}
\def\binrel@#1{\setboxz@h{\thinmuskip0mu
  \medmuskip\m@ne mu\thickmuskip\@ne mu$#1\m@th$}%
 \setbox@ne\hbox{\thinmuskip0mu\medmuskip\m@ne mu\thickmuskip
  \@ne mu${}#1{}\m@th$}%
 \setbox\tw@\hbox{\hskip\wd@ne\hskip-\wdz@}}
\def\overset#1\to#2{\binrel@{#2}\ifdim\wd\tw@<\z@
 \mathbin{\mathop{\kern\z@#2}\limits^{#1}}\else\ifdim\wd\tw@>\z@
 \mathrel{\mathop{\kern\z@#2}\limits^{#1}}\else
 {\mathop{\kern\z@#2}\limits^{#1}}{}\fi\fi}
\def\underset#1\to#2{\binrel@{#2}\ifdim\wd\tw@<\z@
 \mathbin{\mathop{\kern\z@#2}\limits_{#1}}\else\ifdim\wd\tw@>\z@
 \mathrel{\mathop{\kern\z@#2}\limits_{#1}}\else
 {\mathop{\kern\z@#2}\limits_{#1}}{}\fi\fi}
\def\oversetbrace#1\to#2{\overbrace{#2}^{#1}}
\def\undersetbrace#1\to#2{\underbrace{#2}_{#1}}
\def\sideset#1\and#2\to#3{%
 \setbox@ne\hbox{$\dsize{\vphantom{#3}}#1{#3}\m@th$}%
 \setbox\tw@\hbox{$\dsize{#3}#2\m@th$}%
 \hskip\wd@ne\hskip-\wd\tw@\mathop{\hskip\wd\tw@\hskip-\wd@ne
  {\vphantom{#3}}#1{#3}#2}}
\def\rightarrowfill@#1{\setboxz@h{$#1-\m@th$}\ht\z@\z@
  $#1\m@th\copy\z@\mkern-6mu\cleaders
  \hbox{$#1\mkern-2mu\box\z@\mkern-2mu$}\hfill
  \mkern-6mu\mathord\rightarrow$}
\def\leftarrowfill@#1{\setboxz@h{$#1-\m@th$}\ht\z@\z@
  $#1\m@th\mathord\leftarrow\mkern-6mu\cleaders
  \hbox{$#1\mkern-2mu\copy\z@\mkern-2mu$}\hfill
  \mkern-6mu\box\z@$}
\def\leftrightarrowfill@#1{\setboxz@h{$#1-\m@th$}\ht\z@\z@
  $#1\m@th\mathord\leftarrow\mkern-6mu\cleaders
  \hbox{$#1\mkern-2mu\box\z@\mkern-2mu$}\hfill
  \mkern-6mu\mathord\rightarrow$}
\def\overrightarrow{\mathpalette\overrightarrow@}
\def\overrightarrow@#1#2{\vbox{\ialign{##\crcr\rightarrowfill@#1\crcr
 \noalign{\kern-\ex@\nointerlineskip}$\m@th\hfil#1#2\hfil$\crcr}}}

\def\overleftarrow{\mathpalette\overleftarrow@}
\def\overleftarrow@#1#2{\vbox{\ialign{##\crcr\leftarrowfill@#1\crcr
 \noalign{\kern-\ex@\nointerlineskip}$\m@th\hfil#1#2\hfil$\crcr}}}
\def\overleftrightarrow{\mathpalette\overleftrightarrow@}
\def\overleftrightarrow@#1#2{\vbox{\ialign{##\crcr\leftrightarrowfill@#1\crcr
 \noalign{\kern-\ex@\nointerlineskip}$\m@th\hfil#1#2\hfil$\crcr}}}
\def\underrightarrow{\mathpalette\underrightarrow@}
\def\underrightarrow@#1#2{\vtop{\ialign{##\crcr$\m@th\hfil#1#2\hfil$\crcr
 \noalign{\nointerlineskip}\rightarrowfill@#1\crcr}}}

\def\underleftarrow{\mathpalette\underleftarrow@}
\def\underleftarrow@#1#2{\vtop{\ialign{##\crcr$\m@th\hfil#1#2\hfil$\crcr
 \noalign{\nointerlineskip}\leftarrowfill@#1\crcr}}}
\def\underleftrightarrow{\mathpalette\underleftrightarrow@}
\def\underleftrightarrow@#1#2{\vtop{\ialign{##\crcr$\m@th\hfil#1#2\hfil$\crcr
 \noalign{\nointerlineskip}\leftrightarrowfill@#1\crcr}}}
\message{various kinds of dots,}
\let\DOTSI\relax
\let\DOTSB\relax

\newif\ifmath@
{\uccode`7=`\\ \uccode`8=`m \uccode`9=`a \uccode`0=`t \uccode`!=`h
 \uppercase{\gdef\math@#1#2#3#4#5#6\math@{\global\math@false\ifx 7#1\ifx 8#2%
 \ifx 9#3\ifx 0#4\ifx !#5\xdef\meaning@{#6}\global\math@true\fi\fi\fi\fi\fi}}}
\newif\ifmathch@
{\uccode`7=`c \uccode`8=`h \uccode`9=`\"
 \uppercase{\gdef\mathch@#1#2#3#4#5#6\mathch@{\global\mathch@false
  \ifx 7#1\ifx 8#2\ifx 9#5\global\mathch@true\xdef\meaning@{9#6}\fi\fi\fi}}}
\newcount\classnum@
\def\getmathch@#1.#2\getmathch@{\classnum@#1 \divide\classnum@4096
 \ifcase\number\classnum@\or\or\gdef\thedots@{\dotsb@}\or
 \gdef\thedots@{\dotsb@}\fi}
\newif\ifmathbin@
{\uccode`4=`b \uccode`5=`i \uccode`6=`n
 \uppercase{\gdef\mathbin@#1#2#3{\relaxnext@
  \DNii@##1\mathbin@{\ifx\space@\next\global\mathbin@true\fi}%
 \global\mathbin@false\DN@##1\mathbin@{}%
 \ifx 4#1\ifx 5#2\ifx 6#3\DN@{\FN@\nextii@}\fi\fi\fi\next@}}}
\newif\ifmathrel@
{\uccode`4=`r \uccode`5=`e \uccode`6=`l
 \uppercase{\gdef\mathrel@#1#2#3{\relaxnext@
  \DNii@##1\mathrel@{\ifx\space@\next\global\mathrel@true\fi}%
 \global\mathrel@false\DN@##1\mathrel@{}%
 \ifx 4#1\ifx 5#2\ifx 6#3\DN@{\FN@\nextii@}\fi\fi\fi\next@}}}
\newif\ifmacro@
{\uccode`5=`m \uccode`6=`a \uccode`7=`c
 \uppercase{\gdef\macro@#1#2#3#4\macro@{\global\macro@false
  \ifx 5#1\ifx 6#2\ifx 7#3\global\macro@true
  \xdef\meaning@{\macro@@#4\macro@@}\fi\fi\fi}}}
\def\macro@@#1->#2\macro@@{#2}
\newif\ifDOTS@
\newcount\DOTSCASE@
{\uccode`6=`\\ \uccode`7=`D \uccode`8=`O \uccode`9=`T \uccode`0=`S
 \uppercase{\gdef\DOTS@#1#2#3#4#5{\global\DOTS@false\DN@##1\DOTS@{}%
  \ifx 6#1\ifx 7#2\ifx 8#3\ifx 9#4\ifx 0#5\let\next@\DOTS@@\fi\fi\fi\fi\fi
  \next@}}}
{\uccode`3=`B \uccode`4=`I \uccode`5=`X
 \uppercase{\gdef\DOTS@@#1{\relaxnext@
  \DNii@##1\DOTS@{\ifx\space@\next\global\DOTS@true\fi}%
  \DN@{\FN@\nextii@}%
  \ifx 3#1\global\DOTSCASE@\z@\else
  \ifx 4#1\global\DOTSCASE@\@ne\else
  \ifx 5#1\global\DOTSCASE@\tw@\else\DN@##1\DOTS@{}%
  \fi\fi\fi\next@}}}
\newif\ifnot@
{\uccode`5=`\\ \uccode`6=`n \uccode`7=`o \uccode`8=`t
 \uppercase{\gdef\not@#1#2#3#4{\relaxnext@
  \DNii@##1\not@{\ifx\space@\next\global\not@true\fi}%
 \global\not@false\DN@##1\not@{}%
 \ifx 5#1\ifx 6#2\ifx 7#3\ifx 8#4\DN@{\FN@\nextii@}\fi\fi\fi
 \fi\next@}}}
\newif\ifkeybin@
\def\keybin@{\keybin@true
 \ifx\next+\else\ifx\next=\else\ifx\next<\else\ifx\next>\else\ifx\next-\else
 \ifx\next*\else\ifx\next:\else\keybin@false\fi\fi\fi\fi\fi\fi\fi}
\def\dots{\RIfM@\expandafter\mdots@\else\expandafter\tdots@\fi}
\def\tdots@{\unskip\relaxnext@
 \DN@{$\m@th\mathinner{\ldotp\ldotp\ldotp}\,
   \ifx\next,\,$\else\ifx\next.\,$\else\ifx\next;\,$\else\ifx\next:\,$\else
   \ifx\next?\,$\else\ifx\next!\,$\else$ \fi\fi\fi\fi\fi\fi}%
 \ \FN@\next@}
\def\mdots@{\FN@\mdots@@}
\def\mdots@@{\gdef\thedots@{\dotso@}%                                       %1
 \ifx\next\boldkey\gdef\thedots@\boldkey{\boldkeydots@}\else                %2
 \ifx\next\boldsymbol\gdef\thedots@\boldsymbol{\boldsymboldots@}\else       %3
 \ifx,\next\gdef\thedots@{\dotsc}%                                          %4
 \else\ifx\not\next\gdef\thedots@{\dotsb@}%                                 %5
 \else\keybin@
 \ifkeybin@\gdef\thedots@{\dotsb@}%                                         %6
 \else\xdef\meaning@{\meaning\next..........}\xdef\meaning@@{\meaning@}%    %7
  \expandafter\math@\meaning@\math@
  \ifmath@
   \expandafter\mathch@\meaning@\mathch@
   \ifmathch@\expandafter\getmathch@\meaning@\getmathch@\fi                 %8
  \else\expandafter\macro@\meaning@@\macro@                                 %9
  \ifmacro@                                                                %10
   \expandafter\not@\meaning@\not@\ifnot@\gdef\thedots@{\dotsb@}%          %11
  \else\expandafter\DOTS@\meaning@\DOTS@
  \ifDOTS@
   \ifcase\number\DOTSCASE@\gdef\thedots@{\dotsb@}%
    \or\gdef\thedots@{\dotsi}\else\fi                                      %12
  \else\expandafter\math@\meaning@\math@                                   %13
  \ifmath@\expandafter\mathbin@\meaning@\mathbin@
  \ifmathbin@\gdef\thedots@{\dotsb@}%                                      %14
  \else\expandafter\mathrel@\meaning@\mathrel@
  \ifmathrel@\gdef\thedots@{\dotsb@}%                                      %15
  \fi\fi\fi\fi\fi\fi\fi\fi\fi\fi\fi\fi
 \thedots@}
\def\plainldots@{\mathinner{\ldotp\ldotp\ldotp}}
\def\plaincdots@{\mathinner{\cdotp\cdotp\cdotp}}
\def\dotsi{\!\plaincdots@}
\let\dotsb@\plaincdots@
\newif\ifextra@
\newif\ifrightdelim@
\def\rightdelim@{\global\rightdelim@true                                    %1
 \ifx\next)\else                                                            %2
 \ifx\next]\else
 \ifx\next\rbrack\else
 \ifx\next\}\else
 \ifx\next\rbrace\else
 \ifx\next\rangle\else
 \ifx\next\rceil\else
 \ifx\next\rfloor\else
 \ifx\next\rgroup\else
 \ifx\next\rmoustache\else
 \ifx\next\right\else
 \ifx\next\bigr\else
 \ifx\next\biggr\else
 \ifx\next\Bigr\else                                                        %3
 \ifx\next\Biggr\else\global\rightdelim@false
 \fi\fi\fi\fi\fi\fi\fi\fi\fi\fi\fi\fi\fi\fi\fi}
\def\extra@{%
 \global\extra@false\rightdelim@\ifrightdelim@\global\extra@true            %1
 \else\ifx\next$\global\extra@true                                          %2
 \else\xdef\meaning@{\meaning\next..........}%                              %3
 \expandafter\macro@\meaning@\macro@\ifmacro@                               %4
 \expandafter\DOTS@\meaning@\DOTS@
 \ifDOTS@
 \ifnum\DOTSCASE@=\tw@\global\extra@true                                    %5
 \fi\fi\fi\fi\fi}
\newif\ifbold@
\def\dotso@{\relaxnext@
 \ifbold@
  \let\next\delayed@
  \DNii@{\extra@\plainldots@\ifextra@\,\fi}%
 \else
  \DNii@{\DN@{\extra@\plainldots@\ifextra@\,\fi}\FN@\next@}%
 \fi
 \nextii@}
\def\extrap@#1{%
 \ifx\next,\DN@{#1\,}\else
 \ifx\next;\DN@{#1\,}\else
 \ifx\next.\DN@{#1\,}\else\extra@
 \ifextra@\DN@{#1\,}\else
 \let\next@#1\fi\fi\fi\fi\next@}
\def\ldots{\DN@{\extrap@\plainldots@}%
 \FN@\next@}
\def\cdots{\DN@{\extrap@\plaincdots@}%
 \FN@\next@}

\def\dotsc{\relaxnext@
 \DN@{\ifx\next;\plainldots@\,\else
  \ifx\next.\plainldots@\,\else\extra@\plainldots@
  \ifextra@\,\fi\fi\fi}%
 \FN@\next@}
\def\cdot{\mathchar"2201 }

\message{special superscripts,}
\def\dddot#1{{\mathop{#1}\limits^{\vbox to-1.4\ex@{\kern-\tw@\ex@
 \hbox{\rm...}\vss}}}}
\def\ddddot#1{{\mathop{#1}\limits^{\vbox to-1.4\ex@{\kern-\tw@\ex@
 \hbox{\rm....}\vss}}}}
\def\sphat{^{\mathchoice{}{}%
 {\,\,\botsmash{\hbox{\lower4\ex@\hbox{$\m@th\widehat{\null}$}}}}%
 {\,\botsmash{\hbox{\lower3\ex@\hbox{$\m@th\hat{\null}$}}}}}}

\def\spacute{^{\!\botsmash{\hbox{\lower\@ne ex\hbox{\'{}}}}}}
\def\spgrave{^{\mathchoice{}{}{}{\!}%
 \botsmash{\hbox{\lower\@ne ex\hbox{\`{}}}}}}
\def\spdot{^{\hbox{\raise\ex@\hbox{\rm.}}}}
\def\spddot{^{\hbox{\raise\ex@\hbox{\rm..}}}}
\def\spdddot{^{\hbox{\raise\ex@\hbox{\rm...}}}}
\def\spddddot{^{\hbox{\raise\ex@\hbox{\rm....}}}}
\def\spbreve{^{\!\botsmash{\hbox{\lower4\ex@\hbox{\u{}}}}}}

\message{\string\text,}
\def\textonlyfont@#1#2{\def#1{\RIfM@
 \Err@{Use \string#1\space only in text}\else#2\fi}}
\textonlyfont@\rm\tenrm
\textonlyfont@\it\tenit
\textonlyfont@\sl\tensl
\textonlyfont@\bf\tenbf
\def\oldnos#1{\RIfM@{\mathcode`\,="013B \fam\@ne#1}\else
 \leavevmode\hbox{$\m@th\mathcode`\,="013B \fam\@ne#1$}\fi}
\def\text{\RIfM@\expandafter\text@\else\expandafter\text@@\fi}
\def\text@@#1{\leavevmode\hbox{#1}}
\def\mathhexbox@#1#2#3{\text{$\m@th\mathchar"#1#2#3$}}
\def\dag{{\mathhexbox@279}}
\def\ddag{{\mathhexbox@27A}}
\def\S{{\mathhexbox@278}}
\def\P{{\mathhexbox@27B}}
\newif\iffirstchoice@
\firstchoice@true
\def\text@#1{\mathchoice
 {\hbox{\everymath{\displaystyle}\def\textfonti{\the\textfont\@ne}%
  \def\textfontii{\the\textfont\tw@}\textdef@@ T#1}}
 {\hbox{\firstchoice@false
  \everymath{\textstyle}\def\textfonti{\the\textfont\@ne}%
  \def\textfontii{\the\textfont\tw@}\textdef@@ T#1}}
 {\hbox{\firstchoice@false
  \everymath{\scriptstyle}\def\textfonti{\the\scriptfont\@ne}%
  \def\textfontii{\the\scriptfont\tw@}\textdef@@ S\rm#1}}
 {\hbox{\firstchoice@false
  \everymath{\scriptscriptstyle}\def\textfonti
  {\the\scriptscriptfont\@ne}%
  \def\textfontii{\the\scriptscriptfont\tw@}\textdef@@ s\rm#1}}}
\def\textdef@@#1{\textdef@#1\rm\textdef@#1\bf\textdef@#1\sl\textdef@#1\it}
\def\rmfam{0}
\def\textdef@#1#2{%
 \DN@{\csname\expandafter\eat@\string#2fam\endcsname}%
 \if S#1\edef#2{\the\scriptfont\next@\relax}%
 \else\if s#1\edef#2{\the\scriptscriptfont\next@\relax}%
 \else\edef#2{\the\textfont\next@\relax}\fi\fi}
\scriptfont\itfam\tenit \scriptscriptfont\itfam\tenit
\scriptfont\slfam\tensl \scriptscriptfont\slfam\tensl
\newif\iftopfolded@
\newif\ifbotfolded@
\def\topfoldedtext{\topfolded@true\botfolded@false\foldedtext@}
\def\botfoldedtext{\botfolded@true\topfolded@false\foldedtext@}
\def\foldedtext{\topfolded@false\botfolded@false\foldedtext@}
\Invalid@\foldedwidth
\def\foldedtext@{\relaxnext@
 \DN@{\ifx\next\foldedwidth\let\next@\nextii@\else
  \DN@{\nextii@\foldedwidth{.3\hsize}}\fi\next@}%
 \DNii@\foldedwidth##1##2{\setbox\z@\vbox
  {\normalbaselines\hsize##1\relax
  \tolerance1600 \noindent\ignorespaces##2}\ifbotfolded@\boxz@\else
  \iftopfolded@\vtop{\unvbox\z@}\else\vcenter{\boxz@}\fi\fi}%
 \FN@\next@}
\message{math font commands,}
\def\bold{\RIfM@\expandafter\bold@\else
 \expandafter\nonmatherr@\expandafter\bold\fi}
\def\bold@#1{{\bold@@{#1}}}
\def\bold@@#1{\fam\bffam\relax#1}
\def\slanted{\RIfM@\expandafter\slanted@\else
 \expandafter\nonmatherr@\expandafter\slanted\fi}
\def\slanted@#1{{\slanted@@{#1}}}
\def\slanted@@#1{\fam\slfam\relax#1}
\def\roman{\RIfM@\expandafter\roman@\else
 \expandafter\nonmatherr@\expandafter\roman\fi}
\def\roman@#1{{\roman@@{#1}}}
\def\roman@@#1{\fam\rmfam\relax#1}
\def\italic{\RIfM@\expandafter\italic@\else
 \expandafter\nonmatherr@\expandafter\italic\fi}
\def\italic@#1{{\italic@@{#1}}}
\def\italic@@#1{\fam\itfam\relax#1}
\def\Cal{\RIfM@\expandafter\Cal@\else
 \expandafter\nonmatherr@\expandafter\Cal\fi}
\def\Cal@#1{{\Cal@@{#1}}}
\def\Cal@@#1{\noaccents@\fam\tw@#1}
\mathchardef\Gamma="0000
\mathchardef\Delta="0001
\mathchardef\Theta="0002
\mathchardef\Lambda="0003
\mathchardef\Xi="0004
\mathchardef\Pi="0005
\mathchardef\Sigma="0006
\mathchardef\Upsilon="0007
\mathchardef\Phi="0008
\mathchardef\Psi="0009
\mathchardef\Omega="000A
\mathchardef\varGamma="0100
\mathchardef\varDelta="0101
\mathchardef\varTheta="0102
\mathchardef\varLambda="0103
\mathchardef\varXi="0104
\mathchardef\varPi="0105
\mathchardef\varSigma="0106
\mathchardef\varUpsilon="0107
\mathchardef\varPhi="0108
\mathchardef\varPsi="0109
\mathchardef\varOmega="010A
\let\alloc@@\alloc@
\def\hexnumber@#1{\ifcase#1 0\or 1\or 2\or 3\or 4\or 5\or 6\or 7\or 8\or
 9\or A\or B\or C\or D\or E\or F\fi}
\def\loadmsam{%
 \font@\tenmsa=msam10
 \font@\sevenmsa=msam7
 \font@\fivemsa=msam5
 \alloc@@8\fam\chardef\sixt@@n\msafam
 \textfont\msafam=\tenmsa
 \scriptfont\msafam=\sevenmsa
 \scriptscriptfont\msafam=\fivemsa
 \edef\next{\hexnumber@\msafam}%
 \mathchardef\dabar@"0\next39
 \edef\dashrightarrow{\mathrel{\dabar@\dabar@\mathchar"0\next4B}}%
 \edef\dashleftarrow{\mathrel{\mathchar"0\next4C\dabar@\dabar@}}%
 \let\dasharrow\dashrightarrow
 \edef\ulcorner{\delimiter"4\next70\next70 }%
 \edef\urcorner{\delimiter"5\next71\next71 }%
 \edef\llcorner{\delimiter"4\next78\next78 }%
 \edef\lrcorner{\delimiter"5\next79\next79 }%
 \edef\yen{{\noexpand\mathhexbox@\next55}}%
 \edef\checkmark{{\noexpand\mathhexbox@\next58}}%
 \edef\circledR{{\noexpand\mathhexbox@\next72}}%
 \edef\maltese{{\noexpand\mathhexbox@\next7A}}%
 \global\let\loadmsam\empty}%
\def\loadmsbm{%
 \font@\tenmsb=msbm10 \font@\sevenmsb=msbm7 \font@\fivemsb=msbm5
 \alloc@@8\fam\chardef\sixt@@n\msbfam
 \textfont\msbfam=\tenmsb
 \scriptfont\msbfam=\sevenmsb \scriptscriptfont\msbfam=\fivemsb
 \global\let\loadmsbm\empty
 }
\def\widehat#1{\ifx\undefined\msbfam \DN@{362}%
  \else \setboxz@h{$\m@th#1$}%
    \edef\next@{\ifdim\wdz@>\tw@ em%
        \hexnumber@\msbfam 5B%
      \else 362\fi}\fi
  \mathaccent"0\next@{#1}}
\def\widetilde#1{\ifx\undefined\msbfam \DN@{365}%
  \else \setboxz@h{$\m@th#1$}%
    \edef\next@{\ifdim\wdz@>\tw@ em%
        \hexnumber@\msbfam 5D%
      \else 365\fi}\fi
  \mathaccent"0\next@{#1}}
\message{\string\newsymbol,}
\def\newsymbol#1#2#3#4#5{\define#1{}%
  \count@#2\relax \advance\count@\m@ne % to push case 0 to the \else clause
 \ifcase\count@
   \ifx\undefined\msafam\loadmsam\fi \let\next@\msafam
 \or \ifx\undefined\msbfam\loadmsbm\fi \let\next@\msbfam
 \else  \Err@{\Invalid@@\string\newsymbol}\let\next@\tw@\fi
 \mathchardef#1="#3\hexnumber@\next@#4#5\space}

\def\Bbb{\RIfM@\expandafter\Bbb@\else
 \expandafter\nonmatherr@\expandafter\Bbb\fi}
\def\Bbb@#1{{\Bbb@@{#1}}}
\def\Bbb@@#1{\noaccents@\fam\msbfam\relax#1}
\message{bold Greek and bold symbols,}
\def\loadbold{%
 \font@\tencmmib=cmmib10 \font@\sevencmmib=cmmib7 \font@\fivecmmib=cmmib5
 \skewchar\tencmmib'177 \skewchar\sevencmmib'177 \skewchar\fivecmmib'177
 \alloc@@8\fam\chardef\sixt@@n\cmmibfam
 \textfont\cmmibfam\tencmmib
 \scriptfont\cmmibfam\sevencmmib \scriptscriptfont\cmmibfam\fivecmmib
 \font@\tencmbsy=cmbsy10 \font@\sevencmbsy=cmbsy7 \font@\fivecmbsy=cmbsy5
 \skewchar\tencmbsy'60 \skewchar\sevencmbsy'60 \skewchar\fivecmbsy'60
 \alloc@@8\fam\chardef\sixt@@n\cmbsyfam
 \textfont\cmbsyfam\tencmbsy
 \scriptfont\cmbsyfam\sevencmbsy \scriptscriptfont\cmbsyfam\fivecmbsy
 \let\loadbold\empty
}
\def\boldnotloaded#1{\Err@{\ifcase#1\or First\else Second\fi
       bold symbol font not loaded}}
\def\mathchari@#1#2#3{\ifx\undefined\cmmibfam
    \boldnotloaded@\@ne
  \else\mathchar"#1\hexnumber@\cmmibfam#2#3\space \fi}
\def\mathcharii@#1#2#3{\ifx\undefined\cmbsyfam
    \boldnotloaded\tw@
  \else \mathchar"#1\hexnumber@\cmbsyfam#2#3\space\fi}
\edef\bffam@{\hexnumber@\bffam}
\def\boldkey#1{\ifcat\noexpand#1A%
  \ifx\undefined\cmmibfam \boldnotloaded\@ne
  \else {\fam\cmmibfam#1}\fi
 \else
 \ifx#1!\mathchar"5\bffam@21 \else
 \ifx#1(\mathchar"4\bffam@28 \else\ifx#1)\mathchar"5\bffam@29 \else
 \ifx#1+\mathchar"2\bffam@2B \else\ifx#1:\mathchar"3\bffam@3A \else
 \ifx#1;\mathchar"6\bffam@3B \else\ifx#1=\mathchar"3\bffam@3D \else
 \ifx#1?\mathchar"5\bffam@3F \else\ifx#1[\mathchar"4\bffam@5B \else
 \ifx#1]\mathchar"5\bffam@5D \else
 \ifx#1,\mathchari@63B \else
 \ifx#1-\mathcharii@200 \else
 \ifx#1.\mathchari@03A \else
 \ifx#1/\mathchari@03D \else
 \ifx#1<\mathchari@33C \else
 \ifx#1>\mathchari@33E \else
 \ifx#1*\mathcharii@203 \else
 \ifx#1|\mathcharii@06A \else
 \ifx#10\bold0\else\ifx#11\bold1\else\ifx#12\bold2\else\ifx#13\bold3\else
 \ifx#14\bold4\else\ifx#15\bold5\else\ifx#16\bold6\else\ifx#17\bold7\else
 \ifx#18\bold8\else\ifx#19\bold9\else
  \Err@{\string\boldkey\space can't be used with #1}%
 \fi\fi\fi\fi\fi\fi\fi\fi\fi\fi\fi\fi\fi\fi\fi
 \fi\fi\fi\fi\fi\fi\fi\fi\fi\fi\fi\fi\fi\fi}
\def\boldsymbol#1{%
 \DN@{\Err@{You can't use \string\boldsymbol\space with \string#1}#1}%
 \ifcat\noexpand#1A%
   \let\next@\relax
   \ifx\undefined\cmmibfam \boldnotloaded\@ne
   \else {\fam\cmmibfam#1}\fi
 \else
  \xdef\meaning@{\meaning#1.........}%
  \expandafter\math@\meaning@\math@
  \ifmath@
   \expandafter\mathch@\meaning@\mathch@
   \ifmathch@
    \expandafter\boldsymbol@@\meaning@\boldsymbol@@
   \fi
  \else
   \expandafter\macro@\meaning@\macro@
   \expandafter\delim@\meaning@\delim@
   \ifdelim@
    \expandafter\delim@@\meaning@\delim@@
   \else
    \boldsymbol@{#1}%
   \fi
  \fi
 \fi
 \next@}
\def\mathhexboxii@#1#2{\ifx\undefined\cmbsyfam
    \boldnotloaded\tw@
  \else \mathhexbox@{\hexnumber@\cmbsyfam}{#1}{#2}\fi}
\def\boldsymbol@#1{\let\next@\relax\let\next#1%
 \ifx\next\cdot\mathcharii@201 \else
 \ifx\next\prime{{\null\mathcharii@030 \null}}\else
 \ifx\next\lbrack\mathchar"4\bffam@5B \else
 \ifx\next\rbrack\mathchar"5\bffam@5D \else
 \ifx\next\{\mathcharii@466 \else
 \ifx\next\lbrace\mathcharii@466 \else
 \ifx\next\}\mathcharii@567 \else
 \ifx\next\rbrace\mathcharii@567 \else
 \ifx\next\surd{{\mathcharii@170}}\else
 \ifx\next\S{{\mathhexboxii@78}}\else
 \ifx\next\P{{\mathhexboxii@7B}}\else
 \ifx\next\dag{{\mathhexboxii@79}}\else
 \ifx\next\ddag{{\mathhexboxii@7A}}\else
 \DN@{\Err@{You can't use \string\boldsymbol\space with \string#1}#1}%
 \fi\fi\fi\fi\fi\fi\fi\fi\fi\fi\fi\fi\fi}
\def\boldsymbol@@#1.#2\boldsymbol@@{\classnum@#1 \count@@@\classnum@        %1
 \divide\classnum@4096 \count@\classnum@                                    %2
 \multiply\count@4096 \advance\count@@@-\count@ \count@@\count@@@           %3
 \divide\count@@@\@cclvi \count@\count@@                                    %4
 \multiply\count@@@\@cclvi \advance\count@@-\count@@@                       %5
 \divide\count@@@\@cclvi                                                    %6
 \multiply\classnum@4096 \advance\classnum@\count@@                         %7
 \ifnum\count@@@=\z@                                                        %8
  \count@"\bffam@ \multiply\count@\@cclvi
  \advance\classnum@\count@
  \DN@{\mathchar\number\classnum@}%
 \else
  \ifnum\count@@@=\@ne                                                      %9
   \ifx\undefined\cmmibfam \DN@{\boldnotloaded\@ne}%
   \else \count@\cmmibfam \multiply\count@\@cclvi
     \advance\classnum@\count@
     \DN@{\mathchar\number\classnum@}\fi
  \else
   \ifnum\count@@@=\tw@                                                    %10
     \ifx\undefined\cmbsyfam
       \DN@{\boldnotloaded\tw@}%
     \else
       \count@\cmbsyfam \multiply\count@\@cclvi
       \advance\classnum@\count@
       \DN@{\mathchar\number\classnum@}%
     \fi
  \fi
 \fi
\fi}
\newif\ifdelim@
\newcount\delimcount@
{\uccode`6=`\\ \uccode`7=`d \uccode`8=`e \uccode`9=`l
 \uppercase{\gdef\delim@#1#2#3#4#5\delim@
  {\delim@false\ifx 6#1\ifx 7#2\ifx 8#3\ifx 9#4\delim@true
   \xdef\meaning@{#5}\fi\fi\fi\fi}}}
\def\delim@@#1"#2#3#4#5#6\delim@@{\if#32%
\let\next@\relax
 \ifx\undefined\cmbsyfam \boldnotloaded\@ne
 \else \mathcharii@#2#4#5\space \fi\fi}
\def\vert{\delimiter"026A30C }
\def\Vert{\delimiter"026B30D }
\let\|\Vert
\def\backslash{\delimiter"026E30F }
\def\boldkeydots@#1{\bold@true\let\next=#1\let\delayed@=#1\mdots@@
 \boldkey#1\bold@false}  % = required!
\def\boldsymboldots@#1{\bold@true\let\next#1\let\delayed@#1\mdots@@
 \boldsymbol#1\bold@false}
\message{Euler fonts,}

\def\frak{\mathfont@\frak}

\def\loadmathfont#1{%
   \expandafter\font@\csname ten#1\endcsname=#110
   \expandafter\font@\csname seven#1\endcsname=#17
   \expandafter\font@\csname five#1\endcsname=#15
   \edef\next{\noexpand\alloc@@8\fam\chardef\sixt@@n
     \expandafter\noexpand\csname#1fam\endcsname}%
   \next
   \textfont\csname#1fam\endcsname \csname ten#1\endcsname
   \scriptfont\csname#1fam\endcsname \csname seven#1\endcsname
   \scriptscriptfont\csname#1fam\endcsname \csname five#1\endcsname
   \expandafter\def\csname #1\expandafter\endcsname\expandafter{%
      \expandafter\mathfont@\csname#1\endcsname}%
 \expandafter\gdef\csname load#1\endcsname{}%
}
\def\mathfont@#1{\RIfM@\expandafter\mathfont@@\expandafter#1\else
  \expandafter\nonmatherr@\expandafter#1\fi}
\def\mathfont@@#1#2{{\mathfont@@@#1{#2}}}
\def\mathfont@@@#1#2{\noaccents@
   \fam\csname\expandafter\eat@\string#1fam\endcsname
   \relax#2}
\message{math accents,}
\def\accentclass@{7}
\def\noaccents@{\def\accentclass@{0}}
\def\makeacc@#1#2{\def#1{\mathaccent"\accentclass@#2 }}
\makeacc@\hat{05E}
\makeacc@\check{014}
\makeacc@\tilde{07E}
\makeacc@\acute{013}
\makeacc@\grave{012}
\makeacc@\dot{05F}
\makeacc@\ddot{07F}
\makeacc@\breve{015}
\makeacc@\bar{016}

\newcount\skewcharcount@
\newcount\familycount@
\def\theskewchar@{\familycount@\@ne
 \global\skewcharcount@\the\skewchar\textfont\@ne                           %1
 \ifnum\fam>\m@ne\ifnum\fam<16
  \global\familycount@\the\fam\relax
  \global\skewcharcount@\the\skewchar\textfont\the\fam\relax\fi\fi          %2
 \ifnum\skewcharcount@>\m@ne
  \ifnum\skewcharcount@<128
  \multiply\familycount@256
  \global\advance\skewcharcount@\familycount@
  \global\advance\skewcharcount@28672
  \mathchar\skewcharcount@\else
  \global\skewcharcount@\m@ne\fi\else
 \global\skewcharcount@\m@ne\fi}                                            %3
\newcount\pointcount@
\def\getpoints@#1.#2\getpoints@{\pointcount@#1 }
\newdimen\accentdimen@
\newcount\accentmu@
\def\dimentomu@{\multiply\accentdimen@ 100
 \expandafter\getpoints@\the\accentdimen@\getpoints@
 \multiply\pointcount@18
 \divide\pointcount@\@m
 \global\accentmu@\pointcount@}
\def\Makeacc@#1#2{\def#1{\RIfM@\DN@{\mathaccent@
 {"\accentclass@#2 }}\else\DN@{\nonmatherr@{#1}}\fi\next@}}
\def\unbracefonts@{\let\Cal@\Cal@@\let\roman@\roman@@\let\bold@\bold@@
 \let\slanted@\slanted@@}
\def\mathaccent@#1#2{\ifnum\fam=\m@ne\xdef\thefam@{1}\else
 \xdef\thefam@{\the\fam}\fi                                                 %1
 \accentdimen@\z@                                                           %2
 \setboxz@h{\unbracefonts@$\m@th\fam\thefam@\relax#2$}%                     %3
 \ifdim\accentdimen@=\z@\DN@{\mathaccent#1{#2}}%                            %4
  \setbox@ne\hbox{\unbracefonts@$\m@th\fam\thefam@\relax#2\theskewchar@$}% %5a
  \setbox\tw@\hbox{$\m@th\ifnum\skewcharcount@=\m@ne\else
   \mathchar\skewcharcount@\fi$}%                                          %5b
  \global\accentdimen@\wd@ne\global\advance\accentdimen@-\wdz@
  \global\advance\accentdimen@-\wd\tw@                                     %5c
  \global\multiply\accentdimen@\tw@
  \dimentomu@\global\advance\accentmu@\@ne                                 %5d
 \else\DN@{{\mathaccent#1{#2\mkern\accentmu@ mu}%
    \mkern-\accentmu@ mu}{}}\fi                                             %6
 \next@}\Makeacc@\Hat{05E}
\Makeacc@\Check{014}
\Makeacc@\Tilde{07E}
\Makeacc@\Acute{013}
\Makeacc@\Grave{012}
\Makeacc@\Dot{05F}
\Makeacc@\Ddot{07F}
\Makeacc@\Breve{015}
\Makeacc@\Bar{016}
\def\Vec{\RIfM@\DN@{\mathaccent@{"017E }}\else
 \DN@{\nonmatherr@\Vec}\fi\next@}
\def\accentedsymbol#1#2{\csname newbox\expandafter\endcsname
  \csname\expandafter\eat@\string#1@box\endcsname
 \expandafter\setbox\csname\expandafter\eat@
  \string#1@box\endcsname\hbox{$\m@th#2$}\define
  #1{\copy\csname\expandafter\eat@\string#1@box\endcsname{}}}
\message{roots,}
\def\sqrt#1{\radical"270370 {#1}}
\let\underline@\underline
\let\overline@\overline
\def\underline#1{\underline@{#1}}
\def\overline#1{\overline@{#1}}
\Invalid@\leftroot
\Invalid@\uproot
\newcount\uproot@
\newcount\leftroot@
\def\root{\relaxnext@
  \DN@{\ifx\next\uproot\let\next@\nextii@\else
   \ifx\next\leftroot\let\next@\nextiii@\else
   \let\next@\plainroot@\fi\fi\next@}%
  \DNii@\uproot##1{\uproot@##1\relax\FN@\nextiv@}%
  \def\nextiv@{\ifx\next\space@\DN@. {\FN@\nextv@}\else
   \DN@.{\FN@\nextv@}\fi\next@.}%
  \def\nextv@{\ifx\next\leftroot\let\next@\nextvi@\else
   \let\next@\plainroot@\fi\next@}%
  \def\nextvi@\leftroot##1{\leftroot@##1\relax\plainroot@}%
   \def\nextiii@\leftroot##1{\leftroot@##1\relax\FN@\nextvii@}%
  \def\nextvii@{\ifx\next\space@
   \DN@. {\FN@\nextviii@}\else
   \DN@.{\FN@\nextviii@}\fi\next@.}%
  \def\nextviii@{\ifx\next\uproot\let\next@\nextix@\else
   \let\next@\plainroot@\fi\next@}%
  \def\nextix@\uproot##1{\uproot@##1\relax\plainroot@}%
  \bgroup\uproot@\z@\leftroot@\z@\FN@\next@}
\def\plainroot@#1\of#2{\setbox\rootbox\hbox{$\m@th\scriptscriptstyle{#1}$}%
 \mathchoice{\r@@t\displaystyle{#2}}{\r@@t\textstyle{#2}}
 {\r@@t\scriptstyle{#2}}{\r@@t\scriptscriptstyle{#2}}\egroup}
\def\r@@t#1#2{\setboxz@h{$\m@th#1\sqrt{#2}$}%
 \dimen@\ht\z@\advance\dimen@-\dp\z@
 \setbox@ne\hbox{$\m@th#1\mskip\uproot@ mu$}\advance\dimen@ 1.667\wd@ne
 \mkern-\leftroot@ mu\mkern5mu\raise.6\dimen@\copy\rootbox
 \mkern-10mu\mkern\leftroot@ mu\boxz@}
\def\boxed#1{\setboxz@h{$\m@th\displaystyle{#1}$}\dimen@.4\ex@
 \advance\dimen@3\ex@\advance\dimen@\dp\z@
 \hbox{\lower\dimen@\hbox{%
 \vbox{\hrule height.4\ex@
 \hbox{\vrule width.4\ex@\hskip3\ex@\vbox{\vskip3\ex@\boxz@\vskip3\ex@}%
 \hskip3\ex@\vrule width.4\ex@}\hrule height.4\ex@}%
 }}}
\message{commutative diagrams,}
\let\ampersand@\relax
\newdimen\minaw@
\minaw@11.11128\ex@
\newdimen\minCDaw@
\minCDaw@2.5pc
\def\minCDarrowwidth#1{\RIfMIfI@\onlydmatherr@\minCDarrowwidth
 \else\minCDaw@#1\relax\fi\else\onlydmatherr@\minCDarrowwidth\fi}
\newif\ifCD@
\def\CD{\bgroup\vspace@\relax\iffalse{\fi\let\ampersand@&\iffalse}\fi
 \CD@true\vcenter\bgroup\Let@\tabskip\z@skip\baselineskip20\ex@
 \lineskip3\ex@\lineskiplimit3\ex@\halign\bgroup
 &\hfill$\m@th##$\hfill\crcr}
\def\endCD{\crcr\egroup\egroup\egroup}
\newdimen\bigaw@
\atdef@>#1>#2>{\ampersand@                                                  %1
 \setboxz@h{$\m@th\ssize\;{#1}\;\;$}%                                       %2
 \setbox@ne\hbox{$\m@th\ssize\;{#2}\;\;$}%                                  %3
 \setbox\tw@\hbox{$\m@th#2$}%                                               %4
 \ifCD@\global\bigaw@\minCDaw@\else\global\bigaw@\minaw@\fi                 %5
 \ifdim\wdz@>\bigaw@\global\bigaw@\wdz@\fi
 \ifdim\wd@ne>\bigaw@\global\bigaw@\wd@ne\fi                                %6
 \ifCD@\enskip\fi                                                           %7
 \ifdim\wd\tw@>\z@
  \mathrel{\mathop{\hbox to\bigaw@{\rightarrowfill@\displaystyle}}%
    \limits^{#1}_{#2}}%                                                     %8
 \else\mathrel{\mathop{\hbox to\bigaw@{\rightarrowfill@\displaystyle}}%
    \limits^{#1}}\fi                                                        %9
 \ifCD@\enskip\fi                                                          %10
 \ampersand@}                                                              %11
\atdef@<#1<#2<{\ampersand@\setboxz@h{$\m@th\ssize\;\;{#1}\;$}%
 \setbox@ne\hbox{$\m@th\ssize\;\;{#2}\;$}\setbox\tw@\hbox{$\m@th#2$}%
 \ifCD@\global\bigaw@\minCDaw@\else\global\bigaw@\minaw@\fi
 \ifdim\wdz@>\bigaw@\global\bigaw@\wdz@\fi
 \ifdim\wd@ne>\bigaw@\global\bigaw@\wd@ne\fi
 \ifCD@\enskip\fi
 \ifdim\wd\tw@>\z@
  \mathrel{\mathop{\hbox to\bigaw@{\leftarrowfill@\displaystyle}}%
       \limits^{#1}_{#2}}\else
  \mathrel{\mathop{\hbox to\bigaw@{\leftarrowfill@\displaystyle}}%
       \limits^{#1}}\fi
 \ifCD@\enskip\fi\ampersand@}
\begingroup
 \catcode`\~=\active \lccode`\~=`\@
 \lowercase{%
  \global\atdef@)#1)#2){~>#1>#2>}
  \global\atdef@(#1(#2({~<#1<#2<}}
\endgroup
\atdef@ A#1A#2A{\llap{$\m@th\vcenter{\hbox
 {$\ssize#1$}}$}\Big\uparrow\rlap{$\m@th\vcenter{\hbox{$\ssize#2$}}$}&&}
\atdef@ V#1V#2V{\llap{$\m@th\vcenter{\hbox
 {$\ssize#1$}}$}\Big\downarrow\rlap{$\m@th\vcenter{\hbox{$\ssize#2$}}$}&&}
\atdef@={&\enskip\mathrel
 {\vbox{\hrule width\minCDaw@\vskip3\ex@\hrule width
 \minCDaw@}}\enskip&}
\atdef@|{\Big\Vert&&}
\atdef@\vert{\Big\Vert&&}
\def\pretend#1\haswidth#2{\setboxz@h{$\m@th\scriptstyle{#2}$}\hbox
 to\wdz@{\hfill$\m@th\scriptstyle{#1}$\hfill}}
\message{poor man's bold,}
\def\pmb{\RIfM@\expandafter\mathpalette\expandafter\pmb@\else
 \expandafter\pmb@@\fi}
\def\pmb@@#1{\leavevmode\setboxz@h{#1}%
   \dimen@-\wdz@
   \kern-.5\ex@\copy\z@
   \kern\dimen@\kern.25\ex@\raise.4\ex@\copy\z@
   \kern\dimen@\kern.25\ex@\box\z@
}
\def\binrel@@#1{\ifdim\wd2<\z@\mathbin{#1}\else\ifdim\wd\tw@>\z@
 \mathrel{#1}\else{#1}\fi\fi}
\newdimen\pmbraise@
%      Note: because of the use of \mathpalette, if \pmb@ is
%      applied to a single math italic character (or a single
%      character from some other slanted math font), the italic
%      correction will be added.  This will cause subscripts
%      to fall too far away from the character in some
%      cases, e.g., $\pmb{T}_1$ or $\pmb{\Cal T}_1$.
\def\pmb@#1#2{\setbox\thr@@\hbox{$\m@th#1{#2}$}%
 \setbox4\hbox{$\m@th#1\mkern.5mu$}\pmbraise@\wd4\relax
 \binrel@{#2}%
 \dimen@-\wd\thr@@
   \binrel@@{%
   \mkern-.8mu\copy\thr@@
   \kern\dimen@\mkern.4mu\raise\pmbraise@\copy\thr@@
   \kern\dimen@\mkern.4mu\box\thr@@
}}
\def\documentstyle#1{\W@{}\input #1.sty\relax}
\message{syntax check,}
\font\dummyft@=dummy
\fontdimen1 \dummyft@=\z@
\fontdimen2 \dummyft@=\z@
\fontdimen3 \dummyft@=\z@
\fontdimen4 \dummyft@=\z@
\fontdimen5 \dummyft@=\z@
\fontdimen6 \dummyft@=\z@
\fontdimen7 \dummyft@=\z@
\fontdimen8 \dummyft@=\z@
\fontdimen9 \dummyft@=\z@
\fontdimen10 \dummyft@=\z@
\fontdimen11 \dummyft@=\z@
\fontdimen12 \dummyft@=\z@
\fontdimen13 \dummyft@=\z@
\fontdimen14 \dummyft@=\z@
\fontdimen15 \dummyft@=\z@
\fontdimen16 \dummyft@=\z@
\fontdimen17 \dummyft@=\z@
\fontdimen18 \dummyft@=\z@
\fontdimen19 \dummyft@=\z@
\fontdimen20 \dummyft@=\z@
\fontdimen21 \dummyft@=\z@
\fontdimen22 \dummyft@=\z@
\def\fontlist@{\\{\tenrm}\\{\sevenrm}\\{\fiverm}\\{\teni}\\{\seveni}%
 \\{\fivei}\\{\tensy}\\{\sevensy}\\{\fivesy}\\{\tenex}\\{\tenbf}\\{\sevenbf}%
 \\{\fivebf}\\{\tensl}\\{\tenit}}
\def\font@#1=#2 {\rightappend@#1\to\fontlist@\font#1=#2 }
\def\dodummy@{{\def\\##1{\global\let##1\dummyft@}\fontlist@}}
\def\nopages@{\output{\setbox\z@\box\@cclv \deadcycles\z@}%
 \alloc@5\toks\toksdef\@cclvi\output}
\let\galleys\nopages@
\newif\ifsyntax@
\newcount\countxviii@
\def\syntax{\syntax@true\dodummy@\countxviii@\count18
 \loop\ifnum\countxviii@>\m@ne\textfont\countxviii@=\dummyft@
 \scriptfont\countxviii@=\dummyft@\scriptscriptfont\countxviii@=\dummyft@
 \advance\countxviii@\m@ne\repeat                                           %1
 \dummyft@\tracinglostchars\z@\nopages@\frenchspacing\hbadness\@M}
\def\first@#1#2\end{#1}
\def\printoptions{\W@{Do you want S(yntax check),
  G(alleys) or P(ages)?}%
 \message{Type S, G or P, followed by <return>: }%
 \begingroup % to localize the following change to \endlinechar:
 \endlinechar\m@ne % to prevent a space or \par in \ans@ from ^^M
 \read\m@ne to\ans@
%  Define \ans@ to uppercase itself, and default to P if the user
%  just pressed <return>.
 \edef\ans@{\uppercase{\def\noexpand\ans@{%
   \expandafter\first@\ans@ P\end}}}%
%  Cast the new definition of \ans@ outside the group
 \expandafter\endgroup\ans@
 \if\ans@ P% fine, no action needs to be taken
 \else \if\ans@ S\syntax
 \else \if\ans@ G\galleys
 \else\message{? Unknown option: \ans@; using the `pages' option.}%
 \fi\fi\fi}
\def\alloc@#1#2#3#4#5{\global\advance\count1#1by\@ne
 \ch@ck#1#4#2\allocationnumber=\count1#1
 \global#3#5=\allocationnumber
 \ifalloc@\wlog{\string#5=\string#2\the\allocationnumber}\fi}
\def\document{\def\alloclist@{}\def\fontlist@{}}
\let\enddocument\bye

\let\proclaim\undefined
\let\footnote\undefined
\let\=\undefined
\let\>\undefined

\catcode`\@=\active
\message{... finished}

%  *** end including amstex.tex *** 
% % \input mathdefs
%  *** start including mathdefs.tex *** 
\expandafter\ifx\csname mathdefs.tex\endcsname\relax
  \expandafter\gdef\csname mathdefs.tex\endcsname{}
\else \message{Hey!  Apparently you were trying to
  \string\input{mathdefs.tex} twice.   This does not make sense.} 
\errmessage{Please edit your file (probably \jobname.tex) and remove
any duplicate ``\string\input'' lines}\endinput\fi

%mathdefs.tex v1.3.2

%%% Changes from v1.0: footnote macros, warning for duplicated tags,
%%%   control sequences \( and \verbatimtags.
%%% From v1.2: \pretags, redefinition of \( using \ifinner, multi-part
%%%   equation numbering, control sequences \[, \references, and
%%%   \resetbracket. 
%%% From v1.3: \rm in \lastpart; write root of multi-part tag to .tgs 

%See file texdefs.doc for documentation.

\catcode`\X=12\catcode`\@=11

%Minor control sequences:
\def\n@wcount{\alloc@0\count\countdef\insc@unt}
\def\n@wwrite{\alloc@7\write\chardef\sixt@@n}
\def\n@wread{\alloc@6\read\chardef\sixt@@n}
\def\r@s@t{\relax}\def\v@idline{\par}\def\@mputate#1/{#1}
\def\l@c@l#1X{\firstpart.#1}\def\gl@b@l#1X{#1}\def\t@d@l#1X{{}}

%Creation of tag families and output of assignments and citations:
\def\crossrefs#1{\ifx\all#1\let\tr@ce=\all\else\def\tr@ce{#1,}\fi
   \n@wwrite\cit@tionsout\openout\cit@tionsout=\jobname.cit 
   \write\cit@tionsout{\tr@ce}\expandafter\setfl@gs\tr@ce,}
\def\setfl@gs#1,{\def\@{#1}\ifx\@\empty\let\next=\relax
   \else\let\next=\setfl@gs\expandafter\xdef
   \csname#1tr@cetrue\endcsname{}\fi\next}
\def\m@ketag#1#2{\expandafter\n@wcount\csname#2tagno\endcsname
     \csname#2tagno\endcsname=0\let\tail=\all\xdef\all{\tail#2,}
   \ifx#1\l@c@l\let\tail=\r@s@t\xdef\r@s@t{\csname#2tagno\endcsname=0\tail}\fi
   \expandafter\gdef\csname#2cite\endcsname##1{\expandafter
     \ifx\csname#2tag##1\endcsname\relax?\else\csname#2tag##1\endcsname\fi
     \expandafter\ifx\csname#2tr@cetrue\endcsname\relax\else
     \write\cit@tionsout{#2tag ##1 cited on page \folio.}\fi}
   \expandafter\gdef\csname#2page\endcsname##1{\expandafter
     \ifx\csname#2page##1\endcsname\relax?\else\csname#2page##1\endcsname\fi
     \expandafter\ifx\csname#2tr@cetrue\endcsname\relax\else
     \write\cit@tionsout{#2tag ##1 cited on page \folio.}\fi}
   \expandafter\gdef\csname#2tag\endcsname##1{\expandafter
      \ifx\csname#2check##1\endcsname\relax
      \expandafter\xdef\csname#2check##1\endcsname{}%
      \else\immediate\write16{Warning: #2tag ##1 used more than once.}\fi
      \multit@g{#1}{#2}##1/X%
      \write\t@gsout{#2tag ##1 assigned number \csname#2tag##1\endcsname\space
      on page \number\count0.}%
   \csname#2tag##1\endcsname}}

\def\multit@g#1#2#3/#4X{\def\t@mp{#4}\ifx\t@mp\empty%
      \global\advance\csname#2tagno\endcsname by 1 
      \expandafter\xdef\csname#2tag#3\endcsname
      {#1\number\csname#2tagno\endcsnameX}%
   \else\expandafter\ifx\csname#2last#3\endcsname\relax
      \expandafter\n@wcount\csname#2last#3\endcsname
      \global\advance\csname#2tagno\endcsname by 1 
      \expandafter\xdef\csname#2tag#3\endcsname
      {#1\number\csname#2tagno\endcsnameX}
      \write\t@gsout{#2tag #3 assigned number \csname#2tag#3\endcsname\space
      on page \number\count0.}\fi
   \global\advance\csname#2last#3\endcsname by 1
   \def\t@mp{\expandafter\xdef\csname#2tag#3/}%
   \expandafter\t@mp\@mputate#4\endcsname
   {\csname#2tag#3\endcsname\lastpart{\csname#2last#3\endcsname}}\fi}
\def\t@gs#1{\def\all{}\m@ketag#1e\m@ketag#1s\m@ketag\t@d@l p
\let\realscite\scite
\let\realstag\stag
   \m@ketag\gl@b@l r \n@wread\t@gsin
   \openin\t@gsin=\jobname.tgs \re@der \closein\t@gsin
   \n@wwrite\t@gsout\openout\t@gsout=\jobname.tgs }
\outer\def\localtags{\t@gs\l@c@l}
\outer\def\globaltags{\t@gs\gl@b@l}
\outer\def\newlocaltag#1{\m@ketag\l@c@l{#1}}
\outer\def\newglobaltag#1{\m@ketag\gl@b@l{#1}}

%Reading in tag information:
\newif\ifpr@ 
\def\m@kecs #1tag #2 assigned number #3 on page #4.%
   {\expandafter\gdef\csname#1tag#2\endcsname{#3}
   \expandafter\gdef\csname#1page#2\endcsname{#4}
   \ifpr@\expandafter\xdef\csname#1check#2\endcsname{}\fi}
\def\re@der{\ifeof\t@gsin\let\next=\relax\else
   \read\t@gsin to\t@gline\ifx\t@gline\v@idline\else
   \expandafter\m@kecs \t@gline\fi\let \next=\re@der\fi\next}
\def\pretags#1{\pr@true\pret@gs#1,,}
\def\pret@gs#1,{\def\@{#1}\ifx\@\empty\let\n@xtfile=\relax
   \else\let\n@xtfile=\pret@gs \openin\t@gsin=#1.tgs \message{#1} \re@der 
   \closein\t@gsin\fi \n@xtfile}

%Sections and subsections; local numbering:
\newcount\sectno\sectno=0\newcount\subsectno\subsectno=0
\newif\ifultr@local \def\ultralocal{\ultr@localtrue}
\def\firstpart{\number\sectno}
\def\lastpart#1{\ifcase#1 \or a\or b\or c\or d\or e\or f\or g\or h\or 
   i\or k\or l\or m\or n\or o\or p\or q\or r\or s\or t\or u\or v\or w\or 
   x\or y\or z \fi}

\def\resetall{\global\advance\sectno by 1\subsectno=0
   \gdef\firstpart{\number\sectno}\r@s@t}
\def\resetsub{\global\advance\subsectno by 1
   \gdef\firstpart{\number\sectno.\number\subsectno}\r@s@t}
\def\newsection#1\par{\resetall\vskip0pt plus.3\vsize\penalty-250
   \vskip0pt plus-.3\vsize\bigskip\bigskip
   \message{#1}\leftline{\bf#1}\nobreak\bigskip}
\def\subsection#1\par{\ifultr@local\resetsub\fi
   \vskip0pt plus.2\vsize\penalty-250\vskip0pt plus-.2\vsize
   \bigskip\smallskip\message{#1}\leftline{\bf#1}\nobreak\medskip}

%jj tags:
% On Andrzej's request:  we want to be able 
% to show tags as in noverbatim, with verbatim in the margin,
% and cites as in verbatim, with nonverbatim in the margin
% mg -- July 2000

\newdimen\marginshift

\newdimen\margindelta
\newdimen\marginmax
\newdimen\marginmin

\def\margininit{       
\marginmax=3 true cm                  % how much room, approximately
				      
\margindelta=0.1 true cm              % distance between entries
\marginmin=0.1true cm                 % where will leftmost entry be
\marginshift=\marginmin
}    % we cannot execute this right now, since 
     % there may be a \magnification coming later in the 
     % main file.   So we call \margininit at the end of 
     % alice2jlem

\def\t@gsjj#1,{\def\@{#1}\ifx\@\empty\let\next=\relax\else\let\next=\t@gsjj
   \def\@@{p}\ifx\@\@@\else
   \expandafter\gdef\csname#1cite\endcsname##1{\citejj{##1}}
   \expandafter\gdef\csname#1page\endcsname##1{?}
   \expandafter\gdef\csname#1tag\endcsname##1{\tagjj{##1}}\fi\fi\next}
\newif\ifshowstuffinmargin
\showstuffinmarginfalse
\def\jjtags{\ifx\shlhetal\relax 
      % so this is a public version --> no-op 
  \else
\ifx\shlhetal\undefinedcontrolseq
      % again, this is a public version --> no-op 
\else
\showstuffinmargintrue
\ifx\all\relax\else\expandafter\t@gsjj\all,\fi\fi \fi
}

\def\tagjj#1{\realstag{#1}\oldmginpar{\zeigen{#1}}}
\def\citejj#1{\rechnen{#1}\mginpar{\zeigen{#1}}}     % modified Sep 02, saharon's suggestion

\def\rechnen#1{\expandafter\ifx\csname stag#1\endcsname\relax ??\else
                           \csname stag#1\endcsname\fi}

\newdimen\theight

\def\marginfont{\sevenrm}

\def\trymarginbox#1{\setbox0=\hbox{\marginfont\hskip\marginshift #1}%
		\global\marginshift\wd0 
		\global\advance\marginshift\margindelta}

\def \oldmginpar#1{%
\ifvmode\setbox0\hbox to \hsize{\hfill\rlap{\marginfont\quad#1}}%
\ht0 0cm
\dp0 0cm
\box0\vskip-\baselineskip
\else 
             \vadjust{\trymarginbox{#1}%
		\ifdim\marginshift>\marginmax \global\marginshift\marginmin
			\trymarginbox{#1}%
                \fi
             \theight=\ht0
             \advance\theight by \dp0    \advance\theight by \lineskip
             \kern -\theight \vbox to \theight{\rightline{\rlap{\box0}}%
\vss}}\fi}

\newdimen\upordown
\global\upordown=8pt
\font\tinyfont=cmtt8 % scaled 700
\def\mginpar#1{\smash{\hbox to 0cm{\kern-10pt\raise7pt\hbox{\tinyfont #1}\hss}}}
% testing, october 2005, mg
\def\mginpar#1{{\hbox to 0cm{\kern-10pt\raise\upordown\hbox{\tinyfont #1}\hss}}\global\upordown-\upordown}

% \def\mginpar#1{mg-#1-mg }

%Verbatim tags:
\def\t@gsoff#1,{\def\@{#1}\ifx\@\empty\let\next=\relax\else\let\next=\t@gsoff
   \def\@@{p}\ifx\@\@@\else
   \expandafter\gdef\csname#1cite\endcsname##1{\zeigen{##1}}
   \expandafter\gdef\csname#1page\endcsname##1{?}
   \expandafter\gdef\csname#1tag\endcsname##1{\zeigen{##1}}\fi\fi\next}
\def\verbatimtags{\showstuffinmarginfalse
\ifx\all\relax\else\expandafter\t@gsoff\all,\fi}
\def\zeigen#1{\hbox{$\scriptstyle\langle$}#1\hbox{$\scriptstyle\rangle$}}

% % \def\margincite#1{\ifshowstuffinmargin\mginpar{\rechnen{#1}}\fi}
%  changed, april 2003, mg: we now have always the verbatim tag in the margin
\def\margincite#1{\ifshowstuffinmargin\mginpar{\zeigen{#1}}\fi}

\def\margintag#1{\ifshowstuffinmargin\oldmginpar{\zeigen{#1}}\fi}

%Equation numbering:
\def\(#1){\edef\dot@g{\ifmmode\ifinner(\hbox{\noexpand\etag{#1}})
   \else\noexpand\eqno(\hbox{\noexpand\etag{#1}})\fi
   \else(\noexpand\ecite{#1})\fi}\dot@g}

%Reference numbering:
\newif\ifbr@ck
\def\eat#1{}
\def\[#1]{\br@cktrue[\br@cket#1'X]}
\def\br@cket#1'#2X{\def\temp{#2}\ifx\temp\empty\let\next\eat
   \else\let\next\br@cket\fi
   \ifbr@ck\br@ckfalse\br@ck@t#1,X\else\br@cktrue#1\fi\next#2X}
\def\br@ck@t#1,#2X{\def\temp{#2}\ifx\temp\empty\let\neext\eat
   \else\let\neext\br@ck@t\def\temp{,}\fi
   \def\teemp{#1}\ifx\teemp\empty\else\rcite{#1}\fi\temp\neext#2X}
\def\resetbr@cket{\gdef\[##1]{[\rtag{##1}]}}
\def\references{\resetbr@cket\newsection References\par}

%Footnotes:
\newtoks\symb@ls\newtoks\s@mb@ls\newtoks\p@gelist\n@wcount\ftn@mber
    \ftn@mber=1\newif\ifftn@mbers\ftn@mbersfalse\newif\ifbyp@ge\byp@gefalse
\def\defm@rk{\ifftn@mbers\n@mberm@rk\else\symb@lm@rk\fi}
\def\n@mberm@rk{\xdef\m@rk{{\the\ftn@mber}}%
    \global\advance\ftn@mber by 1 }
\def\rot@te#1{\let\temp=#1\global#1=\expandafter\r@t@te\the\temp,X}
\def\r@t@te#1,#2X{{#2#1}\xdef\m@rk{{#1}}}
\def\b@@st#1{{$^{#1}$}}\def\str@p#1{#1}
\def\symb@lm@rk{\ifbyp@ge\rot@te\p@gelist\ifnum\expandafter\str@p\m@rk=1 
    \s@mb@ls=\symb@ls\fi\write\f@nsout{\number\count0}\fi \rot@te\s@mb@ls}
\def\byp@ge{\byp@getrue\n@wwrite\f@nsin\openin\f@nsin=\jobname.fns 
    \n@wcount\currentp@ge\currentp@ge=0\p@gelist={0}
    \re@dfns\closein\f@nsin\rot@te\p@gelist
    \n@wread\f@nsout\openout\f@nsout=\jobname.fns }
\def\m@kelist#1X#2{{#1,#2}}
\def\re@dfns{\ifeof\f@nsin\let\next=\relax\else\read\f@nsin to \f@nline
    \ifx\f@nline\v@idline\else\let\t@mplist=\p@gelist
    \ifnum\currentp@ge=\f@nline
    \global\p@gelist=\expandafter\m@kelist\the\t@mplistX0
    \else\currentp@ge=\f@nline
    \global\p@gelist=\expandafter\m@kelist\the\t@mplistX1\fi\fi
    \let\next=\re@dfns\fi\next}
\def\symbols#1{\symb@ls={#1}\s@mb@ls=\symb@ls} 
\def\bigsymbol{\textstyle}
\symbols{\bigsymbol\ast,\dagger,\ddagger,\sharp,\flat,\natural,\star}
\def\ftnumbers{\ftn@mberstrue} \def\ftsymbols{\ftn@mbersfalse}
\def\paginal{\byp@ge} \def\resetftnumbers{\ftn@mber=1}
\def\ftnote#1{\defm@rk\expandafter\expandafter\expandafter\footnote
    \expandafter\b@@st\m@rk{#1}}

%Miscellaneous macros:
\long\def\jump#1\endjump{}
\def\ssum{\mathop{\lower .1em\hbox{$\textstyle\Sigma$}}\nolimits}

\def\qed{\nobreak\kern 1em \vrule height .5em width .5em depth 0em}
\def\newneq{\hbox{\rlap{\hbox to 1\wd9{\hss$=$\hss}}\raise .1em 
   \hbox to 1\wd9{\hss$\scriptscriptstyle/$\hss}}}
\def\subsetne{\setbox9 = \hbox{$\subset$}\mathrel{\hbox{\rlap
   {\lower .4em \newneq}\raise .13em \hbox{$\subset$}}}}
\def\supsetne{\setbox9 = \hbox{$\subset$}\mathrel{\hbox{\rlap
   {\lower .4em \newneq}\raise .13em \hbox{$\supset$}}}}

%Blackboard bold:
\def\vbar{\mathchoice{\vrule height6.3ptdepth-.5ptwidth.8pt\kern-.8pt}
   {\vrule height6.3ptdepth-.5ptwidth.8pt\kern-.8pt}
   {\vrule height4.1ptdepth-.35ptwidth.6pt\kern-.6pt}
   {\vrule height3.1ptdepth-.25ptwidth.5pt\kern-.5pt}}
\def\f@dge{\mathchoice{}{}{\mkern.5mu}{\mkern.8mu}}
\def\b@c#1#2{{\rm \mkern#2mu\vbar\mkern-#2mu#1}}
\def\b@b#1{{\rm I\mkern-3.5mu #1}}
\def\b@a#1#2{{\rm #1\mkern-#2mu\f@dge #1}}
\def\bb#1{{\count4=`#1 \advance\count4by-64 \ifcase\count4\or\b@a A{11.5}\or
   \b@b B\or\b@c C{5}\or\b@b D\or\b@b E\or\b@b F \or\b@c G{5}\or\b@b H\or
   \b@b I\or\b@c J{3}\or\b@b K\or\b@b L \or\b@b M\or\b@b N\or\b@c O{5} \or
   \b@b P\or\b@c Q{5}\or\b@b R\or\b@a S{8}\or\b@a T{10.5}\or\b@c U{5}\or
   \b@a V{12}\or\b@a W{16.5}\or\b@a X{11}\or\b@a Y{11.7}\or\b@a Z{7.5}\fi}}

\catcode`\X=11 \catcode`\@=12

% Sep 2003, mg:

% % \input 300stuff
%  *** start including 300stuff.tex *** 

% definitions needed to process the 300x papers  (plus 88 etc) 

\let\thischap\jobname

\def\partof#1{\csname returnthe#1part\endcsname}
\def\chapof#1{\csname returnthe#1chap\endcsname}

\def\setchapter#1,#2,#3;{% 
  \expandafter\def\csname returnthe#1part\endcsname{#2}%
  \expandafter\def\csname returnthe#1chap\endcsname{#3}%
}

\setchapter 300a,A,II.A;
\setchapter 300b,A,II.B;
\setchapter 300c,A,II.C;
\setchapter 300d,A,II.D;
\setchapter 300e,A,II.E;
\setchapter 300f,A,II.F;
\setchapter 300g,A,II.G;
\setchapter  E53,B,N;
\setchapter  88r,B,I;
\setchapter  600,B,III;
\setchapter  705,B,IV;
\setchapter  734,B,V;

\def\cprefix#1{%     \cprefix{300b} generates II.  or "nothing"
 %\leavevmode\vrule width 2cm height 1cm depth 0cm
\edef\theotherpart{\partof{#1}}\edef\theotherchap{\chapof{#1}}%
\ifx\theotherpart\thispart
   \ifx\theotherchap\thischap % nothing
    \else % same part, different chap
     \theotherchap%
    \fi
   \else % different part
     \theotherchap\fi}

 % \sectioncite[\S4]{300b}  ->  II.\S4   or A.II.\S4
\def\sectioncite[#1]#2{%
     \cprefix{#2}#1}

% now define \thispart
\edef\thispart{\partof{\thischap}}
\edef\thischap{\chapof{\thischap}}

\def\lastpage of '#1' is #2.{\expandafter\def\csname lastpage#1\endcsname{#2}}

%  *** end including 300stuff.tex *** 

\def\spuriousreset{}

%  *** end including mathdefs.tex *** 
% % \input citeadd
%  *** start including citeadd.tex *** 
%   citeadd -- a few additions for 
% files from alice that were procesed with "citealice"

\expandafter\ifx\csname citeadd.tex\endcsname\relax
\expandafter\gdef\csname citeadd.tex\endcsname{}
\else \message{Hey!  Apparently you were trying to
\string\input{citeadd.tex} twice.   This does not make sense.} 
\errmessage{Please edit your file (probably \jobname.tex) and remove
any duplicate ``\string\input'' lines}\endinput\fi

%  *** end including citeadd.tex *** 
\sectno=-1   % start with sect 0
\localtags
\jjtags
\NoBlackBoxes
\define\mr{\medskip\roster}
\define\sn{\smallskip\noindent}
\define\mn{\medskip\noindent}
\define\bn{\bigskip\noindent}
\define\ub{\underbar}
\define\wilog{\text{without loss of generality}}
\define\ermn{\endroster\medskip\noindent}

\define\dbcu{\dsize\bigcup}
\define \nl{\newline}
\magnification=\magstep 1
\documentstyle{amsppt}
% % \input alice2000
%  *** start including alice2000.tex *** 
% This file should be inputted whenever we use amsppt.sty and 
% old tex.  
%  Here we redefine \subjclass (use 1991 instead of 2000, otherwise 
% the following definition comes directly from 
%% 
%%              `amsppt.sty', generated 
%% on <1997/2/2> with the docstrip utility (2.2i).
%% 
%% The original source files were:
%% 
%% amsppt.doc 
%%% ====================================================================
%%% @AMSTeX-style-file{
%%%   filename  = "amsppt.sty",
%%%   version   = "2.1h",
%%%   date      = "1997/02/02",
%%%   time      = "09:27:44 EST",
%%%   checksum  = "56844 3264 16617 137829",
%%%   author    = "American Mathematical Society",
%%%   address   = "PO Box 6248, Providence, RI 02940-6248, USA",
%%%   telephone = "401-455-4080 or (in the USA) 800-321-4AMS",

{    % the braces make the catcode-change local. 
\catcode`@11

\ifx\alicetwothousandloaded@\relax
  \endinput\else\global\let\alicetwothousandloaded@\relax\fi

\gdef\subjclass{\let\savedef@\subjclass
 \def\subjclass##1\endsubjclass{\let\subjclass\savedef@
   \toks@{\def\usualspace{{\rm\enspace}}\eightpoint}%
   \toks@@{##1\unskip.}%
   \edef\thesubjclass@{\the\toks@
     \frills@{{\noexpand\rm2000 {\noexpand\it Mathematics Subject
       Classification}.\noexpand\enspace}}%
     \the\toks@@}}%
  \nofrillscheck\subjclass}
} 

%  *** end including alice2000.tex *** 
% % \input alice2jlem
%  *** start including alice2jlem.tex *** 
%% # Keywords  Input file to be used for texing Alice's files

\expandafter\ifx\csname alice2jlem.tex\endcsname\relax
  \expandafter\xdef\csname alice2jlem.tex\endcsname{\the\catcode`@}
\else \message{Hey!  Apparently you were trying to
\string\input{alice2jlem.tex}  twice.   This does not make sense.}
\errmessage{Please edit your file (probably \jobname.tex) and remove
any duplicate ``\string\input'' lines}\endinput\fi

% % \input bib4plain
%  *** start including bib4plain.tex *** 
\expandafter\ifx\csname bib4plain.tex\endcsname\relax
  \expandafter\gdef\csname bib4plain.tex\endcsname{}
\else \message{Hey!  Apparently you were trying to \string\input
  bib4plain.tex twice.   This does not make sense.}
\errmessage{Please edit your file (probably \jobname.tex) and remove
any duplicate ``\string\input'' lines}\endinput\fi

%  This file should be inputted if you want to use 
%  bibtex fom within plain TeX. 
      % Not really need for standard
       % bibtex files, but these commands
\def\renewcommand{\newcommand}	       % are used in our literal-unsrt.bst
\edef\cite{\the\catcode`@}%
\catcode`@ = 11
\let\@oldatcatcode = \cite
\chardef\@letter = 11
\chardef\@other = 12
%
%
% Next come some things that will be useful later.
%
% Make an outer definition into an inner one (due to Chris Thompson).
% The arguments should be the control sequence to be defined, and the
% new of the \outer control sequence, as characters; the control
% sequence #1 is defined to be just the same as \csname#2\endcsname, but
% not \outer.  For example, \@innerdef\innernewcount{newcount} would
% define \innernewcount to be a non-outer version of \newcount.
%
\def\@innerdef#1#2{\edef#1{\expandafter\noexpand\csname #2\endcsname}}%
%
% We use \@innerdef to make some of our allocations local, because
% Eplain includes our code inside a conditional.  We put @'s in the
% names to minimize the (already small) chance of conflicts.
%
\@innerdef\@innernewcount{newcount}%
\@innerdef\@innernewdimen{newdimen}%
\@innerdef\@innernewif{newif}%
\@innerdef\@innernewwrite{newwrite}%
%
%
% Swallow one parameter.
%
\def\@gobble#1{}%
%
%
% Use TeX 3.0's \inputlineno to get the line number, for better error
% messages, but if we're using an old version of TeX, don't do anything.
%
\ifx\inputlineno\@undefined
   \let\@linenumber = \empty % Pre-3.0.
\else
   \def\@linenumber{\the\inputlineno:\space}%
\fi
%
%
% The following macro \@futurenonspacelet (from the TeXbook) behaves
% essentially like \futurelet except that it discards any implicit or
% explicit space tokens that intervene before a nonspace is scanned:
%
\def\@futurenonspacelet#1{\def\cs{#1}%
   \afterassignment\@stepone\let\@nexttoken=
}%
\begingroup % The grouping here avoids stepping on an outside use of `\\'.
\def\\{\global\let\@stoken= }%
\\ % now \@stoken is a space token (\\ is a control symbol, so that
   % space after it is seen).
\endgroup
\def\@stepone{\expandafter\futurelet\cs\@steptwo}%
\def\@steptwo{\expandafter\ifx\cs\@stoken\let\@@next=\@stepthree
   \else\let\@@next=\@nexttoken\fi \@@next}%
\def\@stepthree{\afterassignment\@stepone\let\@@next= }%
%
%
% \@getoptionalarg\CS gets an optional argument from the input, enclosed
% in brackets, then expands \CS.  We set \@optionalarg to \empty if we
% don't find one, otherwise to the text of the argument.  This assumes
% the brackets don't have a funny category code.
%
\def\@getoptionalarg#1{%
   \let\@optionaltemp = #1%
   \let\@optionalnext = \relax
   \@futurenonspacelet\@optionalnext\@bracketcheck
}%
%
% The \expandafter's in this macro let us avoid the use of \aftergroup,
% which is somewhat more expensive.
%
\def\@bracketcheck{%
   \ifx [\@optionalnext
      \expandafter\@@getoptionalarg
   \else
      \let\@optionalarg = \empty
      % We can't do the \temp after the \fi, because then the \temp gets
      % in the way of reading the optional argument from the input, if
      % we do have one.
      \expandafter\@optionaltemp
   \fi
}%
\def\@@getoptionalarg[#1]{%
   \def\@optionalarg{#1}%
   \@optionaltemp
}%
%
%
% From LaTeX.
%
\def\@nnil{\@nil}%
\def\@fornoop#1\@@#2#3{}%
\def\@for#1:=#2\do#3{%
   \edef\@fortmp{#2}%
   \ifx\@fortmp\empty \else
      \expandafter\@forloop#2,\@nil,\@nil\@@#1{#3}%
   \fi
}%
\def\@forloop#1,#2,#3\@@#4#5{\def#4{#1}\ifx #4\@nnil \else
       #5\def#4{#2}\ifx #4\@nnil \else#5\@iforloop #3\@@#4{#5}\fi\fi
}%
\def\@iforloop#1,#2\@@#3#4{\def#3{#1}\ifx #3\@nnil
       \let\@nextwhile=\@fornoop \else
      #4\relax\let\@nextwhile=\@iforloop\fi\@nextwhile#2\@@#3{#4}%
}%
%
%
% This macro tests if a file \jobname.#1 exists, and sets \if@fileexists
% appropriately.  If an optional argument is given, it is used as the
% root part of the filename instead of \jobname.
%
\@innernewif\if@fileexists
\def\@testfileexistence{\@getoptionalarg\@finishtestfileexistence}%
\def\@finishtestfileexistence#1{%
   \begingroup
      \def\extension{#1}%
      \immediate\openin0 =
         \ifx\@optionalarg\empty\jobname\else\@optionalarg\fi
         \ifx\extension\empty \else .#1\fi
         \space
      \ifeof 0
         \global\@fileexistsfalse
      \else
         \global\@fileexiststrue
      \fi
      \immediate\closein0
   \endgroup
}%
%
%
%% [[[start of BibTeX-specific stuff]]]
%
% Now come the four main LaTeX commands and their associated .aux
% commands.  Just as in LaTeX, \bibliographystyle defines the BibTeX
% style name (.bst file, that is), and \bibliography defines the
% database (.bib) file(s).  The corresponding .aux-file commands are
% \bibstyle and \bibdata, which are there only for BibTeX's (but not
% LaTeX's) use.
%
\def\bibliographystyle#1{%
   \@readauxfile
   \@writeaux{\string\bibstyle{#1}}%
}%
\let\bibstyle = \@gobble
%
% As well as writing the \bibdata command to tell BibTeX which .bib
% files to read, we read the .bbl file that BibTeX (or a person,
% conceivably) has produced.  We use \bblfilebasename as the root of the
% filename to read; this defaults to \jobname.
%
\let\bblfilebasename = \jobname
\def\bibliography#1{%
   \@readauxfile
   \@writeaux{\string\bibdata{#1}}%
   \@testfileexistence[\bblfilebasename]{bbl}%
   \if@fileexists
      % We just output a non-discardable item (the `whatsit' with the
      % \bibdata command).  This means that the glue that will be
      % inserted next (\parskip or \baselineskip, most likely) will be a
      % legal breakpoint.  Most likely, this is after some kind of
      % heading, where we don't want to allow a page break.  So:
      \nobreak
      \@readbblfile
   \fi
}%
\let\bibdata = \@gobble
%
% The \nocite{label,label,...} command writes its argument to \@auxfile,
% unless instructed not to, but produces no text in the document.  Both
% the \nocite and \cite commands produce \citation commands in the .aux file.
%
\def\nocite#1{%
   \@readauxfile
   \@writeaux{\string\citation{#1}}%
}%
\@innernewif\if@notfirstcitation
%
% \cite[note]{label,label,...} produces the citations for the labels as
% well.  If the optional argument `note' is present, it's added after
% the labels.  Since \cite calls \nocite to do its .aux-file writing,
% \cite doesn't need to call \@readauxfile (\nocite does).
%
\def\cite{\@getoptionalarg\@cite}%
%
% Typeset the citations for the labels in #1, followed by the note, if
% it exists.  To change the citation's format in the text, redefine one
% or more `\print...' macros, whose defaults appear later in this file.
%
\def\@cite#1{%
   % Remember the optional argument, in case one of the macros we call
   % below ends up looking for an optional argument itself.  For
   % example, if a \cite[note] triggers reading the .aux file, then the
   % [note] would be clobbered, since \@testfileexistence looks for an
   % optional arg.
   \let\@citenotetext = \@optionalarg
   % Start printing the text, beginning with a left bracket by default.
   \printcitestart
   % It's complicated, but because \nocite puts a `whatsit' onto the list,
   % \nocite should follow \printcitestart.  It's conceivable, but very
   % unlikely, that this `whatsit' will cause a problem (glue that doesn't
   % disappear when you want it to is the most likely symptom), requiring
   % a change either to \printcitestart or to the label that the .bst file
   % produces.
   \nocite{#1}%
   \@notfirstcitationfalse
   \@for \@citation :=#1\do
   {%
      \expandafter\@onecitation\@citation\@@
   }%
   \ifx\empty\@citenotetext\else
      \printcitenote{\@citenotetext}%
   \fi
   \printcitefinish
}%
\newif\ifweareinprivate
\weareinprivatetrue
\ifx\shlhetal\undefinedcontrolseq\weareinprivatefalse\fi
\ifx\shlhetal\relax\weareinprivatefalse\fi
\def\@onecitation#1\@@{%
   \if@notfirstcitation
      \printbetweencitations
   \fi
   \expandafter \ifx \csname\@citelabel{#1}\endcsname \relax
      \if@citewarning
         \message{\@linenumber Undefined citation `#1'.}%
      \fi
      % Give it a dummy definition:
     \ifweareinprivate
      \expandafter\gdef\csname\@citelabel{#1}\endcsname{%
% Change: marginal remark added, goldstrn@math.huji.ac.il, 
% goldstern@tuwien.ac.at, May 1996 mg
%  !!! change !!!
\strut 
\vadjust{\vskip-\dp\strutbox
\vbox to 0pt{\vss\parindent0cm \leftskip=\hsize 
\advance\leftskip3mm
\advance\hsize 4cm\strut\openup-4pt 
\rightskip 0cm plus 1cm minus 0.5cm ?  #1 ?\strut}}
         {\tt
            \escapechar = -1
            \nobreak\hskip0pt\pfeilsw%\special{ps:100 0 0 setrgbcolor }%
            \expandafter\string\csname#1\endcsname
                     %\special{ps:0 0 0 setrgbcolor }
             \pfeilso
            \nobreak\hskip0pt
         }%
      }%
     \else  % ifweareinprivate = false
      \expandafter\gdef\csname\@citelabel{#1}\endcsname{%
            {\tt\expandafter\string\csname#1\endcsname}
      }%
     \fi  % ifweareinprivate
   \fi
   % Now produce the text, whether it was undefined or not.
   \csname\@citelabel{#1}\endcsname
   \@notfirstcitationtrue
}%
%
% Given a label `foo', the macro `\b@foo' is supposed to
% hold the text that should be produced.
%
\def\@citelabel#1{b@#1}%
%
% So, how does a citation label get defined?  When we read the .bbl file
% (below), a \bibitem writes out a \@citedef command.  And when we read
% the \@citedef, we define \@citelabel{#1}, where #1 is the user's
% label.
%
\def\@citedef#1#2{\expandafter\gdef\csname\@citelabel{#1}\endcsname{#2}}%
%
%
% Reading the .bbl file also produces the typeset bibliography.  Please
% notice, however, that we do not produce the title for the references
% (e.g., `References'), as LaTeX does.  The formatting and spacing of
% that title, whether it should go into the headline, and so on, are all
% things determined by your format.  We cannot know those things in
% advance.  If you wish, you can define \bblhook to produce the title.
% Or just do it before the \bibliography command.
%
\def\@readbblfile{%
   % Define a counter to tell us which item number we are on, unless
   % we've already defined it (because the document has more than one
   % bibliography).
   \ifx\@itemnum\@undefined
      \@innernewcount\@itemnum
   \fi
   \begingroup
      \def\begin##1##2{%
         % ##1 is just `thebibliography'.
         % ##2 is the widest label.
         % We set (new dimen) \biblabelwidth based on the widest label
         \setbox0 = \hbox{\biblabelcontents{##2}}%
         \biblabelwidth = \wd0
      }%
      \def\end##1{}% ##1 is `thebibliography' again.
      %
      % Here we have two possibilities:
      % \bibitem[typesetlabel]{citationlabel}
      % \bibitem{citationlabel}
      % If we have the second of these, the citations are numbered, starting
      % from one; we use our own count register \@itemnum for this.
      %
      \@itemnum = 0
      \def\bibitem{\@getoptionalarg\@bibitem}%
      \def\@bibitem{%
         \ifx\@optionalarg\empty
            \expandafter\@numberedbibitem
         \else
            \expandafter\@alphabibitem
         \fi
      }%
      \def\@alphabibitem##1{%
         % Need \xdef here for various reasons.
         \expandafter \xdef\csname\@citelabel{##1}\endcsname {\@optionalarg}%
         % Left-justify alpha labels, unless \biblabel{pre,post}contents
         % are already defined.
         \ifx\biblabelprecontents\@undefined
            \let\biblabelprecontents = \relax
         \fi
         \ifx\biblabelpostcontents\@undefined
            \let\biblabelpostcontents = \hss
         \fi
         \@finishbibitem{##1}%
      }%
      \def\@numberedbibitem##1{%
         \advance\@itemnum by 1
         \expandafter \xdef\csname\@citelabel{##1}\endcsname{\number\@itemnum}%
         % Right-justify numeric labels, unless \biblabel{pre,post}contents
         % are already defined.
         \ifx\biblabelprecontents\@undefined
            \let\biblabelprecontents = \hss
         \fi
         \ifx\biblabelpostcontents\@undefined
            \let\biblabelpostcontents = \relax
         \fi
         \@finishbibitem{##1}%
      }%
      \def\@finishbibitem##1{%
         \biblabelprint{\csname\@citelabel{##1}\endcsname}%
         \@writeaux{\string\@citedef{##1}{\csname\@citelabel{##1}\endcsname}}%
         \ignorespaces
      }%
      %
      % Do the printing (we're producing the bibliography, remember).
      %
      \let\em = \bblem
      \let\newblock = \bblnewblock
      \let\sc = \bblsc
      % Punctuation won't affect spacing;
      \frenchspacing
      % the penalties below are from LaTeX's [article,book,report].sty;
      \clubpenalty = 4000 \widowpenalty = 4000
      % the next two values come from LaTeX's \sloppy command;
      \tolerance = 10000 \hfuzz = .5pt
      \everypar = {\hangindent = \biblabelwidth
                      \advance\hangindent by \biblabelextraspace}%
      \bblrm
      % the \parskip is a guess at what looks good;
      \parskip = 1.5ex plus .5ex minus .5ex
      % and the space between label and text comes from LaTeX's \labelsep.
      \biblabelextraspace = .5em
      \bblhook
      \input \bblfilebasename.bbl
   \endgroup
}%
%
% The widest label's width is useful for redefining \biblabelprint;
% you redefine \biblabelwidth, in effect, by redefining the
% \biblabelcontents macro that appears below.  And \biblabelextraspace,
% which is redefinable inside \bblhook, is added to \biblabelwidth to
% determine the amount of hanging indentation.
%
\@innernewdimen\biblabelwidth
\@innernewdimen\biblabelextraspace
%
% Now come the main macros that are related to the printing of the
% bibliography.  Since you might want to redefine them, they are given
% default definitions outside of \@readbblfile.
%
% The first one controls the printing of a bibliography entry's label.
% If you change it, make sure that it starts with something like
% \noindent or \indent or \leavevmode that puts TeX into horizontal mode
% (even if the label itself is empty); otherwise, the hanging
% indentation will get messed up in certain circumstances.
%
\def\biblabelprint#1{%
   \noindent
   \hbox to \biblabelwidth{%
      \biblabelprecontents
      \biblabelcontents{#1}%
      \biblabelpostcontents
   }%
   \kern\biblabelextraspace
}%
%
% If you are using numeric labels, and you want them left-justified
% (numeric labels by default are right-justified), do something like:
%     \def\biblabelprecontents{\relax}
%     \def\biblabelpostcontents{\hss}
%
% By default the labels are typeset in \bblrm, and enclosed in brackets.
%
\def\biblabelcontents#1{{\bblrm [#1]}}%
%
% The main text, too, is typeset using \bblrm, which is \rm by default.
%
\def\bblrm{\rm}%
%
% Emphasis for producing, e.g., titles, is done with \it by default.
%
\def\bblem{\it}%
%
% Some styles use a caps-and-small-caps font for author names.  LaTeX
% defines an \sc command but plain TeX doesn't, so we need one here.
% The definition below doesn't load the font unless it's needed, but it
% tries to load only the 10pt version, because it might not exist at
% other point sizes.
%
\def\bblsc{\ifx\@scfont\@undefined
              \font\@scfont = cmcsc10
           \fi
           \@scfont
}%
%
% The major parts of an entry are separated with \bblnewblock.  The
% numbers below are taken from LaTeX's `article' style.
%
\def\bblnewblock{\hskip .11em plus .33em minus .07em }%
%
% Here's where you stick any other bibliography-formatting goodies, or
% redefine the values above.
%
\let\bblhook = \empty
%
%
% Here are the four default definitions for formatting the in-text
% citations.  These are what you redefine (after your \input btxmac but
% before your \bibliography) to get parens instead of brackets, or
% superscripts, or footnotes, or whatever.
%
\def\printcitestart{[}%         left bracket
\def\printcitefinish{]}%        right bracket
\def\printbetweencitations{, }% comma, space
\def\printcitenote#1{, #1}%     comma, space, note (if it exists)
%
% That scheme is pretty flexible.  For example you could use
%     \def\printcitestart{\unskip $^\bgroup}
%     \def\printcitefinish{\egroup$}
%     \def\printbetweencitations{,}
%     \def\printcitenote#1{\hbox{\sevenrm\space (#1)}}
%     \font\eighttt = cmtt8
%     \scriptfont\ttfam = \eighttt
% to get superscripted in-text citations.  (The scriptfont stuff
% exists only to print an undefined citation; it's in cmtt8 because
% there is no cmtt7.)  To get something radically different, however,
% you'll have to define your own \cite command.
%
% When we read `\citation' from the .aux file, it means nothing.
%
\let\citation = \@gobble
%
%
% Now comes the stuff for dealing with LaTeX's \newcommand.  As
% mentioned earlier, this \newcommand will redefine a preexisting
% command; that's different from how LaTeX's \newcommand behaves.
%
\@innernewcount\@numparams
%
% \newcommand{\foo}[n]{text} defines the control sequence \foo to have
% n parameters, and replacement text `text'.
%
\def\newcommand#1{%
   \def\@commandname{#1}%
   \@getoptionalarg\@continuenewcommand
}%
%
% Figure out if this definition has parameters.
%
\def\@continuenewcommand{%
   % If no optional argument, we have zero parameters.  Otherwise, we
   % have that many.
   \@numparams = \ifx\@optionalarg\empty 0\else\@optionalarg \fi \relax
   \@newcommand
}%
%
% \@numparams is how many arguments this command has.  The name of the
% command is \@commandname.  The replacement text for the new macro is #1.
%
\def\@newcommand#1{%
   \def\@startdef{\expandafter\edef\@commandname}%
   \ifnum\@numparams=0
      \let\@paramdef = \empty
   \else
      \ifnum\@numparams>9
         \errmessage{\the\@numparams\space is too many parameters}%
      \else
         \ifnum\@numparams<0
            \errmessage{\the\@numparams\space is too few parameters}%
         \else
            \edef\@paramdef{%
               % This is disgusting, but \loop doesn't work inside \edef,
               % because \body isn't defined.
               \ifcase\@numparams
                  \empty  No arguments.
               \or ####1%
               \or ####1####2%
               \or ####1####2####3%
               \or ####1####2####3####4%
               \or ####1####2####3####4####5%
               \or ####1####2####3####4####5####6%
               \or ####1####2####3####4####5####6####7%
               \or ####1####2####3####4####5####6####7####8%
               \or ####1####2####3####4####5####6####7####8####9%
               \fi
            }%
         \fi
      \fi
   \fi
   \expandafter\@startdef\@paramdef{#1}%
}%
%
%% [[[end of BibTeX-specific stuff]]]
%
%
% Names of references (arguments given in the \cite and \nocite
% commands) and file names (arguments given in the \bibliography and
% \bibliographystyle commands) are recorded in \jobname.aux, called the
% \@auxfile in these macros.  Here's how they get read in.
%
\def\@readauxfile{%
   \if@auxfiledone \else % remember: \@auxfiledonetrue if \noauxfile is defined
      \global\@auxfiledonetrue
      \@testfileexistence{aux}%
      \if@fileexists
         \begingroup
            % Because we might be in horizontal mode when \@readauxfile
            % is called (if it's in response to a \cite or \nocite), we
            % want to ignore all the would-be spaces at the ends of
            % lines in the aux file.  Fortunately, it's highly unlikely
            % an end-of-line might actually be desired.
            % And because we don't change the category code of anything
            % but @, primitives like \gdef can't be used to define labels
            % in the aux file.  The solution adopted by btxmac.tex is to
            % write `\@citedef{LABEL}{DEFINITION}' to the aux file, and
            % use \csname on LABEL.
            \endlinechar = -1
            \catcode`@ = 11
            \input \jobname.aux
         \endgroup
      \else
         \message{\@undefinedmessage}%
         \global\@citewarningfalse
      \fi
      \immediate\openout\@auxfile = \jobname.aux
   \fi
}%
%
% The \@readauxfile macro does all that work the first time it's called.
% Since it's called once for every \cite, \nocite, \bibliography, and
% \bibliographystyle command that the user issues, we need to remember
% whether the work's been done.  It's considered done if we're not to do
% it---that is, if \noauxfile is defined.
%
\newif\if@auxfiledone
\ifx\noauxfile\@undefined \else \@auxfiledonetrue\fi
%
% It's conceivable you'd want to change how other characters are read;
% to do that, change their category code before doing \input btxmac.
%
%
% After reading the .aux file, \@readauxfile opens it for writing.
% The \@writeaux macro does the actual writing (as long as
% \noauxfile is undefined).
%
\@innernewwrite\@auxfile
\def\@writeaux#1{\ifx\noauxfile\@undefined \write\@auxfile{#1}\fi}%
%
%
% A macro package that uses btxmac.tex might define
% \@undefinedmessage (before doing an \input btxmac).
%
\ifx\@undefinedmessage\@undefined
   \def\@undefinedmessage{No .aux file; I won't give you warnings about
                          undefined citations.}%
\fi
%
% Even if citations are undefined, we want to complain only if
% \@citewarningtrue.  The default is to set \@citewarningtrue unless
% \noauxfile is defined.  Again, a macro package that uses
% btxmac.tex might want to redefine this.
%
\@innernewif\if@citewarning
\ifx\noauxfile\@undefined \@citewarningtrue\fi
%
%
% Finally, before leaving we restore @'s old category code.
%
\catcode`@ = \@oldatcatcode

\def\pfeilso{\leavevmode
            \vrule width 1pt height9pt depth 0pt\relax
           \vrule width 1pt height8.7pt depth 0pt\relax
           \vrule width 1pt height8.3pt depth 0pt\relax
           \vrule width 1pt height8.0pt depth 0pt\relax
           \vrule width 1pt height7.7pt depth 0pt\relax
            \vrule width 1pt height7.3pt depth 0pt\relax
            \vrule width 1pt height7.0pt depth 0pt\relax
            \vrule width 1pt height6.7pt depth 0pt\relax
            \vrule width 1pt height6.3pt depth 0pt\relax
            \vrule width 1pt height6.0pt depth 0pt\relax
            \vrule width 1pt height5.7pt depth 0pt\relax
            \vrule width 1pt height5.3pt depth 0pt\relax
            \vrule width 1pt height5.0pt depth 0pt\relax
            \vrule width 1pt height4.7pt depth 0pt\relax
            \vrule width 1pt height4.3pt depth 0pt\relax
            \vrule width 1pt height4.0pt depth 0pt\relax
            \vrule width 1pt height3.7pt depth 0pt\relax
            \vrule width 1pt height3.3pt depth 0pt\relax
            \vrule width 1pt height3.0pt depth 0pt\relax
            \vrule width 1pt height2.7pt depth 0pt\relax
            \vrule width 1pt height2.3pt depth 0pt\relax
            \vrule width 1pt height2.0pt depth 0pt\relax
            \vrule width 1pt height1.7pt depth 0pt\relax
            \vrule width 1pt height1.3pt depth 0pt\relax
            \vrule width 1pt height1.0pt depth 0pt\relax
            \vrule width 1pt height0.7pt depth 0pt\relax
            \vrule width 1pt height0.3pt depth 0pt\relax}

\def\pfeilsw{ \leavevmode 
            \vrule width 1pt height0.3pt depth 0pt\relax
            \vrule width 1pt height0.7pt depth 0pt\relax
            \vrule width 1pt height1.0pt depth 0pt\relax
            \vrule width 1pt height1.3pt depth 0pt\relax
            \vrule width 1pt height1.7pt depth 0pt\relax
            \vrule width 1pt height2.0pt depth 0pt\relax
            \vrule width 1pt height2.3pt depth 0pt\relax
            \vrule width 1pt height2.7pt depth 0pt\relax
            \vrule width 1pt height3.0pt depth 0pt\relax
            \vrule width 1pt height3.3pt depth 0pt\relax
            \vrule width 1pt height3.7pt depth 0pt\relax
            \vrule width 1pt height4.0pt depth 0pt\relax
            \vrule width 1pt height4.3pt depth 0pt\relax
            \vrule width 1pt height4.7pt depth 0pt\relax
            \vrule width 1pt height5.0pt depth 0pt\relax
            \vrule width 1pt height5.3pt depth 0pt\relax
            \vrule width 1pt height5.7pt depth 0pt\relax
            \vrule width 1pt height6.0pt depth 0pt\relax
            \vrule width 1pt height6.3pt depth 0pt\relax
            \vrule width 1pt height6.7pt depth 0pt\relax
            \vrule width 1pt height7.0pt depth 0pt\relax
            \vrule width 1pt height7.3pt depth 0pt\relax
            \vrule width 1pt height7.7pt depth 0pt\relax
            \vrule width 1pt height8.0pt depth 0pt\relax
            \vrule width 1pt height8.3pt depth 0pt\relax
            \vrule width 1pt height8.7pt depth 0pt\relax
            \vrule width 1pt height9pt depth 0pt\relax
      }

%  *** end including bib4plain.tex *** 
  % This will define \cite and make sure it works as in latex

\def\widestnumber#1#2{}
  % Our amstex-ppt style does not know about \widestnumber

\def\citewarning#1{\ifx\shlhetal\relax 
   % normal mode, do not write anything
    \else
   % private mode
    \par{#1}\par
    \fi
}

\def\rm{\fam0 \tenrm}

\def\fakesubhead#1\endsubhead{\bigskip\noindent{\bf#1}\par}

% % \input rsfs
%  *** start including rsfs.tex *** 

% # Keywords: Script or Calligraphic (Caligraphic) letters with the RSFS Font

% The story so far:    July 1998 -- Saharon would like to have a
% ``nicer'' calligraphic font. In particualr, the leters S and P in
% the usual calligraphic font do not look ``special'' enough. 
% 
% I found out that ``rsfs'' (``Ralph Smith Formal Script'') may be
% what he wants.   I installed the mf file, the .tfm file, as well as
% a few pk files in ~/TeX/rsfs.    Let's hope that this is enough.
% Using amstex, all you have to do is to \input rsfs.tex 
% Files prepared with citealice willdothis automatically. 
%
%  Note:  for some reason xdvi calls MakeTeXpk, then Maketexpk
%  complains about wrong resolution, but still writes commands to
%  missfont.log...  
%

% we redefine a macro inside amstex's \Cal command , so that it calls
% our nice font ``rsfs'' rather than the usual calligraphic font. 
% Note thisworks for amstex only.   
% In plain tex, would have to add definitions of \Cal
% in latex... we should insteaduse mathrsfs.sty
% 

\font\textrsfs=rsfs10
\font\scriptrsfs=rsfs7
\font\scriptscriptrsfs=rsfs5

\newfam\rsfsfam
\textfont\rsfsfam=\textrsfs
\scriptfont\rsfsfam=\scriptrsfs
\scriptscriptfont\rsfsfam=\scriptscriptrsfs

\edef\oldcatcodeofat{\the\catcode`\@}
\catcode`\@11

\def\Cal@@#1{\noaccents@ \fam \rsfsfam #1}

\catcode`\@\oldcatcodeofat

%  *** end including rsfs.tex *** 

\expandafter\ifx \csname margininit\endcsname \relax\else\margininit\fi

\long\def\red#1\endred{}
\long\def\green#1\endgreen{}
\long\def\blue#1\endblue{}
\long\def\private#1\endprivate{}

\def\endred{ \unmatched endred! }
\def\endgreen{ \unmatched endgreen! }
\def\endblue{ \unmatched endblue! }
\def\endprivate{ \unmatched endprivate! }

\ifx\latexcolors\undefinedcs\def\latexcolors{}\fi

\def\emptycs{}
\def\evaluatelatexcolors{%
        \ifx\latexcolors\emptycs\else
        \expandafter\xxevaluate\latexcolors\xxfertig\evaluatelatexcolors\fi}
\def\xxevaluate#1,#2\xxfertig{\setupthiscolor{#1}%
        \def\latexcolors{#2}}

 % \def\definedprivate{rgb 0.5 0 0.5} 
       % black 

\font\smallfont=cmsl7
\def\rutgerscolor{\ifmmode\else\endgraf\fi\smallfont% \vrule width 3cm height 1cm depth 0cm
\advance\leftskip0.5cm\relax}
\def\setupthiscolor#1{\edef\tmptmpcs{\noexpand\bgroup\noexpand\rutgerscolor
\noexpand\def\noexpand\currentcolor{#1}%
\noexpand}%
\expandafter\let\csname#1\endcsname\tmptmpcs
\def\tmptmpcs{\checkColorUnmatched{#1}\popthecolor}
\expandafter\let\csname end#1\endcsname\tmptmpcs}

\def\checkColorUnmatched#1{\def\expectcolor{#1}%
    \ifx\expectcolor\currentcolor   % OK! 
    \else \edef\failhere{\noexpand\tryingToClose '\currentcolor' with end\expectcolor}\failhere\fi}

\def\currentcolor{???}

\def\popthecolor{\ifmmode\else\endgraf\fi\egroup}

\expandafter\def\csname#1\endcsname{}

\evaluatelatexcolors

 \let\outerhead\head
 \def\head{\innerhead}
 \let\innerhead\outerhead

 \let\outersubhead\subhead
 \def\subhead{\innersubhead}
 \let\innersubhead\outersubhead

 \let\outersubsubhead\subsubhead
 \def\subsubhead{\innersubsubhead}
 \let\innersubsubhead\outersubsubhead

 \def\proclaim{\innerproclaim}
 \let\innerproclaim\outerproclaim

 % 
 % \newdimen\wzero 
 % \newdimen\hzero 
 % \newdimen\dzero 
 % 
 % \def\whiteblack#1{\ifx \shlhetal\relax      
 %    \hbox{#1}% normal mode
 %     \else\ifx\shlhetal\undefinedcontrolsequence
 %    \hbox{#1}% normal mode
 %         \else
 %    % private mode
 %     \setbox0=\hbox{#1}\leavevmode
 %     \wzero=\wd0\hzero=\ht0\dzero=\dp0
 %         \advance\wzero1mm
 %         \advance\dzero1mm
 %         \advance\hzero1mm
 %     \hbox{\vrule width \wzero height \hzero depth\dzero\hskip-\wzero
 %         \special{ps: currentgray 1  setgray}%
 %         \hbox to \wzero{\hss\copy0\hss}%
 %            \special{ps:setgray}}\fi\fi
 % }
 % 
 % 

\def\demo#1{\medskip\noindent{\it #1.\/}}
\def\enddemo{\smallskip}

%  *** end including alice2jlem.tex *** 
\pageheight{8.5truein}
%\pageheight{48.5pc}
\topmatter
\title{Power set modulo small, the singular of uncountable cofinality} \endtitle
\rightheadtext{Power set}
\author {Saharon Shelah \thanks {\null\newline 
This research was supported by the United States-Israel Binational
Science Foundation. Publication 861.  \null\newline
I would like to thank 
Alice Leonhardt for the beautiful typing.
} \endthanks} \endauthor 
% First Typed - 04/Aug/25
% Previous version - 05/Mar/25
% Formerly F475 - changed 04/Sept/28

\affil{The Hebrew University of Jerusalem \\
Einstein Institute of Mathematics \\
Edmond J. Safra Campus, Givat Ram \\
Jerusalem 91904, Israel
 \medskip
 Department of Mathematics \\
 Hill Center-Busch Campus \\
  Rutgers, The State University of New Jersey \\
 110 Frelinghuysen Road \\
 Piscataway, NJ 08854-8019 USA} \endaffil

\abstract  Let $\mu$ be singular of uncountable cofinality.  If $\mu >
2^{\text{cf}(\mu)}$, we prove that in $\Bbb P =
([\mu]^\mu,\supseteq)$ as a forcing notion we have a natural complete
embedding of Levy$(\aleph_0,\mu^+)$ (so $\Bbb P$ collapses $\mu^+$ to
$\aleph_0$) and even Levy$(\aleph_0,\bold
U_{J^{\text{bd}}_\kappa}(\mu))$.  The ``natural" means that the forcing
$(\{p \in [\mu]^\mu:p$ closed$\},\supseteq)$ is naturally embedded
and is equivalent to the Levy algebra.  Moreover we 
prove more than conjectured: if $\Bbb P$ fails the
$\chi$-c.c. then it collapses $\chi$ to $\aleph_0$. 
We even prove the parallel
results for the case $\mu > \aleph_0$ is regular or of countable
cofinality. We also prove: for regular uncountable $\kappa$, 
there is a family ${\bold P}$ of ${\frak b}_\kappa$ partitions 
$\bar A=\langle A_\alpha:\alpha<\kappa\rangle$ of $\kappa$ 
such that for any $A\in [\kappa]^\kappa$ for some 
$\langle A_\alpha:\alpha<\kappa\rangle \in {\bold P}$ 
we have $\alpha<\kappa\Rightarrow |A_\alpha
\cap A|=\kappa$.
\endabstract
\endtopmatter
\document

\newpage

\head {\S0 Introduction} \endhead  \resetall \sectno=0
 \spuriousreset
\bigskip

This work on the one hand continue the celebrated work of the 
Czech school on the completion of the Boolean algebras
${\Cal P} (\lambda)/[\lambda]^{<\lambda}$ solving some of their 
questions and on the other hand
tries to confirm the ``pcf is effective" thesis.  

We may consider the completions of the Boolean Algebras ${\Cal
 P}(\mu)/\{u \subseteq \mu:|u| < \mu\} = {\Cal P}(\mu)/[\mu]^{<
 \mu}$. This is equivalent to considering the partial orders $\Bbb P_\mu =
([\mu]^\mu,\supseteq)$, viewing them as forcing notions, so actually
looking at their completion $\hat{\Bbb P}_\mu$, which are complete
Boolean Algebras.  Recall that forcing notions
$\Bbb P^1,\Bbb P^2$ are equivalent iff
their completions are isomorphic Boolean Algebras.  The Czech school
has investigated them, in particular, (letting $\ell(\mu)$ be 0 if
cf$(\mu) > \aleph_0$ and 1 if $\mu > \text{\rm cf}(\mu) = \aleph_0$,
(and $\aleph_{\ell(\mu)} = {\frak h}$ if $\mu = \aleph_0$) consider
 the questions:
\mr
\item "{$\otimes_1$}"  $(a) \quad$ is $\hat{\Bbb P}_\mu$ isomorphic to
 the completion of the Levy collapse Levy$(\aleph_{\ell(\mu)},2^\mu)$?
\sn
\item "{${{}}$}"  $(b) \quad$ which cardinals $\chi$ the forcing
notion $\Bbb P_\mu$ collapse to $\aleph_{\ell(\mu)}$ 
in particular is $\mu^+$ 
collapsed 
\sn
\item "{${{}}$}"  $(c) \quad$ is $\Bbb P_\mu$  $(\theta,\chi)$-nowhere
distributive for $\theta = \aleph_{\ell(\mu)}$? This can \nl

\hskip25pt be phrased as: for some $\Bbb P_\mu$-name 
$\underset\tilde {}\to f$ of a function from
$\aleph_{\ell(\mu)}$ to $\chi$, \nl 

\hskip25pt for every $p \in \Bbb P_\mu$ for some
$i < \theta$ the set $\{\alpha < \chi:p \nVdash {\underset\tilde
{}\to f}(i) \ne \alpha\}$ \nl

\hskip25pt has cardinality $\chi$.
\ermn
The first, (a) is a full answer, the second, (b)seems central 
for set theories and essentially give sufficient condition for the
first,
the last is sufficient if the density is right, to get the first.  The
case of collapsing seems central (it also implies clause (c)) so we 
repeat the summary from Balcar, Simon \cite{BaSi95} of what was 
 known of the collapse of cardinals by $\Bbb P_\mu$, i.e.,
$\otimes_1(b)$.   Let $\chi \rightarrow_\mu \theta$ 
denote the fact that $\chi$ is collapsed to $\theta$ by $\Bbb P_\mu$
\mr
\item "{$\boxtimes_1$}"  $(i) \quad$ for 
$\mu = \aleph_0,2^{\aleph_0} \rightarrow_\mu {\frak h}$, (but $\Bbb
P_\mu$ adds no new sequence of length $< {\frak h}$ so \nl

\hskip25pt  we are done), Balcar, Pelant, Simon \cite{BPS}
\sn
\item "{${{}}$}"  $(ii) \quad$ for $\mu$ uncountable and regular, 
${\frak b}_\mu \rightarrow_\mu \aleph_0$, (hence $\mu^+
\rightarrow_\mu \aleph_0$), 
\nl

\hskip25pt Balcar, Simon \cite{BaSi88}
\sn
\item "{${{}}$}"  $(iii) \quad$ for $\mu$ singular with cf$(\mu) =
\aleph_0,2^{\aleph_0} \rightarrow_\mu \aleph_1$, Balcar, Simon \cite{BaSi95}
\sn
\item "{${{}}$}"  $(iv) \quad$ for $\mu$ singular with cf$(\mu) \ne
\aleph_0,{\frak b}_{\text{cf}(\mu)} \rightarrow_\mu
\aleph_0$, 
\nl

\hskip25pt Balcar, Simon \cite{BaSi95}; 
\ermn
under additional assumptions on cardinal arithmetic 
for singular cardinals more is known
\mr 
\item "{${{}}$}"  $(v) \quad$ for $\mu$ singular with cf$(\mu) =
\aleph_0$ and $\mu^{\aleph_0} = 2^\mu,\mu^{\aleph_0}
\rightarrow_\mu \aleph_1$, \nl

\hskip25pt Balcar, Simon \cite{BaSi88}
\sn
\item "{${{}}$}"  $(vi) \quad$ for $\mu$ singular 
with cf$(\mu) \ne \aleph_0$
and $2^\mu = \mu^+,2^\mu \rightarrow_\mu \aleph_0$, \cite{BaSi88}.
\endroster
\bn

Now \cite{BaSi95} finish with the following very reasonable conjecture.
\nl
\margintag{0.0}\ub{\stag{0.0} Conjecture}:  (Balcar and Simon) in ZFC: for a singular
cardinal $\mu$ with countable cofinality, $\mu^{\aleph_0}
\rightarrow_\mu \aleph_1$ and for a singular cardinal $\mu$ with an
uncountable cofinality $\mu^+ \rightarrow_\mu \aleph_0$
(here we concentrate on the case cf$(\mu) > \aleph_0$, see below).

Concerning the other questions they prove
\mr
\item "{$\boxtimes_2$}"  $(i) \quad$ Balcar, Franek \cite{BaFr87}:
\nl

\hskip25pt if $\mu > \text{\rm cf}(\mu) >
\aleph_0,2^{\text{cf}(\mu)} = \text{ cf}(\mu)^+$ then $\Bbb P_\mu$ is
$(\aleph_0,\mu^+)$-nowhere distributive
\sn
\item "{${{}}$}"  $(ii) \quad$ Balcar, Simon \cite[5.20, pg.380]{BaSi89}:
\nl

\hskip25pt if $2^\mu = \mu^+$ and 
$2^{\text{cf}(\mu)} = \text{\rm cf}(\mu)^+$ \ub{then} 
$\Bbb P_\mu$ is equivalent: 
\nl

\hskip25pt to Levy$(\aleph_0,\mu^+)$ if cf$(\mu) > \aleph_0$ and
\nl

\hskip25pt  to Levy$(\aleph_1,\mu^+)$ if cf$(\mu) = \aleph_0$
\sn
\item "{${{}}$}"  $(iii) \quad$ Balcar, Franek \cite{BaFr87}:
\nl

\hskip25pt if $2^\mu = \mu^+,\mu = \text{\rm cf}(\mu) > \aleph_0,J$ a
$\mu$-complete ideal on $\mu$ and
\nl

\hskip25pt $J$ is nowhere precipitous extending $[\mu]^{< \mu}$ then ${\Cal
P}(\mu)/J$ is equivalent 
\nl

\hskip25pt to Levy$(\aleph_0,\mu^+)$; also the parallel of (ii).
\ermn
So under G.C.H. the picture was complete; getting clause $(ii)$ of
$\boxtimes_2$. Also under ZFC for regular cardinals $\mu > \aleph_0$ the
picture is reasonable, particularly if we recall that by 
Baumgartner \cite{Ba}
\mr
\item "{$\boxtimes_3$}"  if $\kappa = \text{\rm cf}(\mu) < \theta
= \theta^{< \theta} < \mu < \chi$ and $\bold V \models$
G.C.H. for simplicity and $\Bbb P$ is forcing for adding $\chi$ Cohen subsets to 
$\theta$ then
{\roster
\itemitem{ $(a)$ }   forcing with $\Bbb P$ collapses no cardinal,
changes no cofinality, adds no new sets of $< \theta$ ordinals
\sn
\itemitem{ $(b)$ }  in $\bold V^{\Bbb P},([\mu]^\mu,\supseteq)$ satisfies
the $\mu^+_1$-c.c. where $\mu_1 = (2^\mu)^{\bold V}$;
hence does not collapse any cardinal $\geq \mu_1^+$.
\endroster}
\endroster
\bn
Lately, Kojman, Shelah \cite{KjSh:720} prove the conjecture \scite{0.0}
for the case when $\mu > \text{ cf}(\mu) = \aleph_0$; morever
\mr
\item "{$\boxtimes_4$}"   $(i) \quad$ if $\mu > \text{\rm cf}(\mu) = \aleph_0$
then Levy$(\aleph_1,\mu^{\aleph_0})$ can be
completely embedded 
\nl

\hskip25pt into the completion of $\Bbb P_\mu$.
Moreover,
\sn
\item "{${{}}$}"  $(ii) \quad$  the embedding is ``natural":
Levy$(\aleph_1,\mu^{\aleph_0})$ is equivalent to ${\Bbb Q}_\mu$ 
which is $\lessdot {\Bbb P}_\mu$ where
\nl

\hskip25pt $\Bbb Q_\mu = (\{A \subseteq \mu:A$ a closed subset of $\mu$ of
cardinality $\mu\},\supseteq)$.
\ermn
Here we continue \cite{KjSh:720} in \S1, \cite{BaSi89} in \S2 but make
it self contained.  
Both sections use results on pcf (in addition to guessing clubs)
Naturally we may add to the questions (answered positively for the
case cf$(\mu) = \aleph_0$ by \cite{KjSh:720}) 
\mr
\item "{$\otimes_2$}"  $(a)\quad$ can we strengthen ``$\Bbb P_\mu$
collapse $\chi$ to $\aleph_{\ell(\mu)}$" to
``Levy$(\aleph_{\ell(\mu)},\chi)$ is 
\nl

\hskip25pt completely embeddable into $\Bbb P_\mu$ (really $\hat{\Bbb P}_\mu$)"
\sn
\item "{${{}}$}"  $(b) \quad$ can we find natural such embeddings.
\ermn
We may add that by \cite{BaSi95} the Baire number of  ${\Cal U}[\mu]$, the
space of all uniform ultrafilters over uncountable $\mu$ is $\aleph_1$, except
when $\mu > \text{\rm cf}(\mu) = \aleph_0$ and in that case it is
$\aleph_2$
under some reasonable assumptions.  By \cite{KjSh:720}
the Baire number of ${\Cal U}[\mu]$ is always $= \aleph_2$ 
when $\mu > \text{\rm cf}(\mu) = \aleph_0$.
\mn
Our original aim in this work has been to deal with $\mu > \text{\rm cf}(\mu) >
\aleph_0$, proving the conjecture of Balcar and Simon above (i.e.,
that $\mu^+$ is collapsed to $\aleph_0$), first of all when
$2^{\text{cf}(\mu)} < \mu$ answering $\otimes_2(a)+(b)$ using pcf
(and replacing $\mu^+$ by pp$_{J^{\text{bd}}_{\text{cf}(\mu)}}(\mu))$.
In fact this seems, at least to me, the best we can reasonably
expect.  But a posteriori we have more to say.

For $\mu = \kappa = \text{\rm cf}(\mu) > \aleph_0$,
though by the above we know that some cardinal $> \mu$ is collapsed  
(that is ${\frak b}_\kappa$), we do not 
know what occurs up to $2^\mu$ or when the
c.c. fails.  This leads to the following conjecture, (stronger than
the Balcar, Simon one mentioned above).
Of course, it naturally breaks to cases according to $\mu$. 

\demo{\stag{0.2} Conjecture}   If $\mu > \aleph_0$ and $\Bbb P_\mu$
does not satisfy the $\chi$-c.c., \ub{then} forcing with $\Bbb P_\mu$
collapse $\chi$ to $\aleph_{\ell(\mu)}$, see Definition \scite{0.5}
below.
\enddemo
\bn
Note that
\demo{\stag{0.3} Observation}
If conjecture \scite{0.2} holds for
$\mu > \aleph_0$ then $\Bbb P_\mu$ is equivalent to a Levy
collapse iff it fails the $d(\Bbb P_\mu)$-c.c. where 
$d(\Bbb P_\mu)$ is the density of $\Bbb P_\mu$.  
\enddemo
\bn
Lastly, we turn to the results; by \scite{1.9}(1):
\proclaim{\stag{0.1} Theorem}  If $\mu > \kappa = \text{\rm cf}(\mu) >
\aleph_0$ and $\mu > 2^\kappa$ \ub{then} $\Bbb Q_\mu$ (a natural
complete subforcing of $\Bbb P_\mu$, forcing with closed sets) is
equivalent to Levy$(\aleph_0,\bold U_{J^{\text{bd}}_\kappa}(\mu))$.
\endproclaim
\bn
By \scite{1.11}, \scite{1.12} and \scite{4b.4} we have
\proclaim{\stag{0.4} Theorem}  Conjecture \scite{0.2} holds except
possibly when $\aleph_0 < \text{\rm cf}(\mu) < \mu < 2^{\text{cf}(\mu)}$.
\endproclaim
\bn
We shall in a
subsequent paper prove the Balcar, Simon conjecture fully, i.e., in all
cases.
\definition{\stag{0.5} Definition}  For $\mu > \aleph_0$ we define
$\ell(\mu) \in \{0,1\}$ by

$\ell(\mu) = 0$ if cf$(\mu) > \aleph_0$

$\ell(\mu) = 1$ if $\mu > \text{\rm cf}(\mu) = \aleph_0$
\nl

and may add

$\ell(\mu) = \alpha$ when $\mu = \aleph_0,{\frak h} = \aleph_\alpha$.
\enddefinition
\bn
We thank Menachem Kojman for discussions on earlier attempts,
Shimoni Garti for corrections and Bohuslav Bakar 
and Pek Simon for improving 
the presentation.
\newpage

\head {\S1 Forcing with closed set is equivalent to the Levy algebra} \endhead  \resetall \sectno=1
 \spuriousreset
\bigskip

\definition{\stag{1.3.7} Definition}  1) For $f \in 
{}^\kappa(\text{\rm Ord} \backslash \{0\})$ and ideal $I$ on $\kappa$ let

$$
\align
\bold U_I(f) = \text{\rm Min}\{|{\Cal P}|:&{\Cal P} \subseteq
[\text{\rm sup Rang}(f)]^{\le \kappa} \\
  &\text{\rm such that for every } g \le f \,{ \text{\rm for some }} u \in
  {\Cal P} \\
  &\text{\rm we have } \{i < \kappa:g(i) \in u\} \in I^+\}.
\endalign
$$
\mn
2) Let $\bold U_I(\lambda)$ means ${\bold U}_I(f)$ where $f$ is the
   function with domain Dom$(I)$ which is constantly $\lambda$
\enddefinition
\bigskip

\demo{\stag{1.1} Hypothesis}
\mr
\item "{$(a)$}"  $\mu$ is a singular cardinal
\sn
\item "{$(b)$}"  $\kappa = \text{\rm cf}(\mu) > \aleph_0$.
\endroster
\enddemo
\bigskip

\definition{\stag{1.2} Definition}  1) $\Bbb P_\mu$ is the following
forcing notion

$$
p \in \Bbb P_\mu \text{ \ub{iff} } p \in [\mu]^\mu
$$

$$
\Bbb P_\mu \models p \le q \text{ iff } p \supseteq q.
$$
\mn
2) $\Bbb P'_\mu$ is the forcing notion with the same set of elements and
with the partial order

$$
\Bbb P'_\mu \models p \le q \text{ iff } |q \backslash p| < \mu.
$$
\mn
3) $\Bbb Q_\mu = \Bbb Q^0_\mu$ is $\Bbb P_\mu \restriction \{p \in
\Bbb P_\mu:p$ is closed in the order topology of $\mu \}$.
\enddefinition
\bigskip

\definition{\stag{1.3.2} Choice/Definition}  1) Let $\langle
\lambda_i:i < \kappa \rangle$ be an increasing sequence of regular
cardinals $> \kappa$ with limit $\mu$.
\nl
2) Let $\lambda^-_i = \cup\{\lambda_j:j <i\}$.
\nl
3) For $p \in \Bbb P_\mu$ let $a(p) = \{i < \kappa:p \cap
[\lambda^-_i,\lambda_i) \ne \emptyset\}$.
\nl
4) $\Bbb Q^1_\mu = \{p \in \Bbb P_\mu:i < \kappa \Rightarrow
|p \cap \lambda_i| < \lambda_i$
and for each $i \in a(p)$ the set $p \cap \lambda_i \backslash
\lambda^-_i$ has no last element, is
closed in its supremum and has cardinality $> |p \cap \lambda^-_i|\}$.
\nl
5) For $p \in \Bbb Q^1_\mu$ let ch$_p \in \dsize \prod_{i \in a(p)}
\lambda_i$ be ch$_p(i) = \cup\{\alpha + 1:\alpha \in p \cap
[\lambda^-_i,\lambda_i)\}$ and  cf$_p \in \dsize \prod_{i \in a(p)}
\lambda_i$ be cf$_p(i) = \text{\rm cf}(\text{\rm ch}_p(i))$.
\nl
6) $\Bbb Q^2_\mu = \{p \in \Bbb Q^1_\mu$:cf$_p(i) > |p \cap
\lambda^-_i|$ for $i \in a(p)\}$.
\comment
\nl
7) $p \in \Bbb P_\mu$ is $\bar\lambda'$-normal when $\bar\lambda' =
\langle \lambda'_i:i \in a(p)\rangle$ and otp$(p \cap
[\lambda^-_i,\lambda_i)) = \lambda'_i$ for $i \in a(p)$.
\endcomment
\enddefinition
\bigskip

\proclaim{\stag{1.3} Claim}  1) $\Bbb Q^0_\mu,\Bbb Q^1_\mu,\Bbb Q^2_\mu$ are 
complete sub-forcings of $\Bbb P_\mu$.
\nl
2) For $\ell=0,1,2$ and $p,q \in \Bbb Q^\ell_\mu$ we have 
$p \Vdash_{\Bbb Q^\ell_\mu} ``q
\in \underset\tilde {}\to G"$ iff $|p\backslash q| < \mu$ and
similarly for $\Bbb P_\mu$.
\nl
3) $\Bbb Q_\mu = \Bbb Q^0_\mu,\Bbb Q^1_\mu,\Bbb Q^2_\mu$ are
equivalent, in fact $\Bbb Q^2_\mu$ is a dense subset of $\Bbb Q^1_\mu$ and
for $\ell=0,1,\{p /\approx:p \in \Bbb Q^\ell_\mu\}$ does not 
depend on $\ell$ where 
$\approx$ is the equivalence relation of $\Bbb P_\mu$, defined by $p_1
\approx p_2$ iff $(\forall q \in \Bbb P_\mu)(q
\Vdash_{{\Bbb P}_\mu} p_1 \in \underset\tilde {}\to G \Leftrightarrow
q \Vdash_{\Bbb P_\mu} p_2 \in \underset\tilde {}\to G)$.
\endproclaim
\bigskip

\demo{Proof}  Easy.
\enddemo
\bn
Recall
\proclaim{\stag{1.3.3} Claim}  1) $\Bbb P_\kappa$ can be completely
embedded into $\Bbb P_\mu$ (naturally).
\nl
2) $\Bbb Q_\mu$ can be completely embedded into $\Bbb P_\mu$
   (naturally).
\nl
3) $\Bbb P_\kappa$ is completely embeddable into $\Bbb Q_\mu$ (naturally).
\endproclaim
\bigskip

\demo{Proof}  1) Known: just $a \in [\kappa]^\kappa$ can be mapped to
$\cup\{[\lambda^-_i,\lambda_i):i \in a\}$.
\nl
2) By \cite[2.2]{KjSh:720}.
\nl
3) Should be clear (map $A \in [\kappa]^\kappa$ to
$\cup\{[\lambda^-_i,\lambda_i]:i \in A\}$).  \hfill$\square_{\scite{1.3.3}}$\margincite{1.3.3}
\enddemo
\bigskip

\definition{\stag{1.4.2} Choice/Definition}   $\lambda_* = 
\bold U_{J^{\text{bd}}_\kappa}(\mu)$. 
\comment
\nl
2) $\chi$ is, e.g., $(\beth_8(\lambda_*))^+,
<^*_\chi$ a well ordering of ${\Cal H}(\chi)$.
\nl
3) ${\frak B}$ is an elementary submodel of $({\Cal H}(\chi),
\in,<^*_\chi)$ of cardinality $\lambda_*$ such that $\lambda_* +1
\subseteq {\frak B}$.
\endcomment
\enddefinition
\bn
Recall
\proclaim{\stag{1.5} Claim}  Assume $\mu > 2^\kappa$.
\nl
1) $\lambda_* = 
\sup\{\text{\rm pp}_{J^{\text{bd}}_\kappa}(\mu'):\kappa < \mu' 
\le \mu$, cf$(\mu') = \kappa\} = 
\sup\{\text{\rm tcf}(\dsize \prod_{i < \kappa}
\lambda'_i/J^{\text{bd}}_\kappa):\lambda'_i \in \text{\rm Reg} \cap
(\kappa,\mu)$ and $\dsize \prod_{i < \kappa}
\lambda'_i/J^{\text{bd}}_\kappa$ has true cofinality$\}$.
\nl
2) For every regular cardinal $\theta \in [\mu,\lambda_*]$, 
for some increasing sequence $\langle
\lambda^*_i:i < \kappa \rangle$ of regulars $\in (\kappa,\mu)$ we
have $\theta = \text{\rm tcf}(\dsize \prod_{i < \kappa} 
\lambda^*_i,<_{J^{\text{bd}}_\kappa})$.
\nl
3) $\lambda_*=|{\Cal P}|$ where ${\Cal P}\subseteq [\mu]^\kappa$ 
is any maximal almost disjoint family. 
\endproclaim
\bigskip

\demo{Proof}  1) Note that $J^{\text{bd}}_\kappa \restriction A
\approx J^{\text{bd}}_\kappa$ if $A \in (J^{\text{bd}}_\kappa)^+$, we
use this freely.  By their definitions the second and third terms are equal.
Also by the definition the second is smaller or equal to the first.

By \cite[1.1]{Sh:589}, the first, $\lambda_* = \bold
U_{J^{\text{bd}}_\kappa}(\mu)$ is $\le$ than the second number (well
it speaks on $T^2_{J^{\text{bd}}_\kappa}(\mu)$ instead of $\bold
U_{J^{\text{bd}}_\kappa}(\mu)$ but as $2^\kappa < \mu$ they are the
same).
\nl
2) By \cite[1.1]{Sh:589} we actually get the
stronger conclusion.  
\nl
3) It follows easily from the definitions \scite {1.3.7} and \scite{1.4.2}, and
from the inequalities $2^\kappa<\mu<\lambda_*$.
\hphantom{s}\hfill$\square_{\scite{1.5}}$\margincite{1.5}
\enddemo
\bigskip

\definition{\stag{1.4.4} Claim/Definition}  Fix a set 
$\Cal P\subseteq [\mu]^\kappa$ exemplifying $\lambda_* = \bold
 U_{J^{\text{bd}}_\kappa}(\mu)$. 
\nl
1) There is $\bar C^* =
\langle C^*_\alpha:\alpha < \mu \rangle$ such that:
\mr
\item "{$(a)$}"  $C^*_\alpha$ is a subset of
$[\lambda^-_i,\lambda_i)$ closed in its supremum when $\alpha \in
(\lambda^-_i,\lambda_i]$
\sn
\item "{$(b)$}"   if $i < \kappa,\gamma < \lambda_i$, 
$\gamma$ is a regular cardinal 
and $C$ is a closed
subset of $[\lambda^-_i,\lambda_i)$ of order type $\gamma^{++}$,
then for some $\alpha \in (\lambda^-_i,\lambda_i),C^*_\alpha \subseteq
C$ and otp$(C^*_\alpha) = \gamma$.
\ermn
2) $\Bbb Q^3_\mu = \{p \in \Bbb Q^2_\mu:$ if $i \in a(p)$ then $p \cap
[\lambda^-_i,\lambda_i) \in \{C^*_\alpha:\alpha \in 
(\lambda^-_i,\lambda_i)\}$ and for some $u\in\Cal P$, 
$\{\alpha<\mu:$ for some $i<\kappa, p \cap
[\lambda^-_i,\lambda_i)=C^*_\alpha\}\subseteq u\}$ 
is a dense subset of $\Bbb Q^1_\mu,\Bbb
Q^2_\mu$, hence of $\Bbb Q_\mu$.
\nl
3) For $p \in \Bbb Q^3_\mu$ let cd$_p \in \dsize \prod_{i \in a(p)}
\lambda_i$ be such that cd$_p(i) \in (\lambda^-_i,\lambda_i)$ is the minimal
$\alpha \in [\lambda^-_i,\lambda_i)$ such that $p \cap
[\lambda^-_i,\lambda_i) = C^*_\alpha$. 
Notice that for every $p\in\Bbb Q^3_\mu$, 
there is some $u\in\Cal P$ with Rang$(\text{cd}_p)\subseteq u$.
\enddefinition
\bigskip

\demo{Proof}  1) It is enough, for any limit $\delta \in
(\lambda^-_i,\lambda_i)$ and regular $\theta,\theta^+ < \text{\rm cf}(\delta)$,
to find a family ${\Cal P}_{\delta,\theta}$ of closed subsets of
$(\lambda^-_i,\delta)$ of order type $\theta$ such that any club of
$\delta$ contains (at least) one of them.  This holds by  
guessing clubs, see \cite[III,\S2]{Sh:g}.
\nl
2), 3) By the definitions. \hfill$\square_{\scite{1.4.4}}$\margincite{1.4.4}
\enddemo
\bigskip

\proclaim{\stag{1.8} Claim}   1) If $\mu > 2^\kappa$ (or just
$\lambda_* \ge 2^\kappa$) \ub{then} $\Bbb Q^2_\mu$ (hence 
$\Bbb Q^1_\mu$) has a dense subset of cardinality $\lambda_*$.
\nl
2) If $\mu > 2^\kappa$ (or just
$\lambda_* \ge 2^\kappa$) \ub{then} $\Bbb Q^3_\mu$ 
is a dense subset of $\Bbb Q^1_\mu$ and has cardinality $\lambda_*$.
\endproclaim
\bigskip

\demo{Proof}  1) By part (2).
\nl
2) By \scite{1.4.4}(2) it suffices to deal with $\Bbb Q^3_\mu$.
The cardinality of the set ${\Cal P}$ from \scite{1.4.4} is $\lambda_*$.
Whenever $p\in\Bbb Q^3_\mu$, then the function cd$_p$ is uniquely determined 
by its range, because 
$i\in$Dom$($cd$_p)$ iff Rang$($cd$_p)\cap[\lambda^-_i,\lambda_i)\ne\emptyset$ 
and the value cd$_p(i)=\alpha$ 
iff $\alpha\in[\lambda^-_i,\lambda_i)\cap$Rang$($cd$_p)$.
Also, the function cd$_p$ 
uniquely determines $p$ by 
$p=\bigcup\{C^*_{\text{cd}_p(i)}:i\in \text{Dom}(p)\}$.
Since Rang$($cd$_p)\subseteq u$, $u\in\Cal P$, we get
$|\Bbb Q^3_\mu|\le 2^\kappa\cdot\lambda_*=\lambda_*$. 
 \hphantom{a}\hfill$\square_{\scite{1.8}}$\margincite{1.8}
\enddemo
\bn
\relax From now on (till the end of this section)
\demo{\stag{1.5.13} Hypothesis}  $2^\kappa < \mu$ (in addition to \scite{1.1}).

Recall (Claim \scite{1.6.5}(1) is Balcar, Simon \cite[1.15]{BaSi89} and
\scite{1.6.5}(2) is a variant).
\enddemo
\bigskip

\definition{\stag{1.6} Definition}   A forcing notion $\Bbb P$ is
$(\theta,\lambda)$-nowhere distributive \ub{when} there are maximal
antichains $\bar p^\varepsilon = \langle p^\varepsilon_\alpha:\alpha <
\alpha_\varepsilon \rangle$ of $\Bbb P$ for $\varepsilon < \theta$ such
that for every $p \in \Bbb P$ for some $\varepsilon < \theta$,
we have $\lambda \le
|\{\alpha < \alpha_\varepsilon:p,p^\varepsilon_\alpha$ are compatible$\}|$.
\enddefinition
\bigskip

\proclaim{\stag{1.6.5} Claim}  1) If
\mr
\item "{$(a)$}"  $\Bbb P$ is a forcing notion,
$(\theta,\lambda)$-nowhere distributive 
\sn
\item "{$(b)$}"  $\Bbb P$ has density $\lambda$
\sn
\item "{$(c)$}"  $\theta > \aleph_0 \Rightarrow \Bbb P$ has a
$\theta$-complete dense subset
\ermn
\ub{then} $\Bbb P$ is equivalent to Levy$(\theta,\lambda)$.
\nl
2) If $\Bbb P$ is a forcing notion of density $\lambda$ collapsing
$\lambda$ to $\aleph_0$ \ub{then} $\Bbb P$ is equivalent to
Levy$(\aleph_0,\lambda)$.
\nl
3) If $\Bbb P$ is a forcing notion of density $\lambda$ and is 
$(\theta,\lambda)$-nowhere
distributive \ub{then} $\Bbb P$ collapses $\lambda$
to $\theta$ (and may or may not collapse $\theta$). 
  \hfill$\square_{\scite{1.6.5}}$\margincite{1.6.5}
\endproclaim
\bigskip

\proclaim{\stag{1.7.3} Claim}    Assume $\langle
b_\varepsilon:\varepsilon < \kappa \rangle$ is a sequence of pairwise
disjoint members of $[\kappa]^\kappa$ with union $b$.
\ub{Then} we can find an antichain ${\Cal I}$ of $\Bbb Q^3_\mu$ such that: 
\mr
\item "{$(*)$}"  if $q \in \Bbb Q^3_\mu$ and 
$(\forall \varepsilon < \kappa)(a(q) \cap
b_\varepsilon \in [\kappa]^\kappa)$, \ub{then} $q$ is compatible with
$\lambda_* =: \bold U_{J^{\text{bd}}_\kappa}(\mu)$ 
of the members of ${\Cal I}$.
\endroster
\endproclaim
\bigskip

\demo{Proof}  Let 

$$
\align
{\Cal I}^* = \{p \in \Bbb Q^3_\mu:& \text{ we can find an
increasing sequence } \langle i_\varepsilon:\varepsilon < \kappa \rangle\\
  & \text{ such that } i_\varepsilon \in b_\varepsilon \backslash
\varepsilon,  a(p) \subseteq \{i_\varepsilon:\varepsilon < \kappa\}\text{ and}\\ 
&  i_\varepsilon \in a(p) \Rightarrow p \cap
  [\lambda^-_{i_\varepsilon},\lambda_{i_\varepsilon})
  \text{ has order type } \lambda_\varepsilon\}.
\endalign
$$
\mn
Let ${\Cal J}^* = \{p \in \Bbb Q^3_\mu$: for every $\varepsilon <
\kappa$ we have $a(p) \cap b_\varepsilon \in [\kappa]^\kappa\}$.

Clearly
\mr
\item "{$(a)$}"  $|{\Cal I}^*| \le \lambda_* =
\bold U_{J^{\text{bd}}_\kappa}(\mu)$
\nl
[Why?  As ${\Cal I}^* \subseteq \Bbb Q^3_\mu$.]
\sn
\item "{$(b)$}"  if ${\Cal I} \subseteq {\Cal I}^*,|{\Cal I}| <
\lambda_*$ and $q
\in {\Cal J}^*$ then there is $r$ such that $q \le r \in {\Cal I}^*$
and $r$ is incompatible with every $p \in {\Cal I}$.
\ermn
[Why?  Let $\theta = |{\Cal I}| + \mu$, it is $< \lambda_*$, hence 
we can find an
increasing sequence $\langle \theta_\varepsilon:\varepsilon < \kappa
\rangle$ of regular cardinals with limit $\mu$ such that $\dsize
\prod_{\varepsilon < \kappa} \theta_\varepsilon/J^{\text{bd}}_\kappa$
has true cofinality $\theta^+$, this by \scite{1.5} + the no hole
lemma \cite[II,\S3]{Sh:g}.   By renaming \wilog \,
$\theta_\varepsilon > \lambda_\varepsilon$.

Let $u = \{\varepsilon < \kappa:a(q) \cap b_\varepsilon \in [\kappa]^\kappa\}$,
so we know that $u$ is $\kappa$.  For each $\varepsilon
\in u$ we know that $a(q) \cap b_\varepsilon \in [\kappa]^\kappa$, and
so for some $\zeta_\varepsilon < \kappa$ we have $\theta_\varepsilon <
\lambda_{\text{otp}(a(q) \cap \zeta_\varepsilon)}$.  Now
choose $i(\varepsilon) \in b_\varepsilon$ such that
$i(\varepsilon) > \varepsilon \wedge i(\varepsilon) >
\zeta_\varepsilon \wedge (\forall \varepsilon_1 <
\varepsilon)(i(\varepsilon_1) < i(\varepsilon))$.  
As $q \in \Bbb Q^3_\mu$ it follows that
$(q \cap [\lambda^-_{i(\varepsilon)},\lambda_{i(\varepsilon)}))$ 
has order type $\ge \lambda_{\text{otp}(a(q)\cap \zeta_\varepsilon)} 
> \theta_\varepsilon$.
Let $C_{q,\varepsilon} = \{\alpha:\alpha \in
q,\alpha \in [\lambda^-_{i(\varepsilon)},\lambda_{i(\varepsilon)})$
and otp$(q \cap [\lambda^-_{i(\varepsilon)},\lambda_{i(\varepsilon)})
\cap \alpha)$ is $< \theta_\varepsilon\}$.  Now for every $p \in
{\Cal I}^*$ the set 
$p \cap [\lambda^-_{i(\varepsilon)},\lambda_{i(\varepsilon)})
\subseteq \cup\{[\lambda^-_i,\lambda_i):i \in b_\varepsilon\}$ if
non-empty has
cardinality $\le \lambda_\varepsilon$ which is $< \theta_\varepsilon$
hence $p \cap C_{q,\varepsilon}$ is
a bounded subset of $C_{q,\varepsilon}$, call the lub
$\alpha_{p,\varepsilon}$.  As $\theta = |{\Cal I}| + \mu < \text{\rm
tcf}(\dsize \prod_{\varepsilon < \kappa}
\theta_\varepsilon/J^{\text{bd}}_\kappa)$ clearly there is $h \in \dsize
\prod_{\varepsilon \in u} C_{q,\varepsilon}$ such that $p \in {\Cal I}^*
\Rightarrow \langle \alpha_{p,\varepsilon}:\varepsilon < \kappa \rangle
<_{J^{\text{bd}}_\kappa} h$ and let

$$
\align
r = \{\alpha:&\text{for some } \varepsilon \in u \text{ we have }
\alpha \in C_{q,\varepsilon} \backslash h(\varepsilon) \\
  &\text{ and } |C_{q,\varepsilon} \cap \alpha \backslash h(\varepsilon)|
  <\lambda_\varepsilon\}.
\endalign
$$
\mn
So $r$ is as required in clause (b).  (We can assume that $r \in
\Bbb Q^3_\mu$,  since by the density propositions of \scite{1.8} we can find $r
\le r' \in \Bbb Q^3_\mu$ as required.)  So clause (b) holds.] 
\nl
As by \scite{1.8}(2) in the 
conclusion of the claim it is enough to deal with $q \in
\Bbb Q^3_\mu$, there are only $\lambda_*$ such $q$'s so
we can finish easily by (clause (b) and) diagonalization.  
\hfill$\square_{\scite{1.7.3}}$\margincite{1.7.3}
\enddemo
\bigskip

\proclaim{\stag{1.7.6} Claim}  The forcing notion $\Bbb Q^3_\mu$ is
$({\frak b}_\kappa,\lambda_*)$-nowhere distributive.
\endproclaim
\bigskip

\demo{Proof}  Let $\langle \bar A_\alpha:\alpha < {\frak
b}_\kappa\rangle$ be such that: $\bar A_\alpha = \langle
 A_{\alpha,i}:i < \kappa \rangle,A_{\alpha,i} \in [\kappa]^\kappa,i < j
 \Rightarrow A_{\alpha,i} \cap A_{\alpha,j} = \emptyset$ and $(\forall
 B \in [\kappa]^\kappa)(\exists \alpha < {\frak b}_\kappa)(\forall i <
 \kappa)[\kappa = |B \cap A_{\alpha,i}|]$, exists by \scite{4b.4}(2) below.
 Hence for each $\alpha < {\frak b}_\kappa,{\Cal I}^*_\alpha 
\subseteq \Bbb Q^3_\mu$
as in \scite{1.7.3} for the sequence $\bar A_\alpha$ exists. 
So $\langle {\Cal I}^*_\alpha:\alpha < {\frak b}_\kappa \rangle$ is a
sequence of ${\frak b}_\kappa$ antichains of $\Bbb Q^3_\mu$ and we
shall show that it witnesses the conclusion.  Now
\mr
\item "{$\circledast$}"     if $q \in \Bbb Q^3_\mu$ then for some 
$\alpha < {\frak b}_\kappa$  the set $\{p \in{\Cal I}^*_\alpha:p$
compatible with $q \in \Bbb Q^3_\mu\}$ has cardinality $\lambda_*$.
\ermn
Why?  By the choice of $\langle \bar A_\alpha:\alpha < {\frak
b}_\kappa \rangle$ there is $\alpha < {\frak b}_\kappa$ such that
\mr
\item "{$(*)$}"   $a(q) \cap A_{\alpha,i} \in [\kappa]^\kappa$ for 
every $i < \kappa$.
\ermn
Hence $q$ fits the 
demand in \scite{1.7.3} with $\bar A_\alpha$ here standing for
$\langle b_\varepsilon:\varepsilon < \kappa \rangle$.  Hence it is
compatible with $\lambda_*$ members of ${\Cal I}^*_\alpha$ which, of course,
shows that we are done.
\nl
${{}}$   \hfill$\square_{\scite{1.7.6}}$\margincite{1.7.6}
\enddemo
\bigskip

\demo{\stag{1.9} Conclusion}    1) If $2^\kappa < \mu$ (and $\aleph_0 < 
\kappa = \text{\rm cf}(\mu) < \mu$, of course) 
\ub{then} $\Bbb Q_\mu$ is equivalent to
Levy$(\aleph_0,\lambda_*)$, i.e., they have isomorphic completions
(recalling $\Bbb Q_\mu$ is naturally completely embeddable into the completion
of $\Bbb P_\mu = ([\mu]^\mu,\supseteq))$.
\nl
2) If $(\forall \alpha < \mu)(|\alpha|^\kappa < \mu)$ \ub{then} $\Bbb Q_\mu$ is
equivalent to Levy$(\aleph_0,\mu^\kappa)$.
\nl
3) If $\mu$ is strong limit (singular of uncountable cofinality
$\kappa$), \ub{then} $\Bbb P_\mu$ is equivalent to
Levy$(\aleph_0,\mu^\kappa) = \text{\rm Levy}(\aleph_0,2^\mu)$.
\enddemo
\bigskip

\demo{Proof}  1) By \scite{1.8}(1), $\Bbb Q^3_\mu$ has density
(even cardinality) $\lambda_*$ and by \scite{1.7.6} it is 
$({\frak b}_\kappa,\lambda_*)$-nowhere distributive 
hence by \scite{1.6.5}(3), we know
that $\Bbb Q^3_\mu$ collapses $\lambda_*$ to ${\frak b}_\kappa$.
But $\Bbb P_\kappa$ is completely embeddable into $\Bbb Q^2_\mu$ (see
\scite{1.3.3}(3)) and $\Bbb P_\kappa$ collapses ${\frak b}_\kappa$ to
$\aleph_0$ (e.g. see \S2) and $\Bbb Q^3_\mu$ is dense in $\Bbb Q^2_\mu$.  
Together forcing with $\Bbb Q^3_\mu$ collapses $\lambda_*$
to $\aleph_0$.  As $\Bbb Q^3_\mu$ has density $\lambda_*$, by
\scite{1.6.5}(2) we get that $\Bbb Q^2_\mu$ is equivalent to
Levy$(\aleph_0,\lambda_*)$. 

Lastly $\Bbb Q_\mu,\Bbb Q^3_\mu$ are equivalent by \scite{1.3}(3) +
\scite{1.4.4}(2) so we are done.
\nl
2) Recalling \scite{1.5}, by \cite[VIII]{Sh:g} we have 
$\lambda_* = \mu^\kappa$ (alternatively directly as in
\cite[\S3]{Sh:506}).  Now apply part (1).
\nl
3) By easy cardinal arithmetic $\mu^\kappa = 2^\mu$.  Enough to check
the demands in \scite{1.6.5}(2).  Now as $\Bbb Q_\mu$ collapses
$\lambda_*$ to $\aleph_0$ by part (1) and $\Bbb Q_\mu$ can be
completely embeddable into $\Bbb P_\mu$ (see \scite{1.3.3}(2)) clearly
$\Bbb P_\mu$ collapses $\lambda_*$ to $\aleph_0$.  But $|\Bbb P_\mu| \le
|[\mu]^\mu| = 2^\mu$, so $\Bbb P_\mu$ has density $\le 2^\mu$.

Lastly $\lambda_* = 2^\mu$ by \cite[VIII]{Sh:g}.  
So we are done.  \hfill$\square_{\scite{1.9}}$\margincite{1.9}
\enddemo
\bigskip

\proclaim{\stag{1.11} Claim}  Assume that $\Bbb P_\mu$ does not satisfy
the $\chi$-c.c.  \ub{Then} 
forcing with $\Bbb P_\mu$ collapses $\chi$ to $\aleph_0$.
\endproclaim
\bigskip

\demo{Proof}  By the nature of the conclusion \wilog \, $\chi$ is
regular.  Now we can find $\bar X$ such that
\mr
\item "{$(*)_1$}"  $(a) \quad \bar X = \langle X_\xi:\xi <
\chi \rangle$
\sn
\item "{${{}}$}"  $(b) \quad X_\xi \in \Bbb P_\mu$
\sn
\item "{${{}}$}"  $(c) \quad X_\zeta \cap X_\xi \in [\mu]^{< \mu}$ for
$\zeta \ne \xi < \chi$.
\ermn
As $\Bbb Q_\mu \lessdot \Bbb P_\mu$, by the earlier proof (e.g.,
\scite{1.9}(1)) it suffices to prove
\comment
that $\Bbb P_\mu$ collapses $\chi$ to $\lambda_*$.  Let $\bold P =
\{\bar A:\bar A=\langle A_\alpha:\alpha < \mu \rangle$, the
$A_\alpha$'s are pairwise disjoint and each $A_\alpha$ belongs to
$[\mu]^\mu\}$.  Now
\endcomment
that $\Bbb P_\mu$ collapses $\chi$ to $\lambda_*$.  
There exists $\bold P \subseteq {\bold P}_*:=
\{\bar A:\bar A=\langle A_\alpha:\alpha < \mu \rangle$, the
$A_\alpha$'s are pairwise disjoint and each $A_\alpha$ belongs to
$[\mu]^\mu\}$ such that $|\bold P|=\lambda_*$ and
\mr
\item "{$(*)_2$}"  for every $p \in \Bbb P_\mu$ 
there is an $\bar A \in \bold P$ such that 
$(\forall \alpha < \mu)[|A_\alpha \cap p| = \mu]$.
\ermn
[Why? For each $i<\kappa$ fix some partition 
$\langle W_{i,\alpha}:\alpha<\lambda_i\rangle$ 
of $\lambda_i$ into $\lambda_i$ (pairwise disjoint) 
sets each of cardinality $\lambda_i$.
Now for each $p\in {\Bbb P}_\mu$ we shall choose 
$\bar A=\bar A^p \in {\bold P}_\lambda$ as required in 
$(*)_2$ such that
${\bold P}:=\{\bar A^p:p\in {\Bbb P}_\mu\}$ has cardinality 
$\leq \lambda_*$ this suffice; so fix $p\in {\Bbb P}_\mu$.
By induction on $\varepsilon < \kappa$ we can find
$\delta_\varepsilon < \mu$ of cofinality $\lambda^{++}_\varepsilon$
such that $p \cap \delta_\varepsilon$ is unbounded in
$\delta_\varepsilon$ and $\delta_\varepsilon >
\cup\{\delta_\zeta:\zeta < \varepsilon\}$.  
There is a club $C^1_\varepsilon$
of $\delta_\varepsilon$ of order type $\lambda^{++}_\varepsilon$ with
min$(C^1_\varepsilon) > \cup\{\delta_\zeta:\zeta < \varepsilon\}$.
Let $C^2_\varepsilon = \{\delta \in C^1_\varepsilon: 
\delta$ is a limit ordinal such that $C^1_\epsilon\cap p$
is unbounded in $\delta$ and has order type 
divisible by $\lambda^+_\epsilon\}$, 
it is a club of $\delta_\varepsilon$.  
But by the club guessing (see
\scite{1.4.4}) there is
$C^3_\varepsilon$ such that: $C^3_\varepsilon \subseteq
C^2_\varepsilon (\subseteq C^1_\varepsilon)$ and otp$(C^3_\varepsilon)
= \lambda_\varepsilon$.  

By the definition of $\Bbb Q^3_\mu$, there is some $a \in
[\kappa]^\kappa$ such that $\bigcup\{
C^3_\varepsilon:\varepsilon \in a\}\in\Bbb Q^3_\mu$.  Lastly, let
us define $\bar A = \langle A_\alpha:\alpha < \mu \rangle$ by

$$
\align
A_\alpha = \cup\{[\beta,\text{min}(C^3_\varepsilon \backslash (\beta
+1)):&\varepsilon \in a \text{ satisfies} \\
  &\alpha < \lambda_\varepsilon \text{ and } \beta \in C^3_\varepsilon
\text{ and} \\
  &\text{ otp}(C^3_\varepsilon \cap \beta) \in
  W_{\varepsilon,\alpha}\}.
\endalign
$$
\mn
Easily $\langle A_\alpha:\alpha < \mu \rangle$ is as
required in $(*)_2$, and since $\bar A$ is determined by an 
element of $\Bbb Q^3_\mu$ (and the constant 
$\langle W_{i,\alpha}:\alpha<\lambda_i:i<\kappa\rangle$),
the cardinality $|\bold P|\le|\Bbb Q^3_\mu|\le\lambda_*$.] 
\nl
Now for $\bar A \in \bold P$ we define a $\Bbb
P_\mu$-name ${\underset\tilde {}\to \tau_{\bar A}}$ as follows: for $\bold G
\subseteq \Bbb P_\mu$ generic over $\bold V$,
\mr
\item "{$(*)_3$}"  ${\underset\tilde {}\to \tau_{\bar A}}[\bold G] =
\xi$ iff $\xi$ is minimal such that 
$\cup\{A_\alpha:\alpha \in X_\xi\} \in \bold G$
\ermn
clearly
\mr
\item "{$(*)_4$}"  ${\underset\tilde {}\to \tau_{\bar A}}[\bold G]$ is
defined in at most one way;
\sn
\item "{$(*)_5$}"  for every $p \in \Bbb P_\mu$ for some $\bar A \in
\bold P$ for every $\xi < \chi$ we have $p \nVdash
``{\underset\tilde {}\to \tau_{\bar A}} \ne \xi"$.
\ermn
[Why?  Let $\bar A \in \bold P$ be such that $(\forall
\alpha < \mu)(\mu = |p \cap A_\alpha|)$, it exists by $(*)_2$.
Now we can find $q$ satisfying 
$p \le q \in \Bbb P_\mu$ such that $(\forall \alpha
< \mu)(q \cap A_\alpha$ is a singleton) and for each $\xi < \chi$ let
$q_\xi = \cup\{A_\alpha \cap q:\alpha \in X_\xi\}$.  Clearly $\zeta <
\xi \Rightarrow |X_\zeta \cap X_\xi| < \mu \Rightarrow \cup
\{A_\alpha:\alpha \in X_\zeta\} \cap q_\xi \subseteq \cup\{A_\alpha
\cap q_\xi:\alpha \in X_\zeta\} = \cup\{A_\alpha \cap q:\alpha \in
X_\zeta \cap X_\xi\} \in [\mu]^{< \mu}$, hence $q_\xi \Vdash ``\xi =
{\underset\tilde {}\to \tau_{\bar A}}"$.]
\nl
So
\mr
\item "{$(*)_6$}"  $\Vdash_{\Bbb P_\mu} ``\chi = \{{\underset\tilde {}\to
\tau_{\bar A}}[{\underset\tilde {}\to {\bold G}}]:\bar A \in \bold P\}"$.
\ermn
Together clearly $\Bbb P_\mu$ collapses $\chi$ to $\lambda_* + |\bold P|$ 
which is $\le \lambda_*$, so as said
above we are done.  \hfill$\square_{\scite{1.11}}$\margincite{1.11}
\enddemo
\bn
Lastly, concerning the singular $\mu_*$ of cofinality $\aleph_0$ so we forget
the hypothesis \scite{1.1}, \scite{1.5.13}.
\proclaim{\stag{1.12} Claim}  If $\mu_* > \text{\rm cf}(\mu_*) =
\aleph_0$ and $\Bbb P_{\mu_*}$ fails the $\chi$-c.c., \ub{then} $\Bbb
P_{\mu_*}$ collapses $\chi$ to $\aleph_1$; note that in this case $\Bbb
Q_{\mu_*}$ is equivalent to 
Levy$(\aleph_1,\mu^{\aleph_0}_*)$ by \cite{KjSh:720}.
\endproclaim
\bigskip

\demo{Proof}  Let $\lambda_* = \mu^{\aleph_0}_*$.

By Kojman, Shelah \cite{KjSh:720}, $\Bbb P_{\mu_*}$ collapses
$\lambda_*$ to $\aleph_1$ hence it suffices to prove that $\Bbb P_\mu$
collapse $\chi$ to $\lambda_*$ assuming $\chi > \lambda_*$ (otherwise
the conclusion is known).  Let $\langle \lambda_n:n <
\omega \rangle$ be a sequence of regular uncountable cardinals with limit
$\mu_*$.  
Now repeat the proof of \scite{1.11}  \hfill$\square_{\scite{1.12}}$\margincite{1.12}  
\enddemo
\newpage

\head {\S2 The regular uncountable case} \endhead  \resetall \sectno=2
 \spuriousreset
\bigskip

We prove that (for $\kappa$ regular uncountable), $\Bbb P_\kappa$
collapse $\lambda$ to $\aleph_0$ iff $\Bbb P_\kappa$ fail the $\lambda$-c.c.
This continues Balcar, Simon \cite[2.8]{BaSi88} so we first re-represent
what they do; the proof of \scite{4b.3} is made to help later.  
In the present notation they let $\lambda = 
{\frak b}_\kappa$ (rather that $\lambda\in {\frak b}^{\text{{\rm spc}}} 
\kappa$ as below, 
let $\langle f_\alpha:\alpha < {\frak b}_\kappa \rangle$ be
a sequence exemplifying it; let $C_\alpha = \{\delta <
\kappa:(\forall \beta < \delta)(f_\alpha(\beta) < \delta),\delta$ a 
limit ordinal$\}$ and let
$B_\alpha = \kappa \backslash C_\alpha$, so $\langle B_\alpha:\alpha
< \lambda \rangle$ is a $(\kappa,\lambda)$-sequence (see
\scite{4a.2}(1)), derive a good $(\kappa,{}^{\omega >}\lambda)$-sequence
from it (see \scite{4a.2}(2)), define $\alpha_n(A),\beta_n(A)$ and used the
$A_{\eta,\delta,i}$'s to define the $\Bbb P_\kappa$-names
${\underset\tilde {}\to \beta_n}$ and prove $\Vdash_{\Bbb P_\kappa}
``\{g^*({\underset\tilde {}\to \beta_n}):n < \omega\} = {\frak
b}_\kappa"$ (see \scite{4b.3}).  We then prove the new result: if $\Bbb
P_\kappa$ fail the $\chi$-c.c. then it collapses $\chi$ to $\aleph_0$.
\bigskip

\demo{\stag{4.0} Context}  $\kappa$ is a fixed regular uncountable cardinal.
\enddemo
\bigskip

\definition{\stag{4.1} Definition}   1) Let ${\frak
b}^{\text{spc}}_\kappa$ be the set of regular
  $\lambda > \kappa$ such that there is a
  $<_{J^{\text{bd}}_\kappa}$-increasing sequence $\langle
  f_\alpha:\alpha < \lambda \rangle$ of members of ${}^\kappa \kappa$ with no
$\le_{J^{\text{bd}}_\kappa}$-upper bound in ${}^\kappa \kappa$.
\nl
2) Let ${\frak b}_\kappa = \text{\rm Min}({\frak b}^{\text{spc}}_\kappa)$.
\enddefinition
\bigskip

\definition{\stag{4a.1} Definition}  1) We say $\bar B$ is a
$(\kappa,\lambda)$-sequence when
\mr
\item "{$(a)$}"  $\bar B = \langle B_\alpha:\alpha < \lambda \rangle$
\sn
\item "{$(b)$}"  $B_\alpha \in [\kappa]^\kappa$ and $\kappa \backslash
B_\alpha \in [\kappa]^\kappa$ and $B_{\alpha +1} \backslash B_\alpha
\in [\kappa]^\kappa$
\sn
\item "{$(c)$}"  for every $B \in [\kappa]^\kappa$ for some
$\alpha,B \cap B_\alpha \in [\kappa]^\kappa$
\sn
\item "{$(d)$}"  $B_\alpha \subseteq^* B_\beta$ when $\alpha < \beta <
\lambda$, i.e., $B_\alpha \backslash B_\beta \in [\kappa]^{< \kappa}$.
\ermn
2) We say that $\bar B$ is a $(\kappa,{}^{\omega >}\lambda)$-sequence when:
\mr
\item "{$(a)$}"  $\bar B = \langle B_\eta:\eta \in {}^{\omega >}
\lambda \rangle$
\sn
\item "{$(b)$}"  $B_\eta \in [\kappa]^\kappa$
\sn
\item "{$(c)$}"  if $\eta_1 \triangleleft \eta_2 \in {}^{\omega >}\lambda$
then $B_{\eta_2} \subseteq^* B_{\eta_1}$ which means $B_{\eta_2}
\backslash B_{\eta_1} \in [\kappa]^{< \kappa}$
\sn
\item "{$(d)$}"  $B_{<>} = \kappa$
\sn
\item "{$(e)$}"  if $\eta \in {}^{\omega >}\lambda$ and $A \in
[B_\eta]^\kappa$ then for some $\alpha < \lambda$ we have $A \cap
B_{\eta {}^\frown <\alpha>} \in [\kappa]^\kappa$
\sn
\item "{$(f)$}"  if $\eta \in {}^{\omega >}\lambda$ and $\alpha <
\beta < \lambda$ then $B_{\eta {}^\frown<\alpha>} \subseteq^* 
B_{\eta {}^\frown <\beta>}$  and $B_\eta \backslash B_{\eta {}^\frown <\alpha>}
\in [\kappa]^\kappa$ and $B_{\eta {}^\frown<\beta>} \backslash B_{\eta
{}^\frown<\alpha>} \in [\kappa]^\kappa$.
\ermn
3) For a $(\kappa,{}^{\omega >}\lambda)$-sequence $\bar B$ and $A \in
[\kappa]^\kappa$ we try to define an ordinal $\alpha_k(A,\bar B)$ by
induction on $k< \omega$.  If $\eta = \langle \alpha_\ell(A,\bar
B):\ell < k\rangle$ is well defined (holds for $k=0$) and there is an
$\alpha <  \lambda$ such that $A \subseteq^* B_{\eta {}^\frown
<\alpha>} \wedge (\forall \beta < \alpha)(A \cap B_{\eta {}^\frown
<\beta>} \in [\kappa]^{< \kappa})$ then we let $\alpha_k(A,\bar B) =
\alpha$; note that $\alpha$, if exists, is unique.  Let $n(A,\bar
B)$ be the $n \le \omega$ such that $\alpha_\ell(A,\bar B)$ is well
defined iff $\ell < n$.
\nl
4) We say that $(\bar B,\bar \nu)$ is a 
$(\kappa,{}^{\omega >}\lambda)$-parameter when:
\mr
\item "{$(a)$}"  $\bar B = \langle B_\eta:\eta \in {}^{\omega
>}\lambda\rangle$ is a $(\kappa,{}^{\omega >}\lambda)$-sequence
\sn
\item "{$(b)$}"  $\bar \nu$ is an $S^\lambda_\kappa$-ladder which
means that $\bar \nu = \langle \nu_\delta:\delta \in
S^\lambda_\kappa \rangle,\nu_\delta$ is an increasing sequence of
ordinals of length $\kappa$ with limit $\delta$, where
$S^\lambda_\kappa = \{\delta < \lambda:\text{\rm cf}(\delta)=\kappa\}$.
\ermn
5) We say $(\bar B,\bar \nu)$ is a 
good $(\kappa,{}^{\omega >}\lambda)$-parameter when (a)+(b) of part
(4) holds and
\mr
\item "{$(c)$}"  if $A \in [\kappa]^\kappa$ then for some $n <
\omega,\eta \in {}^n\lambda$ and $\delta \in S^\lambda_\kappa$ and $A'
\in [A]^\kappa$ we have
{\roster
\itemitem{$ (\alpha)$ }  $\alpha_\ell(A',\bar B) = \eta(\ell)$ for $\ell <
n$
\sn
\itemitem{ $(\beta)$ }  for $\kappa$ many ordinals $\zeta < \kappa$ we have
$(\forall \varepsilon < \zeta)(A' \cap
B_{\eta {}^\frown <\nu_\delta(\zeta)>} \backslash 
B_{\eta {}^\frown <\nu_\delta(\varepsilon)>}$ belongs to
$[\kappa]^\kappa)$.
\endroster}
\ermn
6) $\bar B$ is a good $(\kappa,{}^{\omega >}\lambda)$-sequence if clause (a) of
(4) and clause (c) of (5) holds for some $S^\lambda_\kappa$-ladder
(see above).  We say $\bar B$ is a weakly good sequence if clause (a) of
(4) and clause (c)$^-$ of (5) which means that we ignore subclause
$(\alpha)$ there.  Similarly $(\bar B,\bar \nu)$ is a weakly good
$(\kappa,{}^{\omega >}\lambda)$-parameter.
\enddefinition
\bigskip

\demo{\stag{2.3.1} Observation}  
1)\  In \scite{4a.1}(5)(c)$(\beta)$, the ``for $\kappa$ many ordinal 
$\zeta<\kappa$'' implies ``for club many ordinals $\zeta<\kappa_0$. 
\nl
2)\  In \scite{4a.1}(6) it doesn't matter 
which $S^\lambda_\kappa$-ladder you choose.
\enddemo

\demo{Proof}
If $\nu_1,\nu_2 \in {}^\kappa \delta$ are
increasing and sup$(\nu_1) =
\sup(\nu_2) = \delta$, then $\{i < \kappa:\dbcu_{j<i} \nu_1(j) =
\dbcu_{j<i} \nu_2(j)\}$ is a club of $\kappa$.  
\enddemo
\bn

Note that for \S1 we need no more than Claim \scite{4a.2} (actually
the weakly good version is enough for \S1 except presenting the proof that
${\frak b}_\kappa$ is collapsed).

\proclaim{\stag{4a.2} Claim}  1) Assume $\lambda = {\frak b}_\kappa$ or
just $\lambda \in {\frak b}^{\text{spc}}_\kappa$.  \ub{Then}
$\lambda$ is regular $> \kappa$ and there is a $\subseteq^*$-decreasing
sequence $\langle C_\alpha:\alpha < \lambda \rangle$ of clubs
of $\kappa$ such that for no $A \in [\kappa]^\kappa$ do we have $\alpha
< \lambda \Rightarrow A \subseteq^* C_\alpha$.  Hence $\langle \kappa
\backslash C_\alpha:\alpha < \lambda \rangle$ is a $(\kappa,\lambda)$-sequence.
\nl
2) Assume $\bar C = \langle C_\alpha:\alpha < \lambda \rangle$ is as 
above and $\bar \nu = \langle \nu_\delta:\delta \in S^\lambda_\kappa \rangle$
is an $S^\lambda_\kappa$-ladder, see 
Definition \scite{4a.1}(4), clause (b) (such $\bar \nu$
always exists).  \ub{Then} $\bar B = \bar B_{\bar C},\bar f = \bar f_{\bar
C}$ are well defined and the pair 
$(\bar B,\bar \nu)$ is a good $(\kappa,{}^{\omega
>}\lambda)$-parameter where we define $\bar B$ and $\bar f$ as follows:
\mr
\item "{$\circledast$}"  $(a) \quad \bar B = \langle B_\eta:\eta \in {}^{\omega
>}\lambda \rangle$
\sn
\item "{${{}}$}"  $(b) \quad \bar f = \langle f_\eta:\eta \in {}^{\omega
>}\lambda \rangle$
\sn
\item "{${{}}$}"  $(c) \quad B_{<>} = \kappa,
f_{<>} = \text{\rm id}_\kappa$
\sn
\item "{${{}}$}"  $(d) \quad B_\eta \in [\kappa]^\kappa,
f_\eta$ is a function from $B_\eta$ onto 
$\kappa$, non-decreasing, 
\nl

\hskip25pt and not eventually constant
\sn

\item "{${{}}$}"  $(e) \quad$ if the 
pair $(B_\rho,f_\rho)$ is defined and $\alpha < \lambda$ then we let

$$
B_{\rho{}^\frown <\alpha>} = \{\gamma \in B_\rho:f_\rho(\gamma) \in
\kappa \backslash C_\alpha\}
$$
\sn
\item "{${{}}$}"  $(f) \quad$ if 
$\eta = \rho {}^\frown \langle \alpha \rangle$ and
$B_\rho,f_\rho$ and $B_\eta$ are defined then we let 
\nl

\hskip25pt $f_\eta:B_\eta \rightarrow \kappa$ be defined by
$f_\eta(i)=\text{\rm otp}(C_\alpha\cap f_\rho(i))$
\sn
\hskip25pt for each $i < \kappa$,
\sn
\
\hskip25pt hence
\sn
\item "{${{}}$}"  $(g) \quad$ if $\eta{}^\frown \langle \alpha \rangle
\in {}^{\omega >}\lambda$ then $B_{\eta {}^\frown<\alpha>} \subseteq
B_\eta$ and $i \in B_{\eta {}^\frown <\alpha>} 
\wedge f_\eta(i) > 0 \Rightarrow$
\nl

\hskip25pt $f_\eta(i) > f_{\eta {}^\frown<\alpha>}(i)$.  
\endroster          
\endproclaim
\bigskip

\demo{Proof}  1) Recall $S^\lambda_\kappa := \{\delta < \lambda:\text{\rm
cf}(\delta) = \kappa\}$.

By the definition of ${\frak b}^{\text{spc}}_\kappa$ there is an
$<_{J^{\text{bd}}_\kappa}$-increasing sequence $\langle
f^*_\alpha:\alpha < \lambda \rangle$ of members of ${}^\kappa \kappa$
with no $\le_{J^{\text{bd}}_\kappa}$-upper bound from ${}^\kappa
\kappa$.  Let $C_\alpha := \{\delta < \kappa:\delta$ is a limit
ordinal such that $(\forall \gamma < \delta)(f^*_\alpha(\gamma) <
\delta)\}$.

Clearly
\mr
\item "{$(*)_1$}"  $C_\alpha$ is a club of $\kappa$
\nl
[why?  as $\kappa$ is regular uncountable]
\sn
\item "{$(*)_2$}"   if $\alpha < \beta < \lambda$ then $C_\beta
\subseteq^* C_\alpha$; i.e., $C_\beta \backslash C_\alpha \in [\kappa]^{<
\kappa}$
\nl
[why?  as if $\alpha < \beta$ then $f^*_\alpha <_{J^{\text{bd}}_\kappa}
f^*_\beta$, i.e., for some $\varepsilon < \kappa,(\forall
\zeta)(\varepsilon \le \zeta < \kappa \Rightarrow f^*_\alpha(\zeta) <
f^*_\beta(\zeta))$ hence 
letting $\epsilon_1=\text{sup} (\text{Rang} f^*_\alpha
\restriction \alpha)$, we have
$C_\beta \backslash (\varepsilon_1+1) \subseteq
C_\alpha$ as required] 
\sn
\item "{$(*)_3$}"  for every club $C$ of $\kappa$ for some $\zeta <
\lambda$ we have $C \backslash C_\zeta \in [\kappa]^\kappa$
\nl
[why?  as $\bar f$ has no $\le_{J^{\text{bd}}_\kappa}$-bound in
${}^\kappa \kappa$] 
\ermn
hence
\mr
\item "{$(*)_4$}"  for every unbounded subset $A$ of $\kappa$ for
some $\zeta < \lambda$ we have $A \backslash C_\zeta \in
[\kappa]^\kappa$.
\nl
[Why?  Otherwise the closure of $A$ contradicts $(*)_3$.]
\ermn
Clearly $\langle C_\alpha:\alpha < \lambda \rangle$ is as required.

Lastly, let $B_\alpha = \kappa \backslash C_\alpha$, it is easy to
check that $\langle B_\alpha:\alpha < \lambda \rangle$ is a
$(\kappa,\lambda)$-sequence. 
\medskip
2) Clearly $\bar B_{\bar C},\bar f_{\bar C}$ are well defined and
$(\bar B,\bar \nu)$ is a 
$(\kappa,{}^{\omega >}\lambda)$-parameter and clauses (a)-(g) of
$\circledast$ holds.  Why is it good?  Toward
contradiction assume that it is not, so choose 
$A \in [\kappa]^\kappa$ which exemplify the failure of clause (c) of
Definition \scite{4a.1}(5) and define 

$$
\align
{\Cal T}_0 = {\Cal T}^0_A = \bigl\{ \eta \in {}^{\omega
>}\lambda:&\text{ there is } A' \in [A]^\kappa \text{ such that} \\
  &\langle \alpha_\ell(A',\bar B):\ell < \ell g
(\eta)\rangle \text{ is well defined and equal to } \eta \bigr\}.
\endalign
$$
\mn
and define

$$
\align
{\Cal T}_1 = {\Cal T}^1_A := \bigl\{ \eta \in {\Cal T}^0_A:&\text{ for
  every } k < \ell g(\eta) \text{ there are } < \kappa \\
  &\text{ordinals } \alpha < \eta(k) \text{ such that } (\eta\restriction k)^\frown
  \langle \alpha \rangle \in {\Cal T}_0 \bigr\}.
\endalign
$$
\mn
Clearly
\mr
\item "{$(*)_1$}"  ${\Cal T}_0 \supseteq {\Cal T}_1$ are 
non-empty subsets of ${}^{\omega >}\lambda$ (in fact $<> \in {\Cal
T}_1 \subseteq {\Cal T}_0$)
\sn
\item "{$(*)_2$}"  ${\Cal T}_0,{\Cal T}_1$ are closed under initial segments.
\ermn
For $\eta \in {\Cal T}_\ell$ let Suc$_{{\Cal T}_\ell}(\eta) 
= \{\rho \in {\Cal T}_\ell:\ell g(\rho) = 
\ell g(\eta)+1$ and $\eta \triangleleft \rho\}$.

We define $A_\eta \in [B_\eta]^\kappa$ for $\eta \in {\Cal T}_1$ by
induction on $\ell g(\eta)$:
\mr
\item "{$(*)_3$}"  $(a) \quad A_{<>} = A$
\sn
\item "{${{}}$}"  $(b) \quad$ if $A_\nu$ is defined and $\nu {}^\frown \langle\alpha\rangle
\in {\Cal T}_1$ then we let 
\nl

\hskip25pt $A_{\nu {}^\frown <\alpha>} = A_\nu \cap B_{\nu
{}^\frown <\alpha>} \backslash \bigcup\{B_{\nu {}^\frown <\beta>}:\beta <
\alpha$ and $\nu {}^\frown \langle \beta \rangle \in {\Cal T}_1\}$.
\ermn
Now
\mr
\item "{$(*)_4$}"  if $\nu \in {\Cal T}_1$ then 
{\roster
\itemitem{ $(a)$ }  if $B \in [A]^\kappa$ and $\langle
\alpha_\ell(B,\bar B):\ell < \ell g(\nu)\rangle$ 
is well defined and equal to $\nu$ \ub{then} $B \subseteq^* A_\nu$
\sn
\itemitem{ $(b)$ }  if $\text{{\rm Suc}}_{{\Cal T}_j} (\nu)$ 
has cardinality $<\kappa$ then $A_\nu \backslash
\cup\{A_\rho:\rho \in \text{\rm Suc}_{{\Cal T}_j}(\nu)\}$ 
has cardinality $< \kappa$ for $j=1$ (actually $j=0$ is O.K., too).
\sn
\itemitem{ $(c)$ } If $\text{{\rm Suc}}_{{\Cal T}_1} (\nu)$ has 
cardinality $<\kappa$ then $\text{{\rm Suc}}_{{\Cal T}_0}(\nu)=
\text{{\rm Suc}}_{{\Cal T}_1} (\nu)$
\endroster}
\ermn
[Why?  First we can prove clause (a) by induction on $\ell g(\nu)$
using the definition of ${\Cal T}_1$ and clause (c) of
\scite{4a.1}(2).  Second, we can prove clause 
(b) from it. Third why clause (c) holds? 
\newline
 Otherwise, as ${\Cal T}_1 \subseteq {\Cal T}_0$, there is an $\alpha$ with
 $\nu_n{}^\frown
\langle \alpha \rangle \in \text{ Suc}_{{\Cal T}_0}(\nu_n)
\backslash \text{ Suc}_{{\Cal T}_1}(\nu_n)$. Hence by the definition
of ${\Cal T}_1$ the set $u := \{\beta < \alpha:\nu_n{}^\frown\langle \beta
\rangle \in {\Cal T}_0\}$ has cardinality $\ge \kappa$ but then $\beta
\in u \wedge |\beta \cap u| < \kappa \Rightarrow \nu_n{}^\frown
\langle \beta \rangle \in {\Cal T}_1$ which implies that
$|\text{Suc}_{{\Cal T}_1}(\nu_n)| \ge \kappa$, contradiction to 
the assumption of clause (c).]
\mr
\item "{$(*)_5$}"  $|{\Cal T}_1| \ge \kappa$
\ermn
[Why?  Otherwise by $(*)_4$ the set $A' := \cup\{A_\nu \backslash
\cup\{A_\rho:\rho \in \text{\rm Suc}_{{\Cal T}_0}(\nu)\}:\nu \in {\Cal T}_1\}$ 
is a subset of $\kappa$
of cardinality $< \kappa$ and by clause (d) of $\circledast$ of the
present claim also $A'' = \cup\{f^{-1}_\nu\{0\}:\nu \in {\Cal T}_1\}$
is a subset of $\kappa$ of cardinality $< \kappa$.  So we can
choose $j \in A \backslash (A' \cup A'')$.  Now
we try to choose $\nu_n \in {\Cal T}_1$ 
by induction on $n$ such that $\ell g(\nu_n) = n,\nu_{n+1} \in
\text{\rm Suc}_{{\Cal T}_1}(\nu_n)$ and $j \in A_{\nu_n}$.

So $\nu_0 = <>$ belongs to ${\Cal T}_1$ by $(*)_1 + (*)_3$(a).  Now assume
$\nu_n$ is well defined, then Suc$_{{\Cal T}_0}(\nu_n) = \text{
Suc}_{{\Cal T}_1}(\nu_n)$ by $(*)_4$(2) and our present assumption 
toward contradicting $|{\Cal T}_1|<\kappa$.
\nl

Now $j \notin A',A' \supseteq A_{\nu_n} \backslash
\cup\{A_\rho:\rho \in \text{\rm Suc}_{{\Cal T}_1}(\nu_n)\}$, but $j
\in A_{\nu_n}$ hence clearly $j \in \cup\{A_\rho:\rho \in \text{\rm
Suc}_{{\Cal T}_1}(\nu_n)\}$, so we can choose
$\nu_{n+1}$ as required. So we have carried the definition of 
$\langle \nu_n:n<\omega\rangle$. 

As $j \in A_{\nu_n} \subseteq
B_{\nu_n}$ by $(*)_3(b)$ above, clearly $f_{\nu_n}(j)$ is well
defined (for each $n < \omega$).
As $j \notin A''$ and 
$f^{-1}_{\nu_n}\{0\} \subseteq A''$, so $j \notin
f^{-1}_{\nu_n}\{0\}$, necessarily $f_{\nu_n}(j) \ne 0$ and so 
$f_{\nu_n}(j) > f_{\nu_{n+1}}(j)$ by the choice of 
$f_{\nu_{n+1}}$ in clauses (g)
of $\circledast$.   Hence $\langle f_{\nu_n}(j):n <
\omega\rangle$ is decreasing (sequence of ordinals), 
contradiction.  So $(*)_5$ holds.]
\medskip

Let $n < \omega$ be maximal such that $|{\Cal T}_1 \cap {}^{n \ge}\lambda| <
\kappa$, it exists as $|{\Cal T}_1| \ge \kappa = \text{\rm cf}(\kappa) >
\aleph_0$ and $n=0 \Rightarrow |{\Cal T}_1 \cap {}^{n \ge}\lambda| 
= 1 < \kappa$, and let $\eta \in {\Cal T}_1 \cap {}^n \lambda$ 
be such that Suc$_{{\Cal T}_1}(\eta)$ has $\ge \kappa$ members; it exists
as $\kappa$ is regular.  
We can choose an increasing
sequence $\langle \alpha_i:i < \kappa \rangle$ of ordinals such that
$\alpha_i$ is the $i$-th member of the set $\{\alpha < \lambda:\eta
{}^\frown \langle \alpha \rangle \in {\Cal T}_1\}$
and let $A_i \in [A]^\kappa$ be such that $\langle
\alpha_\ell(A_i,\bar B):\ell \le n \rangle = \eta {}^\frown \langle
\alpha_i\rangle$ and let $\delta = \cup\{\alpha_i:i < \kappa\}$, so
$\delta \in S^\lambda_\kappa$.  Let

$$
\align
A_* = \cup\{A_i:i < \kappa\} \cap B_\eta \backslash \cup\{A_{(\eta
\restriction \ell){}^\frown<\gamma>}:&\ell < \ell g(\eta) \\
  &\text{and } \gamma < \eta(\ell) \text{ and } (\eta \restriction
  \ell) {}^\frown \langle \gamma \rangle \in {\Cal T}_1\}
\endalign
$$
\mn
(note that that number of pairs $(\ell,\gamma)$ as mentioned above is
$<\kappa$).

Clearly $\alpha_\ell(A_*,\bar B) = \eta(\ell)$ for $\ell < \ell
g(\eta)$ hence $\alpha_\ell(A_* \cap A_i,\bar B) 
= \eta(\ell)$ for $i < \kappa,\ell < n$
so clause $(\alpha)$ of (c) of Definition \scite{4a.1}(5)
holds, as well as clause $(\beta)$ because 
$\alpha_n(A_* \cap A_i,\bar B) =
\alpha_i$ for $i < \kappa$ and 
$\langle B_{\eta^\frown\langle \alpha\rangle}:\alpha<\lambda\rangle$
is $\subseteq^*$-increasing.  

\hfill$\square_{\scite{4a.2}}$\margincite{4a.2}
\enddemo
\bigskip

\proclaim{\stag{4b.3} Claim}  If there is a good $(\kappa,{}^{\omega
>}\lambda)$-parameter and $\lambda_1 \in {\frak
 b}^{\text{spc}}_\kappa$ \ub{then} the forcing notion 
$\Bbb P_\kappa$ collapses $\lambda_1$ to $\aleph_0$.
\endproclaim
\bigskip

\demo{Proof}  Let $(\bar B,\bar \nu)$ be a good $(\kappa,{}^{\omega
 >}\lambda)$-parameter. 
\nl
Note
\mr
\item "{$\circledast_1$}"  if $A_1 \subseteq A_2$ are from
$[\kappa]^\kappa$ and $\alpha_\ell(A_2,\bar B)$ is well defined then
$\alpha_\ell(A_1,\bar B)$ is well defined and equal to
$\alpha_\ell(A_2,\bar B)$, recalling Definition \scite{4a.1}(3).
\ermn
Let $\bar h = \langle h_\gamma:\gamma < \lambda_1 \rangle$ exemplify
$\lambda_1 \in {\frak b}^{\text{spc}}_\kappa$, i.e., is as in
Definition \scite{4.1} and \wilog \, $[i < j < \kappa \Rightarrow i <
h_\gamma(i) < h_\gamma(j)]$.  For each $\delta \in
S^\lambda_\kappa$ and $\eta \in {}^{\omega >}\lambda$ let $A_{\eta,\delta,i}
= B_{\eta {}^\frown <\nu_\delta(i+1)>} \backslash \cup 
\{B_{\eta {}^\frown <\nu_\delta(j+1)>}:j < i\}$ for $i < \kappa$ so
$\langle A_{\eta,\delta,i}:i < \kappa \rangle$ are pairwise disjoint
subsets of $\kappa$ (each of cardinality $\kappa$).  For 
$n<\omega$ and $A
\in [\kappa]^\kappa$ we try to define an ordinal $\beta_n(A,\bar
B,\bar \nu,\bar h)$ as follows:
\mr
\item "{$\circledast_2$}"  $\beta_n(A,\bar B,\bar \nu,\bar h) = 
\gamma$ \ub{iff} for some $\eta \in {}^n \lambda$ and 
$\delta \in S^\lambda_\kappa$
we have $\langle \alpha_\ell(A,\bar B):\ell \le
n \rangle = \eta {}^\frown \langle \delta \rangle$ so in particular is
well defined and $A \subseteq^* \cup\{A_{\eta,\delta,i} \cap
h_\gamma(i):i < \kappa\}$ but for every $\beta < \gamma$ we have $A \cap
\cup\{A_{\eta,\delta,i} \cap h_\beta(i):i < \kappa\} \in [\kappa]^{<
\kappa}$.
\ermn
Next we define a $\Bbb P_\kappa$-name ${\underset\tilde {}\to \beta_n}
= {\underset\tilde {}\to \beta_n}(\bar B,\bar \nu,\bar h)$ by:
\mr
\item "{$\circledast_3$}"  for $\bold G \subseteq \Bbb P_\kappa$
generic over $\bold V:{\underset\tilde {}\to \beta_n}[\bold G] =
\gamma$ \ub{iff} for some $A \in \bold G$ we have $\beta_n(A,\bar B,\bar \nu,
\bar h) = \gamma$ or there is no such $A$ and $\gamma=0$.
\ermn
Now
\mr
\item "{$\circledast_4$}"  if $A \in [\kappa]^\kappa$ and (${\Cal
T}^0_A,{\Cal T}^1_A$) $n,\eta,\delta$ are chosen 
as in the proof of \scite{4a.2}(2), then $u := 
\{\beta < \lambda_1:A \nVdash_{\Bbb P_\kappa} 
``{\underset\tilde {}\to \beta_n}(\bar B,\bar \nu,\bar h) \ne 
\beta"\}$ is a $\kappa$-closed unbounded subset of $\lambda_1$.
\ermn
[Why?  We know that $w := \{i < \kappa:A \cap A_{\eta,\delta,i} \in
[\kappa]^\kappa\}$ has cardinality $\kappa$.  Why is $u$ ``unbounded"?  For
any $\gamma_1 < \lambda_1$, we define a function $h \in {}^\kappa
\kappa$ as follows,  $h(i)$ is the minimal $i_1 < \kappa$ such that for some
$i_0,i < i_0 < i_1$ the set $A \cap A_{\eta,\delta,i_0} \cap i_1
\backslash h_{\gamma_1}(i_0)$ is not empty, clearly $h$ is well
defined because $|w| = \kappa$.  So for some $\gamma_2 \in
(\gamma_1,\lambda_1)$ the set $v := \{i < \kappa:h(i) <
h_{\gamma_2}(i)\}$ has cardinality $\kappa$.  Let $C$ be the club $\{\delta <
\kappa:\delta$ is a limit ordinal and $i < \delta \Rightarrow h(i) <
\delta \wedge h_{\gamma_2}(i) <  \delta\}$ and let $\langle
\alpha_\varepsilon:\varepsilon < \kappa\rangle$ list $C \cup \{0\}$
increasing order $\kappa$ and let $A' = \cup\{A \cap A_{\eta,\delta,i} \cap
[\alpha_\varepsilon,\alpha_{\varepsilon +1}):i < \kappa,\varepsilon <
\kappa$ and $\alpha_\varepsilon \le i < \alpha_{\varepsilon +1}\}$,
now $A' \in [\kappa]^\kappa$ (really $i < j < \kappa \Rightarrow i <
h_{\gamma_1}(i) < h_{\gamma_2}(j))$.  So $\Bbb P_\kappa
\models ``A \le A'$" and $A' \Vdash ``{\underset\tilde {}\to
\beta_n}(\bar B,\bar \nu,\bar h) \in (\gamma_1,\gamma_2]"$, recalling
that the $h_\gamma$'s are $<_{J^{\text{bd}}_\kappa}$ 
increasing.  Why ``the set $u$ is $\kappa$-closed"
(that is the limit of any increasing sequence of length $\kappa$ of
members belong to it)?  Easy, too.] 

Let $\langle S_\varepsilon:\varepsilon < \lambda_1\rangle$ be pairwise
disjoint stationary subsets of $S^{\lambda_1}_\kappa$ 
and define $g^*:\lambda_1
\rightarrow \lambda_1$ by $g^*(\gamma) =\varepsilon$ if $\gamma \in
S_\varepsilon \vee (\gamma \in \lambda_1 \backslash 
\dbcu_{\zeta < \lambda_1} S_\zeta \wedge \varepsilon =0)$.  So
\mr
\item "{$\circledast_5$}"  for every $p \in \Bbb P_\kappa$ for some $n$,
for every  $\varepsilon < \lambda_1,\ p \nVdash 
``g^*({\underset\tilde {}\to \beta_n}) \ne \varepsilon"$
\ermn
so we are done. \hfill$\square_{\scite{4b.3}}$\margincite{4b.3} 
\enddemo
\bn
Now we arrive to the main point.
\proclaim{\stag{4b.4} Main Claim}  1) If $\Bbb P_\kappa$ does not satisfy the
$\chi$-c.c. \ub{then} forcing with $\Bbb P_\kappa$ 
collapses $\chi$ to $\aleph_0$.
\nl
2) There is $\langle \bar A_\alpha:\alpha < {\frak b}_\kappa\rangle$
such that $\bar A_\alpha = \langle A_{\alpha,i}:i < \kappa\rangle$ is
a sequence of pairwise disjoint subsets of $\kappa$ each of
cardinality $\kappa$ (\wilog \, each is a partition of $\kappa$) such
that for every $B \in [\kappa]^\kappa$ for some $\alpha < {\frak
b}_\kappa$ we have $i < \kappa \Rightarrow \kappa = |A_{\alpha,i} \cap
B|$; i.e., for every $i < \kappa$ not just for $\kappa$ many $i <
\kappa$.
\endproclaim
\bigskip

\demo{Remark} 
1) In part (2) we can replace ${\frak b}_\kappa$ by any 
$\lambda\in {\frak b}^{\text{spc}}_\kappa$, but this does not 
add information. The proof gives a little more for ``many'' 
$\alpha<\lambda$.
\nl
2) In case 1 we could have assumed
${\frak b}_\kappa>\kappa^+$, 
this suffice
\nl
3) We could have seperated the different roles of $\lambda$
in the proof of case 1. Say
{\roster
\itemitem{ $(a)$ } $(\bar B, \bar \nu^1)$ will be a good 
$(\kappa,{}^{\omega>}(\lambda_1))$-parameter,
\sn
\itemitem{ $(b)$ } $\langle h_\alpha:\alpha<\lambda_2\rangle$ 
exemplify $\lambda_2\in {\frak b}^{\text{spc}}_\kappa$ 
and $\langle \nu^*_\delta:\delta\in S^{\lambda_2}_\kappa\rangle$ 
is an $S^{\lambda_2}_\kappa$-ladder system 
(so $\delta^*\in S^{\lambda_2}_\kappa$ in the proof)
\endroster}
4) Actually, we can revise case 2 to cover Case 1, too: for
$\delta_* \in S^\lambda_{\kappa^+}$ choose $C'_{\delta_*}$ a club of
$\delta_*$ of order type $\kappa^+$.  Now for each $\delta$ we can
repeat the construction of names from the proof of Case 2, for each $p
\in \Bbb P_\kappa$ for some $\delta_*$ 
we succeed to show $\circledast$ below.
\enddemo

\demo{Proof}  The proof is divided to two cases.
\mn
\ub{Case 1}:  $\lambda \in {\frak b}^{\text{spc}}_\kappa,\lambda >
\kappa^+$, e.g. $\lambda={\frak b}_\kappa$.

So $\lambda$ is regular $> \kappa^+$ and a good 
$(\kappa,{}^{\omega >}\lambda)$ sequence $\bar B$ exists (by \scite{4a.2}).

Let $\bar \nu = \langle \nu_\delta:\delta \in S^\lambda_\kappa
\rangle$ be such that $\nu_\delta \in {}^\kappa\delta$ is increasing
continuous with limit $\delta$ and $\bar \nu$ guesses clubs (i.e. for
every club $C$ of $\lambda$, for stationarily many $\delta \in
S^\lambda_\kappa$ we have Rang$(\nu_\delta) \subseteq C)$; exists by
\cite[III,\S2]{Sh:g} because $\lambda = \text{\rm cf}(\lambda) > \kappa^+$.
As $\bar B$ is a good $(\kappa,{}^{\omega >}\lambda)$-sequence,
$(\bar B,\bar \nu)$ is a good $(\kappa,{}^{\omega
>}\lambda)$-parameter by \scite{2.3.1} (or use \scite{4a.2}).  

Let $\langle h_\alpha:\alpha < \lambda\rangle$ exemplify $\lambda
\in {\frak b}^{\text{spc}}_\kappa$ without loss of generality 
$i<j<\kappa\Rightarrow i<h(i)<h(j)$.

For $\eta \in {}^{\omega >} \lambda,\delta \in S^\lambda_\kappa$ and
$i < \kappa$, recall that $A_{\eta,\delta,i} = B_{\eta {}^\frown
<\nu_\delta(i+1)>} \backslash \cup\{B_{\eta {}^\frown
<\nu_\delta(j+1)>}:j<i\}$ and let $\beta_n(A,\bar B,\bar \nu,\bar h),
{\underset\tilde {}\to \beta_n} = {\underset\tilde {}\to \beta_n}(\bar
B,\bar \nu,\bar h)$ be defined as in
the proof of \scite{4b.3}.  
For $\eta \in {}^{\omega >}\lambda,\delta
\in S^\lambda_\kappa$
and $\gamma < \lambda$ let $B^*_{\eta,\delta,\gamma} := \cup\{A_{\eta,\delta,i}
\cap h_\gamma(i):i < \kappa\}$. 
So clearly (for each $\eta \in {}^{\omega >}\lambda,\delta \in
S^\lambda_\kappa$) the sequence $\langle B^*_{\eta,\delta,\gamma}:\gamma <
\lambda\rangle$ is $\subseteq^*$-increasing.
For $\delta^* \in S^\lambda_\kappa$ and $i < \kappa$
let $A^*_{\eta,\delta,\delta^*,i} :=
B^*_{\eta,\delta,\nu_{\delta^*}(i+1)} \backslash \cup
\{B^*_{\eta,\delta,\nu_{\delta^*}(j+1)}:j<i\}$.  So $\langle
A^*_{\eta,\delta,\delta^*,i}:i < \kappa \rangle$ are pairwise disjoint
subsets of $\kappa$.  Note that (by the proof of \scite{4b.3} but not used)
 for each pair $(\eta,\delta)$ as above
for some club $E_{\eta,\delta}$ of $\lambda$, for every $\delta^* \in
S^\lambda_\kappa \cap E_{\eta,\delta}$ and $i <
\kappa,A^*_{\eta,\delta,\delta^*,i}$ has cardinality $\kappa$.  We
shall show during the proof of (1) that $\{\langle
A^*_{\eta,\delta,\delta^*,i}:i < \kappa\rangle:\eta \in {}^{\omega
>}\lambda,\delta \in S^\lambda_\kappa,\delta^* \in S^\lambda_\kappa\}$
is as required in part (2), so this will prove part (2) 
when ${\frak b}_\kappa > \kappa^+$.

Let $\langle X^*_\xi:\xi < \chi\rangle$ be an antichain of
$\Bbb P_\kappa$, it exists by the assumption.  We now for
$\eta,\delta,\delta^*$ as above define $\Bbb
P_\kappa$-names ${\underset\tilde {}\to \gamma_{\eta,\delta,\delta^*}}$: 
for $\bold G \subseteq \Bbb P_\kappa$ generic over
$\bold V$ we let: 
\mr
\item "{$\circledast_0$}"  ${\underset\tilde {}\to
\gamma_{\eta,\delta,\delta^*}}[\bold G] =
\xi$ \ub{iff} for some $A \in \bold G,n < \omega$ and $\eta \in {}^n
\lambda$ and $\delta,\delta^* \in S^\lambda_\kappa$ we have:
{\roster
\itemitem{ $(a)$ }  $\langle \alpha_\ell(A,\bar B):\ell <n \rangle =
\eta$ so in particular is well defined
\sn
\itemitem{ $(b)$ }  $\alpha_n(A,\bar B) = \delta \in S^\lambda_\kappa$
\sn
\itemitem{ $(c)$ }  $\beta_n(A,\bar B,\bar \nu,\bar h) 
= \delta^* \in S^\lambda_\kappa$
\sn
\itemitem{ $(d)$ }  $A \cap A^*_{\eta,\delta,\delta^*,i}$ has at most 
one member for each $i < \kappa$
\sn
\itemitem{ $(e)$ }  $A \subseteq \cup\{A^*_{\eta,\delta,\delta^*,i}:i \in
X^*_\xi\}$
\endroster}
\ermn
Note that demands (a),(b),(c) are natural but actually 
not being used; with them
we could have defined the $\Bbb P_\kappa$-names ${\underset\tilde {}\to
\gamma_n}$ which is ${\underset\tilde {}\to
\gamma_{\eta,\delta,\delta^*}}$ when defined.  Now clearly
\mr
\item "{$\circledast_1$}"  ${\underset\tilde {}\to
\gamma_{\eta,\delta,\delta^*}}$ is a $\Bbb P_\kappa$-name of an ordinal
$< \chi$ (may have no value) 
\sn
\item "{$\circledast_2$}"  for every $p \in \Bbb P_\kappa$ for some
$\eta \in {}^{\omega >} \lambda$ and $\delta,\delta^* \in
S^\lambda_\kappa$, for every 
$\varepsilon < \chi$ there is $q$ such that $p \le q \in \Bbb
P_\kappa$ and $q \Vdash_{\Bbb P_\kappa} ``{\underset\tilde {}\to
\gamma_{\eta,\delta,\delta^*}} = \varepsilon"$.
\ermn
[Why?  We start as in the proof of \scite{4b.3}.  First there are $n <
\omega,\eta \in {}^n \lambda$ and 
$\delta \in S^\lambda_\kappa$ such that $p \cap
A_{\eta,\delta,i} \in [\kappa]^\kappa$ for $\kappa$ many ordinals $i < \kappa$.
\nl
\comment
Second, the set $W_p \subseteq \lambda$ is unbounded in $\lambda$ (by
the proof of \scite{4b.3})
where
$$
\align W_p = \{\beta < \lambda:&\text{for some } \gamma \in
(\beta,\lambda) \text{ we have} \\
  &p \cap B^*_{\eta,\delta,\gamma} \backslash
  B^*_{\eta,\delta,\beta} \text{ is from } [\kappa]^\kappa\}.
\endalign
$$
\endcomment
Second, there is a club $C_p$ of $\lambda$ such that:
\nl
if $\beta < \gamma < \lambda$ are from $C_p$ then $p \cap
B^*_{\eta,\delta,\gamma} \backslash B^*_{\eta,\delta,\beta} \in
[\kappa]^\kappa$. 
Indeed, $C_p=\{\gamma< \lambda$: for every $\beta<\gamma$ 
the set $p\cap
B^*_{\eta,\delta,\gamma} \setminus
  B^*_{\eta,\delta,\beta} \text{ is from } [\kappa]^\kappa\}$ is as required.

Now by the choice of $\bar \nu$, i.e., club guessing, there is $\delta^* \in
\text{\rm acc}(C_p) \cap S^\lambda_\kappa$ such that $(\forall i <
\kappa)(\nu_{\delta^*}(i) \in C_p)$.  So (as we have used
$\nu_{\delta^*}(i+1),\nu_{\delta^*}(j+1)$ in the definition of
$A^*_{\eta,\delta,\delta^*,i}$) 

$$
i < \kappa \Rightarrow p \cap A^*_{\eta,\delta,\delta^*,i} \in
[\kappa]^\kappa.
$$
This fulfills the promise needed for proving part (2) in the present
case 1.
Choose $\zeta_i \in p \cap A^*_{\eta,\delta,\delta^*,i}$ for $i <
\kappa$.  Now for every $\xi < \chi$ let $q_\xi =
\{\zeta_i:i \in X^*_\xi\}$.  
Recall that $\langle X^*_\zeta:\zeta < \chi \rangle$ is an antichain
in $\Bbb P_\kappa$.  Clearly for $\xi < \chi$ we
have $\Bbb P_\kappa \models ``p \le q_\xi"$ and $q_\xi
\Vdash ``{\underset\tilde {}\to \gamma_{\eta,\delta,\delta^*}} =
\xi"$; so we have finished proving $\circledast_2$.] 

This is enough for proving
\mr
\item "{$\circledast_3$}"   forcing with $\Bbb P_\kappa$ collapse
$\chi$ to $\aleph_0$.
\nl
[Why?  By $\circledast_1 + \circledast_2$ we know that $\Vdash_{\Bbb
P_\kappa} ``\chi = \{{\underset\tilde {}\to
\gamma_{\eta,\delta,\delta^*}}:
\eta \in {}^{\omega >}\lambda,\delta \in
S^\lambda_\kappa$ and $\delta^* \in S^\lambda_\kappa\}"$, so it is forced
that $|\chi| \le |\lambda|$.  As we already have by \scite{4b.3} that
$\Vdash_{\Bbb P_\kappa} ``|\lambda| = \aleph_0"$, we are done.]
\endroster
\enddemo
\bn
\ub{Case 2}:  ${\frak b}_\kappa = \kappa^+$.

Let 
$\lambda = \kappa^+$ and $\bar B$ be a good $(\kappa,{}^{\omega
>}\lambda)$-sequence.  Let $\langle S_\varepsilon:\varepsilon < \kappa
\rangle$ be a partition of $S^{\kappa^+}_\kappa$ to (pairwise disjoint)
stationary sets.  For $\alpha < \kappa^+$ let $\langle u^\alpha_i:i <
\kappa \rangle$ be an increasing continuous sequence of subsets of
$\alpha$ each of cardinality $< \kappa$ with union $\alpha$ and \wilog \,
$\alpha < \beta \Rightarrow (\forall^* i < \kappa)(u^\alpha_i =
u^\beta_i \cap \alpha)$.  Let $\bar h = \langle h_\beta:\beta <
\kappa^+\rangle$ exemplifying $\kappa^+ \in {\frak
b}^{\text{spc}}_\kappa$ be such that each $h_\beta$ is strictly increasing,
$(\forall i)h_\beta(i) > i$ and let
$C_\beta = \{\delta < \kappa:\delta$ is a limit ordinal and
for every $i < \delta$ we have $h_\beta(i) <
\delta\}$ and let $(\bar B,\bar \nu)$ be a good $(\kappa,{}^{\omega
>}\lambda)$-parameter; exists by \scite{4a.2}(2).  
Now for $\eta \in {}^{\omega >}\lambda$ and $\delta \in 
S^\lambda_\kappa$ we define
$A_{\eta,\delta,i}(i < \kappa),B^*_{\eta,\delta,\gamma}(\gamma < \lambda)$
as in Case 1.  Now for
$\eta \in {}^{\omega >}\lambda,\delta \in S^\lambda_\kappa,\alpha <
\kappa^+$ and $\beta < \kappa^+$ we define the sequence $\langle
Y_{\eta,\delta,\alpha,\beta,\gamma}:\gamma < \alpha \rangle$ by

$$
Y_{\eta,\delta,\alpha,\beta,\gamma} := \cup\{B^*_{\eta,\delta,\gamma} \cap
[i,\text{Min}(C_\beta \backslash (i+1)) \backslash 
\cup\{B^*_{\eta,\delta,\gamma_1}:
\gamma_1\in \gamma \cap u^\alpha_i\}:i \in C_\beta \text{ satisfy } \gamma
\in u^\alpha_i\}.
$$
\mn
So $\langle Y_{\eta,\delta,\alpha,\beta,\gamma}:\gamma < \alpha
\rangle$ is a sequence of pairwise disjoint subsets of $\kappa$ and 
for $\varepsilon < \kappa$ let

$$
Z_{\eta,\delta,\alpha,\beta,\varepsilon} :=
\cup\{Y_{\eta,\delta,\alpha,\beta,\gamma}:
\gamma \in S_\varepsilon \cap \alpha\}.
$$
\mn
Clearly
\mr
\item "{$\boxdot_1$}"  $\bar Z_{\eta,\delta,\alpha,\beta}=\langle
Z_{\eta,\delta,\alpha,\beta,\varepsilon}:\varepsilon < \kappa \rangle$
is a sequence of pairwise disjoint subsets of $\kappa$.
\ermn

We shall show during the proof of (1) that 
$$
\left< \langle
Z_{\eta,\delta,\alpha,\beta,\varepsilon}:\varepsilon < \kappa
\rangle:\eta \in {}^{\omega >}\lambda,\,\delta \in
S^\lambda_\kappa,\,\alpha < \lambda,\,\beta < \lambda \right>
$$ 
exemplify part (2);
you may wonder: possibly for some quadruple $(\eta,\delta,\beta,\zeta)$
we do not have $(\forall \epsilon<\kappa) 
[|Z_{\eta,\delta,\alpha,\beta,\epsilon}|=\kappa]$, so?
However the quadruple $(\eta,\delta,\alpha,\beta)$ for which 
this fails, cannot satisfy the desired property in part (2), 
so we can just omit them. 

Let $\langle X^*_\xi:\xi < \chi\rangle$ be a family of sets from
$[\kappa]^\kappa$ such that the intersection of any two have
cardinality $< \kappa$, it exists as $\Bbb P_\kappa$ fail the $\chi$-c.c..  
For each $\eta \in {}^{\omega >}\lambda,\delta \in S^\lambda_\kappa,
\alpha < \kappa^+$ and $\beta < \kappa^+$ we define a $\Bbb
P_\kappa$-name ${\underset\tilde {}\to \tau_{\eta,\delta,\alpha,\beta}}$ as
follows: 
\mr
\item "{$\boxdot_2$}"  for $\bold G \subseteq \Bbb P_\kappa$ 
generic over $\bold V,
{\underset\tilde {}\to \tau_{\eta,\delta,\alpha,\beta}}[\bold G]=\xi$
iff
{\roster
\itemitem{ $(\alpha)$ }  for some $A \in \bold G$ we have
\sn
\itemitem{ ${{}}$ }  $(a) \quad \varepsilon < \kappa \Rightarrow A \cap
Z_{\eta,\delta,\alpha,\beta,\varepsilon}$ has at most one member
\sn
\itemitem{ ${{}}$ }  $(b) \quad A \subseteq
\cup\{Z_{\eta,\delta,\alpha,\beta,\varepsilon}:\varepsilon \in 
X^*_\xi\}$
\sn
\itemitem{ $(\beta)$ }  if for no $A \in \bold G$ does (a)+(b) hold 
and $\xi=0$.
\endroster}
\ermn
Clearly
\mr
\item "{$\boxdot_3$}"  ${\underset\tilde {}\to
\tau_{\eta,\gamma,\alpha,\beta}}$ is a well defined ($\Bbb P_\kappa$-name) 
(by $\boxdot_2$).
\ermn
Now
\mr
\item "{$\boxdot_4$}"   for every $p \in \Bbb P_\kappa$, for some
$\eta \in {}^{\omega >}\lambda,\delta \in S^\lambda_\kappa,
\alpha < \kappa^+,\beta < \kappa^+$ we
have: for every $\xi < \chi$ for some $q \in \Bbb P_\kappa$ above $p$ we
have $q \Vdash ``{\underset\tilde {}\to \tau_{\eta,\delta,\alpha,\beta}} =
\xi"$ and $\epsilon<\kappa
\Rightarrow |Z_{\eta,\delta,\alpha,\beta}\cap p|=\kappa$.
\ermn
As in Case 1, this is enough for 
proving that $\Bbb P_\kappa$ collapse $\chi$ to
$\lambda = \kappa^+$.  But by \scite{4b.3} we already know 
that forcing with $\Bbb P_\kappa$ collapses 
$\kappa^+$ to $\aleph_0$ and so we are done.

Note: we can eliminate $\eta$ from the ${\underset\tilde {}\to
\tau_{\eta,\delta,\alpha,\beta}}$, but not worth it.  So we are left with
proving $\boxdot_4$.

Why does $\boxdot_4$ hold?  First, as in the earlier cases, find 
$\eta\in {}^{\omega >}\lambda$ and $\delta \in
S^\lambda_\kappa$ such that $p \cap A_{\eta,\delta,i} \in [\kappa]^\kappa$
for $\kappa$ ordinals $i < \kappa$.  
\comment
Second, as in the previous proof
$W_p \subseteq \lambda =
\text{\rm sup}(W_p)$ where $W_p := \{\beta < \lambda$: for some $\gamma
\in (\beta,\lambda)$ we have $p \cap B^*_{\eta,\delta,\gamma} \backslash
B^*_{\eta,\delta,\beta} \in [\kappa]^\kappa\}$.
\endcomment  
Second, for some club
$C_p$ of $\lambda$ we have $\beta < \gamma \wedge \gamma \in C_p
\Rightarrow p \cap B^*_{\eta,\delta,\gamma} \backslash
B^*_{\eta,\delta,\beta} \in [\kappa]^\kappa$.
As $S_\varepsilon$ (for $\varepsilon
<\kappa$) is a stationary subset of $\lambda$ and $C_p$ a club of
$\lambda$ for each $\varepsilon < \kappa$ we can choose 
$\gamma^*_\varepsilon \in S_\varepsilon \cap C_p$.  Hence
there is $\alpha^* < \kappa^+$ large enough such
that $\varepsilon < \kappa \Rightarrow \gamma^*_\varepsilon < \alpha^*
\in C_p$.  Now define a function $h:\kappa \rightarrow
\kappa$ by induction on $i$, as follows:

$$
\align
h(i) = \text{\rm Min}\{j:&j \in (i,\kappa) \text{ and } i_1 < i
\Rightarrow h(i_1) < j \text{ and } \\
  &\text{ if the pair } 
(\gamma,\epsilon) \text{ is such that } 
\gamma \in u^{\alpha^*}_i \cap S_\varepsilon \text{ then} \\
  &p \cap (i,j) \cap B^*_{\eta,\delta,\gamma} \backslash \cup
  \{B^*_{\eta,\delta,\gamma_1}:\gamma_1\in \gamma \cap u^{\alpha^*}_i\} 
\text{ is not empty}\}.
\endalign
$$
\mn
it is well defined as for a given $i < \kappa$ the number of pairs
$(\gamma,\varepsilon)$ such that $\gamma \in u^{\alpha^*}_i \cap
S_\varepsilon$ is 
$\leq |u^{\alpha^*}_i|< \kappa$ and is increasing; next we define

$$
C = \{j < \kappa:j \text{ is a limit ordinal such that } i <j
\Rightarrow h(i) < j\}.
$$
\mn
Clearly $C$ is a club of $\kappa$ and let $h' \in {}^\kappa \kappa$ be
defined by $h'(i) = h(\text{Min}(C \backslash (i+1))$.
By the choice of $\langle h_\beta:\beta < \lambda\rangle$ there is
$\beta < \lambda$ such that for $\kappa$ many ordinals $i <
\kappa,h'(i) < h_\beta(i)$. Recall that 
$C_\beta = \{\delta < \kappa:\delta$ is a limit ordinal and
for every $i < \delta$ we have $h_\beta(i) <
\delta\}$.

\medskip
So $W_1 = \{i < \kappa:h'(i) < h_\beta(i)\}$, by
the choice of $\beta$ clearly $W_1 \in [\kappa]^\kappa$.  Let
$\langle i^0_j:j < \kappa \rangle$ be an enumeration of 
the club $C\cap \text{acc}(C_\beta)$ of $\kappa$ 
in an increasing order, so clearly ${\Cal U} := \{j < \kappa:W_1 \cap
[i^0_j,i^0_{j+1}) \ne \emptyset\}$ is unbounded in $\kappa$.  For each
$j \in {\Cal U}$ let $i^2_j$ 
be the first member of $W_1 \cap [i^0_j,i^0_{j+1})$, 
then let $i^1_j=\text{sup}(C_\beta\cap (i^2_j+1))$, 
it is well defined as $i^0_j\in C \cap \text{acc}(C_\beta)$, and 
so $i^0_j\leq i^1_j$ and let 
$i^3_j=\text{ min } (C_\beta\setminus (i^2_j+1))$ so 
$i^2_j<i^3_j$ and
\mr
\item "{$(*)$}"  $i^0_j \le i^1_j \leq i^2_j < h(i^2_j) < h'(i^2_j) <
h_\beta(i^2_j) < i^3_{j}$.
\ermn

[why? as said above by the choice of $i^1_j$ by the choice of $h$, 
by the choice of the pair $(C,h')$, by $i^2_j\in W_1$, 
by the choice of $i^3_j$ and $C_\beta$ respectively.]
\mr 
\item ''{$(*)'$}'' $i^1_j<i^3_j$ are successive members of 
$C_\beta$ 

[Why? both are members of $C_\beta$ by their choices hence 
it is enough to prove that $C_\beta \cap (i^1_j,i^3_j)=\emptyset$.
But $C_\beta \cap (i^1_j,i^2_j]=\emptyset$ by the choice of 
$i^1_j$ and $\beta\cap (i^2_j,i^3_j)=\emptyset$ by the choice 
of $i^3_j$]
\ermn

Now for each
$\varepsilon < \kappa$ we know that $\gamma^*_\varepsilon \in \alpha^*
\cap S_\varepsilon \cap C_p \subseteq \alpha^* = \cup\{u^{\alpha^*}_i:i <
\kappa\}$ and $\langle u^{\alpha^*}_i:i < \kappa \rangle$ is
$\subseteq$-increasing hence for some $j(\varepsilon) < \kappa$ if $j
\in {\Cal U}\setminus j(\varepsilon)$ then
$\gamma^*_\varepsilon \in u^{\alpha^*}_{i^0_j}$ hence by the choice of 
$h(i^1_j)$ and $(*)$ we have $p \cap (i^1_j,i^3_{j}) 
\cap B^*_{\eta,\delta,\gamma^*_\varepsilon}
\backslash \cup\{B^*_{\eta,\delta,\gamma_1}:\gamma_1 \in
\gamma^*_\varepsilon \cap u^{\alpha^*}_{i^0_j}\}$ is not empty;
but $i^1_j<i^3_j$ are successive members of $C_\beta$ by $(*)'$ so 
the definition of
$Y_{\eta,\delta,\alpha^*,\beta,\gamma^*_\varepsilon}$ implies that $p
\cap Y_{\eta,\delta,\alpha^*,\beta,\gamma^*_\varepsilon} \cap
(i^1_j,i^3_{j}) \ne \emptyset$.

As this holds for every large enough 
$j\in {\Cal U}$ i.e., for every $j \in {\Cal U} \backslash
j(\varepsilon)$ and ${\Cal U} \in [\kappa]^\kappa$ it follows that $p
\cap Y_{\eta,\delta,\alpha^*,\beta,\gamma^*_\varepsilon} \in
[\kappa]^\kappa$.  By the definition of
$Z_{\eta,\delta,\alpha^*,\beta,\varepsilon}$ it follows that $p \cap
Z_{\eta,\delta,\alpha^*,\beta,\varepsilon} \in [\kappa]^\kappa$.

We have proved this for every $\epsilon<\kappa$.
Choose $\zeta_\varepsilon \in p \cap
Z_{\eta,\delta,\alpha^*,\beta,\varepsilon}$
for every $\epsilon<\kappa$.  Now for each $\xi < \chi$
let

$$
q_\xi = \{\zeta_\varepsilon:\varepsilon \in X^*_\xi\}.
$$
\mn
So clearly:

$$
\xi < \chi \Rightarrow \Bbb P_\kappa \models ``p \le q_\xi" \text{ and }
q_\xi \Vdash_{\Bbb P_\kappa} 
``{\underset\tilde {}\to \tau_{\eta,\delta,\alpha^*,\beta}} = \xi".  
$$
\sn
${{}}$  \hfill$\square_{\scite{4b.4}}$\margincite{4b.4}
\bigskip

\demo{\stag{4b.5} Conclusion}  If $\kappa$ is regular uncountable and
$\Bbb P_\kappa$ fail the $2^\kappa$-c.c. \ub{then} 
comp$(\Bbb P_\kappa)$ is isomorphic to the completion 
of Levy$(\aleph_0,2^\kappa)$.
\enddemo
\newpage

% % PRIVATE excursion starts here
% NO IT DOESNT JK
%     \shlhetal % number 1/1
% % back from PRIVATE part number 1/1

\nocite{ignore-this-bibtex-warning} 
%% you may want to move the following lines up a bit
\newpage
    
REFERENCES.  
\bibliographystyle{lit-plain}
\bibliography{lista,listb,listx,listf,liste}

\def\germ{\frak} \def\scr{\cal} \ifx\documentclass\undefinedcs
  \def\bf{\fam\bffam\tenbf}\def\rm{\fam0\tenrm}\fi % f**k-amstex!
  \def\defaultdefine#1#2{\expandafter\ifx\csname#1\endcsname\relax
  \expandafter\def\csname#1\endcsname{#2}\fi} \defaultdefine{Bbb}{\bf}
  \defaultdefine{frak}{\bf} \defaultdefine{=}{\B} % doublef**k-amstex!!
  \defaultdefine{mathfrak}{\frak} \defaultdefine{mathbb}{\bf}
  \defaultdefine{mathcal}{\cal}
  \defaultdefine{beth}{BETH}\defaultdefine{cal}{\bf} \def\bbfI{{\Bbb I}}
  \def\mbox{\hbox} \def\text{\hbox} \def\om{\omega} \def\Cal#1{{\bf #1}}
  \def\pcf{pcf} \defaultdefine{cf}{cf} \defaultdefine{reals}{{\Bbb R}}
  \defaultdefine{real}{{\Bbb R}} \def\restriction{{|}} \def\club{CLUB}
  \def\w{\omega} \def\exist{\exists} \def\se{{\germ se}} \def\bb{{\bf b}}
  \def\equivalence{\equiv} \let\lt< \let\gt>
  \def\implies{\Rightarrow}\def\mathfrak{\bf}\def\germ{\frak} \def\scr{\cal}
  \ifx\documentclass\undefinedcs
  \def\bf{\fam\bffam\tenbf}\def\rm{\fam0\tenrm}\fi % f**k-amstex!
  \def\defaultdefine#1#2{\expandafter\ifx\csname#1\endcsname\relax
  \expandafter\def\csname#1\endcsname{#2}\fi} \defaultdefine{Bbb}{\bf}
  \defaultdefine{frak}{\bf} \defaultdefine{=}{\B} % doublef**k-amstex!!
  \defaultdefine{mathfrak}{\frak} \defaultdefine{mathbb}{\bf}
  \defaultdefine{mathcal}{\cal}
  \defaultdefine{beth}{BETH}\defaultdefine{cal}{\bf} \def\bbfI{{\Bbb I}}
  \def\mbox{\hbox} \def\text{\hbox} \def\om{\omega} \def\Cal#1{{\bf #1}}
  \def\pcf{pcf} \defaultdefine{cf}{cf} \defaultdefine{reals}{{\Bbb R}}
  \defaultdefine{real}{{\Bbb R}} \def\restriction{{|}} \def\club{CLUB}
  \def\w{\omega} \def\exist{\exists} \def\se{{\germ se}} \def\bb{{\bf b}}
  \def\equivalence{\equiv} \let\lt< \let\gt>
\begin{thebibliography}{KjSh 720}
\makeatletter \renewcommand{\@biblabel}[1]{[#1]} \makeatother
\def\eprintfootnotetext{References of the form {\tt math.XX/$\cdots$}
 refer to {\tt arXiv.org} }
\ifx\documentstyle\undefinedcontrolsequence
   \def\anyfootnote{\footnote{*}}
   \else\def\anyfootnote{\footnote}\fi
\def\eprintfn{\ifEprint\anyfootnote{\eprintfootnotetext}\fi\Eprintfalse }
\newif\ifEprint  \Eprinttrue

\bibitem[BaFr87]{BaFr87}Bohuslav Balcar and Franti{\v s}ek Fran{\v e}k.
\newblock {Completion of factor algebras of ideals}.
\newblock {\em Proceedings of the American Mathematical Society}, {\bf
  100}:205--212, 1987.

\bibitem[BPS]{BPS}Bohuslav Balcar, Jan Pelant, and Petr Simon.
\newblock {The space of ultrafilters on $N$ covered by nowhere dense sets}.
\newblock {\em Fundamenta Mathematicae}, {\bf CX}:11--24, 1980.

\bibitem[BaSi88]{BaSi88}Bohuslav Balcar and Petr Simon.
\newblock {On collections of almost disjoint families}.
\newblock {\em Commentationes Mathematicae Universitatis Carolinae}, {\bf
  29}:631--646, 1988.

\bibitem[BaSi89]{BaSi89}Bohuslav Balcar and Petr Simon.
\newblock {Disjoint refinement}.
\newblock In {\em Handbook of Boolean Algebras}, volume~2, pages 333--388.
  North-Holland, 1989.
\newblock Monk D., Bonnet R. eds.

\bibitem[BaSi95]{BaSi95}Bohuslav Balcar and Petr Simon.
\newblock {Baire number of the spaces of uniform ultrafilters}.
\newblock {\em Israel Journal of Mathematics}, {\bf 92}:263--272, 1995.

\bibitem[Ba]{Ba}James~E. Baumgartner.
\newblock {Almost disjoint sets, the dense set problem and partition calculus}.
\newblock {\em Annals of Math Logic}, {\bf 9}:401--439, 1976.

\bibitem[KjSh 720]{KjSh:720}Menachem Kojman and Saharon Shelah.
\newblock {Fallen Cardinals}.
\newblock {\em Annals of Pure and Applied Logic}, {\bf 109}:117--129, 2001.
\newblock math.LO/0009079.

\bibitem[Sh:g]{Sh:g}Saharon Shelah.
\newblock {\em {Cardinal Arithmetic}}, volume~29 of {\em {Oxford Logic
  Guides}}.
\newblock {Oxford University Press}, 1994.

\bibitem[Sh 506]{Sh:506}Saharon Shelah.
\newblock {The pcf-theorem revisited}.
\newblock In {\em {The Mathematics of Paul Erd\H{o}s, II}}, volume~14 of {\em
  Algorithms and Combinatorics}, pages 420--459. Springer, 1997.
\newblock Graham, Ne\v set\v ril, eds.. math.LO/9502233.

\bibitem[Sh 589]{Sh:589}Saharon Shelah.
\newblock {Applications of PCF theory}.
\newblock {\em {Journal of Symbolic Logic}}, {\bf 65}:1624--1674, 2000.

\end{thebibliography}

\enddocument